\documentclass[twocolumn]{autart}
\usepackage{cite}
\usepackage[cmex10]{amsmath}
\usepackage{graphicx}
\usepackage{subfigure}
\usepackage{amsfonts}
\usepackage{mathrsfs}
\usepackage{graphics}
\usepackage{epsfig}
\usepackage{amsfonts}
\usepackage{color}
\usepackage{times}
\usepackage{mathptmx}
\usepackage{amssymb}
\newcommand*{\QEDA}{\hfill\ensuremath{\blacksquare}}
\newtheorem{remark}{Remark}
\newtheorem{thme}{Theorem}

\newtheorem{ass}{Assumption}
\begin{document}
\begin{frontmatter}
\title{Delay-Compensated Control of Sandwiched ODE-PDE-ODE Hyperbolic Systems for Oil Drilling and Disaster Relief\thanksref{footnoteinfo}}
\thanks[footnoteinfo]{The material in this paper was not presented at any conference.}
\author{Ji Wang}\ead{jiw248@eng.ucsd.edu}~and
\author{Miroslav Krstic}\ead{krstic@ucsd.edu}
\address{Department of Mechanical and Aerospace Engineering, University of California, San Diego, La Jolla, CA 92093-0411, USA}

\begin{keyword}
Distributed parameter system; delay; boundary control; oil drilling; unmanned aerial vehicles.
\end{keyword}

\begin{abstract}
Motivated by engineering applications of subsea installation by deepwater construction vessels in oil drilling, and of aid delivery by unmanned aerial vehicles in disaster relief, we develop output-feedback boundary control of heterodirectional coupled hyperbolic PDEs sandwiched between two ODEs, where the measurement is the output state of one ODE and suffers a time delay. After rewriting the time-delay dynamics as a transport PDE of which the left boundary connects with the sandwiched system, a state observer is built to estimate the states of the overall system of ODE-heterodirectional coupled hyperbolic PDEs-ODE-transport PDE using the right boundary state of the last transport PDE. An observer-based output-feedback controller acting at the first ODE is designed  to stabilize the overall system using backstepping transformations and frequency-domain designs. The exponential stability results of the closed-loop system, boundedness and exponential convergence of the control input are proved. The obtained theoretical result is applied to control of a deepwater oil drilling construction vessel as a simulation case, where the simulation results show the proposed control design reduces cable oscillations and places the oil drilling equipment to be installed in the target area on the sea floor. Performance deterioration under extreme and unmodeled disturbances is also illustrated.
\end{abstract}

\end{frontmatter}

\section{Introduction}\label{sec:Intro}
\subsection{Motivation}
The first motivation of this work arises from off-shore oil drilling, where some equipment, such as a subsea manifold, a subsea pump station, a subsea distribution unit along with associated foundations, flowlines and umbilicals should be installed at designated locations \cite{Stensgaard2010,Standing2002} around the drill center on the seafloor. The installation of the equipment is completed by deepwater construction vessels (DCVs) \cite{Stensgaard2010}, because the installation sites are located outside a radius 45 m of the floating drilling platform (Fig.2 in \cite{Stensgaard2010}) and cannot be accessed by the huge floating drilling platform which has limited access and mobility \cite{Stensgaard2010}, and some of the equipment, such as flowlines, umbilicals, should be installed in advance to prepare to hook up the floating drilling platform when it arrives. The DCV is shown in Fig. \ref{fig:DCV}, where the top of the cable is attached to a crane on a vessel at the ocean surface and the bottom attached to
equipment to be installed at the sea floor, referred to as payloads hereafter. The traditional method in underwater installation by DCVs is regulating the vessel dynamics position and
manipulating the crane to obtain the desired heading for the payload \cite{How2010}. It is not suitable for the deeper water construction in offshore oil drilling (more than a thousand meter) because the cable is very long when the payload is near the seabed, which would increase the natural
period of the cable-payload system and introduce large oscillations \cite{J2017Axial,How2010}. The cable oscillations would cause large offset between the payload and the desired heading position of the crane, namely the designated installation location. In addition to large oscillations of the long cable,
another challenge in the subsea installation is the existence of a sensor delay \cite{How2010} which is due to the fact that the
sensor signal is transmitted over a large distance from the seafloor to the
vessel on the ocean surface through a set of acoustics devices (Ch. 10.6.6 in \cite{Sharma2017}).
It would result in information distortion or even make
the control system lose stability. It is vital to design a delay-compensated control force at the onboard crane to reduce the cable oscillations and then place the equipment in the target area on the sea floor.

The second motivation is aid delivery to dangerous and inaccessible areas, such as flood, earthquake, fire, and industrial disaster victims via unmanned aerial vehicles (UAVs) \cite{Guerrero2015,Palunko2012}, where food, first-aid kits, referred to as suspended objects or payloads hereafter, are tied to the bottom of a cable, of which the other end is hanged to an UAV, i.e., a structure of UAV-cable-payload. The swing/oscillation of cable-payload would appear during the transportation motion due to the properties of the cable and external disturbances, such as wind, which may cause damage to the suspended object, the environment and the people around \cite{Guerrero2015}. At the end of the transport motion, when the UAV arrives at the location directly over to the rescue site and is ready to land the aid supplies, the suspended object naturally continues to swing \cite{Palunko2012} which makes precisely placing these aid supplies at the target position difficult. Therefore, rapid suppression of oscillations of the cable and suspended object through a control force provided by rotor wings of the UAV is required. The measurement can be the oscillation acceleration of the suspended object by an accelerometer placed at the bottom end of the cable. Sensor delay would exist in the process of data acquisition, transmission and integration calculation to obtain the payload oscillation displacement which is used in constructing the observer and controller.  In addition to aid delivery in disaster relief, UAV delivery is also used in some commercial cases to reduce labor cost. For example, some companies use UAVs to transport cargos in storehouses or lift and position building elements in architectural construction \cite{Willmann2012}. Some logistics companies have also begun to use UAVs to deliver packages in a small area \cite{Guerrero2015}.

The vibration/oscillation dynamics of cables are distributed parameter systems modeled by wave PDEs \cite{HeBoundary2017,HeUnified2018,J2017Axial}, and the crane/equipment and UAV/supplies at two ends of the cable can be regarded as lumped tip payloads described by ODEs \cite{He2012,HeCooperative2016}. For the sake of order reduction, the wave PDEs with viscous damping terms describing the cable material damping can be converted to a class of heterodirectional coupled hyperbolic PDE systems \cite{Anfinsen2017Disturbance}, \cite{Deutscher2018Output1}-\cite{Deutscher2017Finite-time} via Riemann transformations \cite{J2018Balancing}.  Therefore, the control problems in the aforementioned two applications come down to a theoretical problem about delay-compensated  boundary control of a sandwiched coupled hyperbolic PDE system.
\subsection{Control of PDE sandwiched system}
Boundary control designs of a transport PDE sandwiched by two ODEs \cite{2008Lyapunov}, \cite{krstic2010compensating}, \cite{Anfinsen2018Stabilization}, viscous Burgers PDE \cite{Liu2000Backstepping}
or heat PDE \cite{J2019heat} sandwiched systems  were developed in the previous research.  Control design of the coupled hyperbolic PDE sandwiched system mentioned in the last section is more challenging because of the in-domain instability which comes from in-domain couplings between PDEs. Recently, some results
about state-feedback control
of a coupled hyperbolic sandwiched system
was proposed in \cite{J2017Control,Saba2017,Saba2019}. Based on observer designs, output-feedback control of the coupled hyperbolic PDE sandwiched system was designed in \cite{Deutscher2018Output,Meglio2019Robust}.
However, the aforementioned research has not investigated delay compensation (\cite{Meglio2019Robust,Saba2019} only achieve robustness to a small delay) in boundary control of sandwiched PDE systems. Actually, time-delay exists frequently in the practical engineering, especially the sensor delay, which exists in most practical sensor-used feedback systems. Considering time-delay compensation in the control design is an important step to apply theoretical results into practice.
\subsection{Sensor Delay Compensation}
The topic of sensor delay compensation has received much attention in the past three decades.
In an advanced result presented in \cite{krstic2008Backstepping,krstic2009delay}, the sensor delay was captured as a transport PDE and then the original plant of ODE with sensor delay was rewritten as an ODE-transport PDE cascaded system without delay, before the observer/controller designs were conducted via backstepping. Therein, the observer was built as a ``full-order'' type which estimated
both plant states and sensor states, compared with some classical results about sensor-delay-compensated observer designs \cite{Ahmed-Ali2013Global}, \cite{Cacace2010}, \cite{Germani2002} which only estimated
plant states, namely ``reduced-order'' type.  Using a model-based predictor, observer design for ODE systems
with a time-varying sensor delay was presented in \cite{2010Lyapunov}. Time-varying sensor delay compensation was also considered in \cite{JiState2019} which designed a delay-compensated observer to estimate vibration states of a wave PDE modeled cable elevator. Boundary stabilization of a wave PDE whose boundary observation suffers a time delay was also proposed in \cite{Guo2008}. In the aforementioned work, the sensor delay was considered in the plant which is an ODE or a simple form PDE  while the sensor delay  in this paper exists in a more complex plant which is a sandwiched PDE system.
\subsection{Main Contribution}
\begin{itemize}
\item Some restrictions on the proximal ODE  structure in the previous results about boundary control of ODE-hyperbolic PDE-ODE sandwiched systems are relieved, such as the first-order and scalar form \cite{Anfinsen2018Stabilization,Meglio2019Robust}, a chain of integrators \cite{J2017Control}, $\det(C_0B_0)\neq 0$ \cite{Deutscher2018Output} and $B_0$ being invertible \cite{Saba2017}.
\item Compared with \cite{Saba2019} which first addressed the left invertible type proximal ODE in sandwiched PDEs, this paper further proposes observer-based output-feedback control design using a delayed measurement of which the delay length is constant and arbitrary. Compensation of sensor delay has not been investigated  in control of sandwiched systems before. This is a more challenging task because the plant is extended to ODE-coupled hyperbolic PDEs-ODE-transport PDE after rewriting the delay as a transport PDE.
\item The obtained theoretical result is applied to oscillation suppression of a DCV with compensating the sensor delay arising from large-distance transmission of the sensing signal via acoustics devices, where only one control force  at the onboard crane is required while one more control force applied at the payload is required in \cite{How2011}.
\end{itemize}
For complete clarity, the comparisons with the recent results of boundary control of sandwiched systems are summarized in Tab. \ref{tab:com}.
\begin{table*}
\caption{Comparisons with recent results of boundary control of linear ODE-PDE-ODE systems.}\label{tab:com}
\begin{tabular}{lcccccc}
  \hline
& Types of ODEs & Types  & Types of  & Delay & Application to \\
  & of actuation& of PDEs &control systems  &compensation&   practical problems\\
  \hline
  \cite{Anfinsen2018Stabilization}& First-order and scalar & Transport PDE & Output-feedback & $\times$& $\times$ \\
  \cite{Meglio2019Robust} & First-order and scalar & $2\times 2$ coupled transport PDEs & Output-feedback & $\times$& $\times$ \\
  \cite{J2017Control} & A chain of integrators & $2\times 2$ coupled transport PDEs & Output-feedback & $\times$& $\times$ \\
   \cite{J2019heat} & A chain of integrators & Heat PDE & Output-feedback & $\times$& $\times$ \\
  \cite{Deutscher2018Output} & $\det(C_0B_0)\neq 0$ & $n$ coupled transport PDEs & Output-feedback & $\times$& $\times$ \\
  \cite{Saba2017} & $B_0$ being invertible & $2\times 2$ coupled transport PDEs & State-feedback & $\times$& $\times$ \\
 \cite{Saba2019} & left invertible & $2\times 2$ coupled transport PDEs & State-feedback & $\times$& $\times$ \\
  This paper & left invertible & $2\times 2$ coupled transport PDEs & Output-feedback & $\surd$&DCV \\
  \hline
\end{tabular}

\footnotesize{$\surd$ denotes ``included'' and $\times$ denotes ``not included''.}
\end{table*}
\subsection{Organization} The concerned model and the control task is described in Section \ref{sec:problem}. Observer design is proposed in Section \ref{sec:observer}. Therein, three transformations are used to convert the observer error system to a target observer error system whose exponential stability is straightforward to obtain, where all the dynamics output injections required in constructing the observer are determined. Observer-based output-feedback control design is proposed in Section \ref{sec:output}, where two transformations are applied to transform the observer to a so-called target system in a ``stable-like'' form  except for the proximal ODE which is influenced by perturbations originating from PDEs and the distal ODE. After representing this target system in the frequency domain to obtain the relationships between the states of the proximal ODE and those perturbation states, the proximal ODE is reformulated as a new ODE without external perturbations in the frequency domain, and then the stabilizing control input is designed. The exponential stability of the closed-loop system and the boundedness and exponential convergence to zero of the control input are proved in Section \ref{sec:stability}. The obtained theoretical result is applied to oscillation suppression and position control of a DCV used for seabed installation as a simulation case in Section \ref{sim}. The conclusions and a discussion of future work are provided in Section \ref{sec:conclusion}.
\section{Problem Formulation}\label{sec:problem}
\subsection{Model description}
\begin{figure}
\centering
\includegraphics[width=9cm]{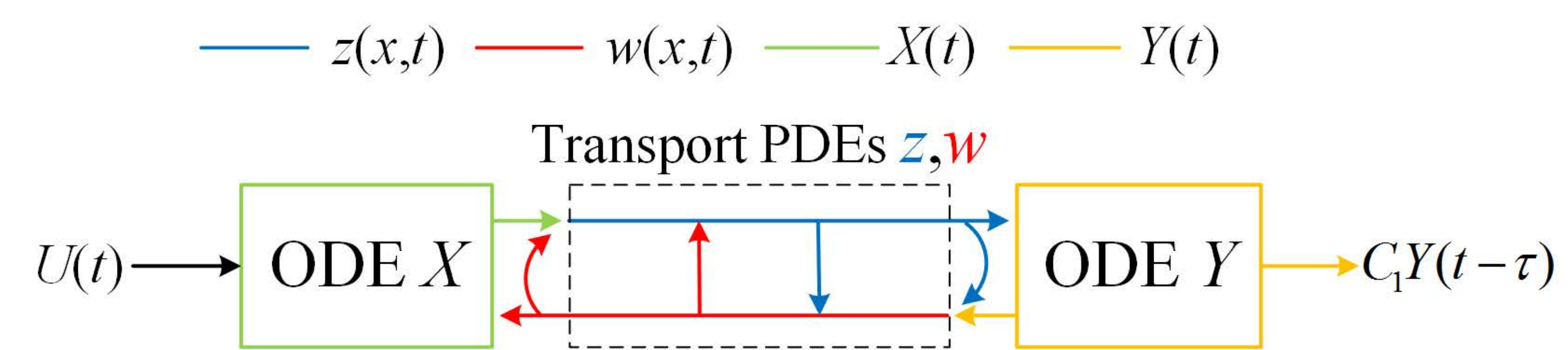}
\caption{Block diagram of the plant \eqref{eq:plant1}-\eqref{eq:plant7}.}
\label{fig:2}
\end{figure}

The plant considered in this paper is
\begin{align}
\dot X(t) &= {A_0}X(t) + {E_0}w(0,t) + {B_0}U(t),\label{eq:plant1}\\
z(0,t) &= pw(0,t) + C_0X(t),\\
{z_t}(x,t) &=  - {q_1}{z_x}(x,t) - {{{c_1}}}w(x,t)- {{{c_1}}}z(x,t),\\
{w_{t}}(x,t) &= {q_2}{w_{x}}(x,t) - {{{c_2}}}w(x,t)- {{{c_2}}}z(x,t),\\
w(1,t) &= qz(1,t) + {C_1}Y(t),\\
\dot Y(t) &= {A_1}Y(t) + {B_1}z(1,t),\\
y_{\rm out}(t)&={C_1}Y(t-\tau)\label{eq:plant7}
\end{align}
$\forall (x,t) \in [0,1]\times[0,\infty)$. The  block diagram of \eqref{eq:plant1}-\eqref{eq:plant7} is shown in Fig. \ref{fig:2}.
$X(t) \in \mathbb{R}^{n\times1}$, $Y(t) \in \mathbb{R}^{m\times1}$ are ODE states. $z(x,t)\in \mathbb{R}, w(x,t)\in \mathbb{R}$ are states of the $2\times2$ coupled hyperbolic PDEs with initial conditions $(z(x,0),w(x,0))\in L^2(0,1)\times L^2(0,1)$. $\tau$ is an arbitrary constant denoting the time delay in the measurement. $U(t)$ is the control input to be designed. $c_1,c_2\in \mathbb{R}$ and $E_0\in \mathbb{R}^{n\times 1}$ are arbitrary.  $q_1$ and $q_2$ are positive transport velocities.  $q,p\in \mathbb{R}$ satisfy Assumption \ref{as:pq}. ${A_0}\in \mathbb{R}^{n\times n}$, ${B_0}\in \mathbb{R}^{n\times 1}$, ${C_0}\in \mathbb{R}^{1\times n}$, ${A_1}\in \mathbb{R}^{m\times m}$, ${B_1}\in \mathbb{R}^{m\times 1}$, ${C_1}\in \mathbb{R}^{1\times m}$
satisfy Assumptions \ref{as:controllable}-\ref{ABC0}.
\begin{ass}\label{as:pq}
$p$, $q$ satisfy
\begin{align}
|pq|<e^{\frac{{{c_2}}}{{{q_2}}} + \frac{{{c_1}}}{{{q_1}}}}
\end{align}
and $q\neq 0$.
\end{ass}
This assumption will be used in the output-feedback control design in Section \ref{sec:output}.
\begin{ass}\label{as:controllable}
The pairs $(A_0,B_0)$, $(A_1,B_1)$ are stabilizable and $(A_0,C_0)$, $(A_1,C_1)$ are detectable.
\end{ass}
According to Assumption \ref{as:controllable}, there exist constant matrices $L_0$, $L_1$, $F_0$, $F_1$ to make the following matrices Hurwitz:
\begin{align}
\bar A_0&=A_0-L_0C_0,\label{eq:barA0}\\
\bar A_1&=A_1-e^{\tau A_1}L_1C_1e^{-\tau A_1}\label{eq:barA1},\\
\hat A_0&=A_0-B_0F_0,\label{eq:hatA0}\\
\hat A_1&=A_1-B_1F_1.\label{eq:bA1}
\end{align}
Note that $A_1-e^{\tau A_1}L_1C_1e^{-\tau A_1}$ has the same eigenvalues as $A_1-L_1C_1$ \cite{krstic2008Backstepping}.
\begin{ass}\label{ABC0}
$(C_0,A_0,B_0)$ satisfy
\begin{align}
\emph{det}\bigg(\left[
            \begin{array}{cc}
              sI-A_0 & B_0 \\
              C_0 & 0 \\
            \end{array}
          \right]
\bigg)\neq 0
\end{align}
for all $s\in \mathbb{C}, {\Re}(s){\ge}0$.
\end{ass}
Assumption \ref{ABC0} is about matrices of the proximal ODE-$X(t)$, namely actuator dynamics. Even though zeros in the closed right-half plane are excluded here while the zeros
are allowed in some previous results on control of sandwiched PDE systems, such as \cite{Deutscher2018Output}, this assumption relieves some restrictions on the
structure of the proximal ODE in the existing literature  (such as $A_0,B_0,C_0$ being scalar in \cite{Anfinsen2018Stabilization,Meglio2019Robust}, $B_0$ being invertible in \cite{Saba2017}, $\det(C_0B_0)\neq 0$ in \cite{Deutscher2018Output}, or a form of a chain of integrators in \cite{J2017Control,J2019heat}) \cite{Saba2019}. This assumption is equal to the existence of a stable left inversion system \cite{Moylan1977} of \eqref{eq:plant1} and is used in the control input design in Section \ref{eq:cdesign}.
\begin{ass}\label{ABC1a}
$(C_1,A_1,B_1)$ satisfy
\begin{align}
{\rm det}\bigg(\left[
            \begin{array}{cc}
              sI-A_1 & B_1 \\
              C_1e^{-\tau A_1} & 0 \\
            \end{array}
          \right]
\bigg)\neq 0
\end{align}
for all $s\in \mathbb{C}, {\Re}(s){\ge}0$.
\end{ass}
Assumption \ref{ABC1a} is about matrices of the distal ODE-$Y(t)$ with a sensor delay $\tau$ in the measurement output state. This assumption also prohibits the zeros of the ODE subsystem $(C_1,A_1,B_1)$ are located in the closed right-half plane. It is not particularly restrictive and $(C_1,A_1,B_1)$ is still quite general covering many application cases.  This assumption is used in the observer design for the overall sandwiched system with the delayed measurement in Section \ref{sec:observer}. Note that if the sensor delay is zero, this assumption has the same form as Assumption \ref{ABC0}.
\begin{remark}
\emph{The design in this paper also can be suitable for collocated control, namely the measurement is the output state of the proximal ODE $X(t)$ with a time delay $\tau$, if $(C_0,A_0,E_0)$ satisfies
\begin{align*}
{\rm det}\bigg(\left[
            \begin{array}{cc}
              sI-A_0 & E_0 \\
              C_0e^{-\tau A_0} & 0 \\
            \end{array}
          \right]
\bigg)\neq 0
\end{align*}}
for all $s\in \mathbb{C}, {\Re}(s){\ge}0$.
\end{remark}
The control objective of this paper: exponentially stabilize the overall sandwiched system, i.e.,the ODE states $Y(t),X(t)$ and the PDE states $u(x,t),v(x,t)$, by constructing an output-feedback control input $U(t)$ applied at the proximal ODE $X(t)$, using the delayed measurement $y_{\rm out}(t)$.
\subsection{Rewrite delay as transport PDE}\label{sec;refeorm}
By defining
\begin{align}
v(x,t) = C_1Y(t-\tau(x-1)),~\forall(x,t)\in[1,2]\times[\tau,\infty)
\end{align}
we obtain a transport PDE
\begin{align}
{v}(1,t) &=  C_1Y(t),\label{eq:delay1}\\
{v_t}(x,t) &=  -\frac{1}{\tau} {v_x}(x,t),\\
y_{\rm out}(t)&=v(2,t)\label{eq:delay3}
\end{align}
$\forall(x,t)\in[1,2]\times[\tau,\infty)$, to describe the time delay in the measurement \eqref{eq:plant7}.

Replacing \eqref{eq:plant7} by \eqref{eq:delay1}-\eqref{eq:delay3}, we obtain a sandwiched hyperbolic PDE-ODE system connecting with another transport PDE, i.e., the following ODE-coupled hyperbolic  PDEs-ODE-transport PDE system:
\begin{align}
\dot X(t) &= {A_0}X(t) + {E_0}w(0,t) + {B_0}U(t),\label{eq:plante1}\\
z(0,t) &= pw(0,t) + C_0X(t),\\
{z_t}(x,t) &=  - {q_1}{z_x}(x,t)\notag\\
 &\quad- {{{c_1}}}w(x,t) - {{{c_1}}}z(x,t),~x\in[0,1]\\
{w_{t}}(x,t) &= {q_2}{w_{x}}(x,t)\notag\\
 &\quad- {{{c_2}}}w(x,t)- {{{c_2}}}z(x,t),~x\in[0,1]\\
w(1,t) &= qz(1,t) + {C_1}Y(t),\\
\dot Y(t) &= {A_1}Y(t) + {B_1}z(1,t),\\
{v}(1,t) &=  C_1Y(t),\\
{v_t}(x,t) &=  -\frac{1}{\tau} {v_x}(x,t),~x\in[1,2]\\
y_{\rm out}(t)&=v(2,t),\label{eq:plante9}
\end{align}
for $t\ge\tau$.

Note that the time delay is ``removed'' at a cost of adding a transport PDE into the plant \eqref{eq:plant1}-\eqref{eq:plant7}. Now, the control task is equivalent to exponentially stabilizing overall system \eqref{eq:plante1}-\eqref{eq:plante9}, i.e., ODE$(X)$-PDE$(z,w)$-ODE$(Y)$-PDE$(v)$, by constructing an output-feedback control input $U(t)$ at the first ODE \eqref{eq:plante1}, using the right boundary state of the last PDE \eqref{eq:plante9}.
\section{Observer design}\label{sec:observer}
In order to build the observer-based output-feedback controller of the plant \eqref{eq:plant1}-\eqref{eq:plant7}, in this section, we  design a state-observer to track the overall system \eqref{eq:plant1}-\eqref{eq:plant7} only using the delayed measurement $y_{\rm out}(t)$. Through the reformulation in Section \ref{sec;refeorm}, the estimation task is equivalent to designing a state-observer to recover the overall system \eqref{eq:plante1}-\eqref{eq:plante9} only
using measurements at the right boundary $x=2$ of the last transport PDE $v$.

The observer is built as a copy of the plant \eqref{eq:plante1}-\eqref{eq:plante9} plus some dynamics output injections:
\begin{align}
&\dot{\hat X}(t) = {A_0}\hat X(t) + {E_0}\hat w(0,t) + {B_0}U(t)\notag\\
&\quad\quad\quad+h_1(y_{\rm out}(t)-\hat v(2,t)),\label{eq:observer1}\\
&\hat z(0,t) = p\hat w(0,t) + C_0\hat X(t),\label{eq:observer2}\\
&{\hat z_t}(x,t) =  - {q_1}{\hat z_x}(x,t) - {{{c_1}}}\hat w(x,t)- {{{c_1}}}\hat z(x,t)\notag\\
&\quad\quad\quad\quad+h_2(y_{\rm out}(t)-\hat v(2,t);x),\label{eq:observer3}\\
&{\hat w_{t}}(x,t) = {q_2}{\hat w_{x}}(x,t)  - {{{c_2}}}\hat w(x,t)- {{{c_2}}}\hat z(x,t)\notag\\
&\quad\quad\quad\quad+h_3(y_{\rm out}(t)-\hat v(2,t);x),\label{eq:observer4}\\
&\hat w(1,t) = q\hat z(1,t) + {C_1}\hat Y(t)+h_4(y_{\rm out}(t)-\hat v(2,t)),\label{eq:observer5}\\
&\dot {\hat Y}(t) = {A_1}\hat Y(t) + {B_1}\hat z(1,t)+\Gamma_1(y_{\rm out}(t)-\hat v(2,t)),\label{eq:observer6}\\
&\hat {v}(1,t) =  C_1\hat Y(t),\label{eq:observer7}\\
&{\hat v_t}(x,t) =  -\frac{1}{\tau} {\hat v_x}(x,t)+h_5(y_{\rm out}(t)-\hat v(2,t);x)\label{eq:observer8}
\end{align}
where a constant matrix $\Gamma_1$ and dynamics $h_1,h_2,h_3,h_4,h_5$  are to be determined. Initial conditions are taken as $(\hat z(x,0),\hat w(x,0),\hat v(x,0))\in L^2(0,1)\times L^2(0,1)\times L^2(1,2)$.
Defining observer error states as
\begin{align}
&(\tilde z(x,t),\tilde w(x,t),\tilde X(t),\tilde Y(t),\tilde v(x,t))\notag\\
=&(z(x,t),w(x,t),X(t),Y(t),v(x,t))\notag\\
&-(\hat z(x,t),\hat w(x,t),\hat X(t),\hat Y(t),\hat v(x,t)),\label{eq:errorstates}
\end{align}
according to \eqref{eq:plante1}-\eqref{eq:plante9} and  \eqref{eq:observer1}-\eqref{eq:observer8}, the observer error system is obtained as
\begin{align}
&\dot{\tilde X}(t) = {A_0}\tilde X(t) + {E_0}\tilde w(0,t)-h_1(\tilde v(2,t)),\label{eq:error1}\\
&\tilde z(0,t) = p\tilde w(0,t) + C_0\tilde X(t),\label{eq:error2}\\
&{\tilde z_t}(x,t) =  - {q_1}{\tilde z_x}(x,t) - {{{c_1}}}\tilde w(x,t)- {{{c_1}}}\tilde z(x,t)\notag\\
&\quad\quad\quad\quad-h_2(\tilde v(2,t);x),\label{eq:error3}\\
&{\tilde w_{t}}(x,t) = {q_2}{\tilde w_{x}}(x,t)- {{{c_2}}}\tilde w(x,t) - {{{c_2}}}\tilde z(x,t)\notag\\
&\quad\quad\quad\quad-h_3(\tilde v(2,t);x),\label{eq:error4}\\
&\tilde w(1,t) = q\tilde z(1,t) + {C_1}\tilde Y(t)-h_4(\tilde v(2,t)),\label{eq:error5}\\
&\dot {\tilde Y}(t) = {A_1}\tilde Y(t) + {B_1}\tilde z(1,t)-\Gamma_1\tilde v(2,t),\label{eq:errorY}\\
&\tilde {v}(1,t) =  C_1\tilde Y(t),\label{eq:errorv1}\\
&{\tilde v_t}(x,t) =  -\frac{1}{\tau} {\tilde v_x}(x,t)-h_5(\tilde v(2,t);x),\label{eq:errorv2}
\end{align}
where $\Gamma_1\tilde v(2,t)$ is an output injection, and $h_1(\tilde v(2,t))$, $h_2(\tilde v(2,t);x)$, $h_3(\tilde v(2,t);x)$, $h_4(\tilde v(2,t))$, $h_5(\tilde v(2,t);x)$ are dynamics output injections which are determined by
\begin{align}
h_1(\tilde v(2,t))&=\mathcal L^{-1} [H_1(s)\tilde v(2,s)],\label{eq:f1}\\
h_2(\tilde v(2,t);x)&=\mathcal L^{-1} [H_2(s;x)\tilde v(2,s)],\label{eq:f2}\\
h_3(\tilde v(2,t);x)&=\mathcal L^{-1} [H_3(s;x)\tilde v(2,s)],\label{eq:f3}\\
h_4(\tilde v(2,t))&=\mathcal L^{-1} [H_4(s)\tilde v(2,s)],\label{eq:f4}\\
h_5(\tilde v(2,t);x)&=\mathcal L^{-1} [H_5(s;x)\tilde v(2,s)],\label{eq:f5}
\end{align}
where $\mathcal L^{-1}$ denotes inverse Laplace transform and transfer functions $H_1(s)$, $H_2(s;x)$, $H_3(s;x)$, $H_4(s)$, $H_5(s;x)$ are to be determined later. Note that $x$ in $H_2,H_3,H_5$ is only a parameter.

Introducing \eqref{eq:f1}-\eqref{eq:f5} is helpful in constructing the dynamics $h_i(\cdot)$ in \eqref{eq:observer1}-\eqref{eq:observer8}, because the algebraic relationships between $\tilde v(2,s)$ and other states can be obtained by using Laplace transform, and the transfer functions in \eqref{eq:f1}-\eqref{eq:f5} can be solved in algebraic equations after rewriting the required conditions of achieving an exponentially stable observer error system in the frequency domain.

Determination of $H_1(s)$, $H_2(s;x)$, $H_3(s;x)$, $H_4(s)$, $H_5(s;x)$ and $\Gamma_1$ in the observer \eqref{eq:observer1}-\eqref{eq:observer8}, will be completed through the following three transformations which convert the observer error system \eqref{eq:error1}-\eqref{eq:errorv2}
to a target observer error system whose exponential stability is straightforward to obtain.
\subsection{First transformation}
Applying the transformation:
\begin{align}
\tilde v(x,t)=\tilde \eta (x,t) + \varphi (x)\tilde Y(t),\label{eq:tran1}
\end{align}
where $\varphi (x)$ is to be determined, we intend to convert \eqref{eq:errorY}-\eqref{eq:errorv2} to a ``stable-like'' form as
\begin{align}
\dot {\tilde Y}(t) &= {{\bar A}_1}\tilde Y(t) + {B_1}\tilde z(1,t) - {\Gamma _1}\tilde \eta (2,t),\label{eq:YHur}\\
\tilde\eta (1,t) &= 0,\label{eq:errortart1}\\
{{\tilde \eta }_t}(x,t) &= -\frac{1}{\tau}{{\tilde \eta }_x}(x,t), ~~x\in[1,2]\label{eq:errortart2}
\end{align}
where $\bar A_1$ is a Hurwitz matrix defined in \eqref{eq:barA1}.
In what follows, $\varphi(x)$, $\Gamma_1$, $H_5(s;x)$ are determined in matching \eqref{eq:errorY}-\eqref{eq:errorv2} and \eqref{eq:YHur}-\eqref{eq:errortart2} via \eqref{eq:tran1}.

Inserting the transformation \eqref{eq:tran1} into \eqref{eq:errorY}, we have
\begin{align}
\dot {\tilde Y}(t) = ({A_1} -\Gamma_1 \varphi (2))\tilde Y(t) + {B_1}\tilde z(1,t)-\Gamma_1\tilde \eta (2,t).
\end{align}
Considering  \eqref{eq:YHur}, \eqref{eq:barA1}, $\Gamma_1$ should satisfy
\begin{align}
\Gamma_1 \varphi (2)=e^{\tau A_1}L_1C_1e^{-\tau A_1}.\label{eq:G1}
\end{align}
Evaluating \eqref{eq:tran1} at $x=1$ and applying \eqref{eq:errorv1}, \eqref{eq:errortart1}, we have
\begin{align}
\tilde v(1,t)&=\tilde \eta (1,t) + \varphi (1)\tilde Y(t)\notag\\
&= \varphi (1)\tilde Y(t)=C_1\tilde Y(t).
\end{align}
Therefore,
\begin{align}
\varphi (1) = C_1.\label{eq:convar1}
\end{align}
Taking the time and spatial derivatives of \eqref{eq:tran1} and submitting the result into \eqref{eq:errorv2}, we have
\begin{align}
&\quad{{\tilde v }_t}(x,t)+\frac{1}{\tau}{{\tilde v }_x}(x,t)+h_5(\tilde v(2,t);x)\notag\\
&={{\tilde \eta}_t}(x,t) + \varphi (x){A_1}\tilde Y(t)+ \varphi (x){B_1}\tilde z(1,t)-\varphi (x)\Gamma_1\tilde v(2,t) \notag\\
 &\quad+ \frac{1}{\tau}{{\tilde \eta}_x}(x,t) +\frac{1}{\tau}\varphi '(x)\tilde Y(t)+h_5(\tilde v(2,t);x)\notag\\
&=\varphi (x){B_1}\tilde z(1,t)-\varphi (x)\Gamma_1\tilde v(2,t)+h_5(\tilde v(2,t);x)\notag\\
 &\quad+ [\varphi (x){A_1}+\frac{1}{\tau}\varphi '(x)]\tilde Y(t)=0,\label{eq:vtvx}
\end{align}
where \eqref{eq:errorY}, \eqref{eq:errortart2} are used, and $\varphi (x)$ should satisfy
\begin{align}
\varphi '(x) &= -\tau{A_1}\varphi (x),\label{eq:convar2}
\end{align}
to make $[\varphi (x){\bar A_1}+\frac{1}{\tau}\varphi '(x)]\tilde Y(t)$ zero, and $H_5(s;x)$ which determines $h_5(\tilde v(2,t);x)$ via \eqref{eq:f5} should be defined to ensure the rest term in \eqref{eq:vtvx} is zero, i.e.,
\begin{align}
\varphi (x){B_1}\tilde z(1,t)-\varphi (x)\Gamma_1\tilde v(2,t)+h_5(\tilde v(2,t);x)=0. \label{eq:rt1}
\end{align}
Before determining $H_5(s;x)$, we solve conditions \eqref{eq:G1}, \eqref{eq:convar1},  \eqref{eq:convar2} to obtain $\varphi (x),\Gamma_1$ as
\begin{align}
 \varphi (x)&=C_1e^{-\tau A_1(x-1)},~~x\in[1,2]\label{eq:var}\\
 \Gamma_1 &=e^{\tau A_1}L_1.\label{eq:Gamma}
\end{align}
Considering \eqref{eq:errortart1}-\eqref{eq:errortart2}, we know
\begin{align}
\tilde \eta(2,t)=0,~t\ge {\tau}.\label{eq:eta2}
\end{align}
Thus \eqref{eq:YHur} can be written as
\begin{align}
\dot {\tilde Y}(t) = {{\bar A}_1}\tilde Y(t) + {B_1}\tilde z(1,t),\label{eq:YHur1}
\end{align}
for $t\ge {\tau}$.
Taking Laplace transform of \eqref{eq:YHur1}, we have
\begin{align}
(sI-\bar A_1){\tilde Y}(s) &=  {B_1}\tilde z(1,s),\label{eq:Ys}
\end{align}
where $I$ is an identity matrix with appropriate dimension. Note: for
brevity, we consider all zero initial conditions while taking Laplace transform (arbitrary initial conditions could be incorporated into the stability statement through an expanded
analysis which is routine).
Recalling $\bar A_1$ being Hurwitz, ${{\rm det}}(sI-\bar A_1)$ does not have any zeros in the closed
right-half plane. Then the matrix $sI-\bar A_1$ is invertible for any $s\in\mathbb{C}$, ${\Re}(s)\ge 0$. Multiplying both sides of \eqref{eq:Ys} by $(sI-{{\bar A}_1})^{-1}$, we have
\begin{align}
{\tilde Y}(s)= ({{sI}}-{{\bar A}_1})^{-1}{B_1}\tilde z(1,s).\label{eq:zY1}%eq:rela
\end{align}
According to \eqref{eq:tran1} and \eqref{eq:eta2}, we have
\begin{align}
\tilde v(2,t)=\varphi (2)\tilde Y(t), ~~t\ge {\tau}.\label{eq:tran10}
\end{align}
Writing \eqref{eq:tran10} in the frequency domain, inserting \eqref{eq:var}, \eqref{eq:zY1} we have
\begin{align}
\tilde v(2,s)&=\varphi (2)\tilde Y(s)=r(s)\tilde z(1,s),\label{eq:tran11}
\end{align}
where
\begin{align}
r(s)=C_1e^{-\tau A_1}(sI-{{\bar A}_1})^{-1}{B_1}.\label{eq:rs}
\end{align}
Notice $r(s)\in\mathbb{R}$  due to $C_1\in \mathbb{R}^{1\times m}$ and $B_1\in\mathbb{R}^{m\times 1}$.
\begin{lem}\label{Lem:r}
$r(s)=C_1e^{-\tau A_1}(sI-{{\bar A}_1})^{-1}{B_1}\in\mathbb{R}$ is nonzero for any $s\in\mathbb{C}$, ${\Re}(s){\ge}0$ under Assumptions \ref{as:controllable} and \ref{ABC1a}.
\end{lem}
\textbf{{Proof.}}
Using \eqref{eq:barA1} in Assumption \ref{as:controllable}, we have
\begin{align}
&\left[
          \begin{array}{cc}
            I & e^{\tau A_1}L_1 \\
            0 & I \\
          \end{array}\right]\left[
                       \begin{array}{cc}
                         sI-A_1 & B_1 \\
                         C_1e^{-\tau A_1} & 0 \\
                       \end{array}
                     \right]\left[
          \begin{array}{cc}
            I & -(sI-\bar A_1)^{-1}B_1 \\
            0 & I \\
          \end{array}\right]\notag\\&=\left[
                       \begin{array}{cc}
                         sI-\bar A_1 & 0 \\
                         C_1e^{-\tau A_1} & -C_1e^{-\tau A_1}(sI-\bar A_1)^{-1}B_1 \\
                       \end{array}
                     \right].
\end{align}
Recalling Assumption \ref{ABC1a}, we know
\begin{align*}
{\rm det}\bigg(\left[
                       \begin{array}{cc}
                         sI-\bar A_1 & 0 \\
                         C_1e^{-\tau A_1} & -C_1e^{-\tau A_1}(sI-\bar A_1)^{-1}B_1 \\
                       \end{array}
                     \right]\bigg)\neq 0
\end{align*}
for any $s\in\mathbb{C}$, ${\Re}(s)\ge 0$. Thereby $C_1e^{-\tau A_1}(sI-{{\bar A}_1})^{-1}{B_1}\neq 0$. The proof of this lemma is completed.\QEDA

According to  Lemma \ref{Lem:r}, we know the existence of $r(s)^{-1}=\frac{1}{r(s)}$.
Let us go back to \eqref{eq:rt1} to determine $H_5(s;x)$ now.
Taking Laplace transform of \eqref{eq:rt1} and recalling \eqref{eq:f5}, inserting \eqref{eq:var} and \eqref{eq:tran11}, we have
\begin{align}
&\quad\varphi (x){B_1}\tilde z(1,s)-\varphi (x)\Gamma_1\tilde v(2,s)+H_5(s;x)\tilde v(2,s)\notag\\
&=\bigg[C_1e^ {-\tau A_1 (x-1)}{B_1}\notag\\
&-\left(C_1e^ {-\tau A_1 (x-1)}\Gamma_1-H_5(s;x)\right)r(s)\bigg]\tilde z(1,s)=0.\label{eq:vtvxs}
\end{align}
$H_5(s;x)$ is chosen as
\begin{align}
H_5(s;x)&=C_1e^ {-\tau A_1 (x-1)}\Gamma_1-C_1e^ {-\tau A_1 (x-1)}{B_1}r(s)^{-1}\notag\\
&=C_1e^ {-\tau A_1 (x-1)}\Gamma_1-\frac{C_1e^ {-\tau A_1 (x-1)}{B_1}}{C_1e^{-\tau A_1}(sI-{{\bar A}_1})^{-1}{B_1}},\label{eq:H5g5a}
\end{align}
where Lemma \ref{Lem:r} is used.

Thereby, \eqref{eq:vtvxs} holds. Then \eqref{eq:rt1} holds by rewriting \eqref{eq:vtvxs} in the time domain. Together with \eqref{eq:convar2}, then \eqref{eq:vtvx} holds for $t\ge {\tau}$. $h_5(\tilde v(2,t);x)$ can then be defined via \eqref{eq:H5g5a} and \eqref{eq:f5}.

In the above, we have completed the conversion between \eqref{eq:errorY}-\eqref{eq:errorv2} and \eqref{eq:YHur}-\eqref{eq:errortart2} through \eqref{eq:tran1} and determined $\Gamma_1,h_5(\tilde v(2,t);x)$ needed in the observer.

In what follows, ${H_4}(s)$ is determined to make the boundary condition \eqref{eq:error5} as zero, i.e.,
\begin{align}
\tilde w(1,t)= q\tilde z(1,t) + {C_1}\tilde Y(t)-h_4(\tilde v(2,t))= 0.\label{eq:w1s}
\end{align}
Taking Laplace transform of \eqref{eq:w1s} and recalling \eqref{eq:f4}, inserting \eqref{eq:zY1} and \eqref{eq:tran11}, we have
\begin{align}
\tilde w(1,s)&=q\tilde z(1,s) + {C_1}\tilde Y(s)-H_4(s)\tilde v(2,s)\notag\\
&=\left(q+C_1{(sI-{{\bar A}_1})^{-1}{B_1}}-{H_4}(s)r(s)\right)\tilde z(1,s)\notag\\
&=0.\label{eq:w2s}
\end{align}
${H_4}(s)$ is chosen as
\begin{align}
{H_4}(s)=&[q+C_1{(sI-{{\bar A}_1})^{-1}{B_1}}]r(s)^{-1}\notag\\
=&\frac{q+C_1{(sI-{{\bar A}_1})^{-1}{B_1}}}{C_1e^{-\tau A_1}(sI-{{\bar A}_1})^{-1}{B_1}}\label{eq:H4s}
\end{align}
to make \eqref{eq:w2s} hold. It follows that $\tilde w(1,t)=0$ in \eqref{eq:w1s} by rewriting $\tilde w(1,s)=0$ in the time domain. $h_4(\tilde v(2,t))$ is thus determined by \eqref{eq:f4}, \eqref{eq:H4s}.

Therefore, through the first transformation \eqref{eq:tran1} with determining the dynamics output injections $h_4(\tilde v(2,t))$, $h_5(\tilde v(2,t);x)$, \eqref{eq:error1}-\eqref{eq:errorv2} can be converted to the first intermediate system as
\begin{align}
\dot{\tilde X}(t) &= {A_0}\tilde X(t) + {E_0}\tilde w(0,t)-h_1(\tilde v(2,t)),\label{eq:iter1}\\
\tilde z(0,t) &= p\tilde w(0,t) + C_0\tilde X(t),\label{eq:iter2}\\
{\tilde z_t}(x,t) &=  - {q_1}{\tilde z_x}(x,t) - {{{c_1}}}\tilde w(x,t)- {{{c_1}}}\tilde z(x,t)\notag\\
&\quad-h_2(\tilde v(2,t);x),\label{eq:iter3}\\
{\tilde w_{t}}(x,t) &= {q_2}{\tilde w_{x}}(x,t)- {{{c_2}}}\tilde w(x,t) - {{{c_2}}}\tilde z(x,t)\notag\\
&\quad-h_3(\tilde v(2,t);x),\label{eq:iter4}\\
\tilde w(1,t) &= 0,\label{eq:iter5}\\
\dot {\tilde Y}(t) &= {\bar A_1}\tilde Y(t) + {B_1}\tilde z(1,t),\label{eq:iter6}\\
\tilde\eta (1,t) &= 0,\label{eq:iter7}\\
{{\tilde \eta }_t}(x,t) &= -\frac{1}{\tau}{{\tilde \eta }_x}(x,t),\label{eq:iter8}
\end{align}
for $t\ge{\tau}$, where \eqref{eq:iter5}-\eqref{eq:iter8} are in a ``stable-like'' form while couplings, i.e., sources terms, exist in the domain $x\in[0,1]$, i.e., \eqref{eq:iter3}-\eqref{eq:iter4}. In the next subsection, we would introduce the second transformation to decouple the couplings in \eqref{eq:iter3}-\eqref{eq:iter4}.
\subsection{Second transformation}
We now apply the second transformation \cite{Bin2017}
\begin{align}
\tilde w(x,t) = \tilde \beta (x,t) - \int_x^1 {\psi (x,y)\tilde \alpha (y,t)dy},\label{eq:trans2a}\\
\tilde z(x,t) = \tilde \alpha (x,t) - \int_x^1 {\phi (x,y)\tilde \alpha (y,t)dy},\label{eq:trans2b}
\end{align}
where the kernels $\psi (x,y),\phi (x,y)$ satisfy
\begin{align}
& \psi (x,x)=\frac{c_2}{{q_1} + {q_2}},\label{eq:t2c1}\\
&p\psi (0,y) - \phi (0,y)=C_0K_1(y),\label{eq:t2c1a}\\
& - {q_1}{{\psi }_y}(x,y) + {q_2}{{\psi }_x}(x,y)-(c_2-{c_1})\psi (x,y) \notag\\
&- {c_2}\phi (x,y)=0,\label{eq:t2c2}\\
& - {q_1}{{\phi }_x}(x,y)- {q_1}{{\phi }_y}(x,y) - {c_1}\psi (x,y) =0\label{eq:t2c3},
\end{align}
in which $K_1(y)$ will be defined in the next subsection where the well-posedness of \eqref{eq:t2c1}-\eqref{eq:t2c3} will be shown, to convert the first intermediate system \eqref{eq:iter1}-\eqref{eq:iter6} to the second intermediate system as
\begin{align}
\dot {\tilde X}(t) =& {{A}_0}\tilde X(t) + {E_0}\tilde \beta (0,t) - {E_0}\int_0^1 {\psi (0,y)\tilde \alpha (y,t)dy} \notag\\
 &+ {h_1}(\tilde v(2,t)),\label{eq:iter21}\\
\tilde \alpha (0,t) =& p\tilde \beta (0,t) + C_0\tilde X(t)- \int_0^1 {C_0K_1(y)\tilde \alpha (y,t)dy}, \label{eq:iter22}\\
{{\tilde \alpha }_t}(x,t) =&  - {q_1}{{\tilde \alpha }_x}(x,t) + \int_x^1 {\bar M(x,y)\tilde \beta (y,t)dy}\notag\\
&- {c_1}\tilde \alpha (x,t)  - {c_1}\tilde \beta (x,t),\label{eq:iter23}\\
{{\tilde \beta }_t}(x,t) =& {q_2}{{\tilde \beta }_x}(x,t) + \int_x^1 {\bar N(x,y)\tilde \beta (y,t)dy}- {c_2}\tilde \beta (x,t),\label{eq:iter24} \\
\tilde \beta (1,t) =& 0,\label{eq:iter25}\\
\dot {\tilde Y}(t) =& {\bar A_1}\tilde Y(t) + {B_1}\tilde \alpha(1,t),\label{eq:iter26}
\end{align}for $t\ge{\tau}$,
with defining
\begin{align}
\bar M(x,y)=  \int_x^y {\phi (x,\delta)\bar M(\delta,y)} d\delta - {c_1}\phi (x,y),\label{eq:barM}\\
\bar N(x,y)= \int_x^y {\psi (x,\delta)\bar M(\delta,y)} d\delta-{c_1}\psi (x,y).\label{eq:barN}
\end{align}
Note that $\tilde\eta(\cdot,t)$ \eqref{eq:iter7}-\eqref{eq:iter8} is removed for brevity because $\tilde\eta(\cdot,t) \equiv 0$, $t\ge \tau$.

In what follows, $H_2(s;x)$, $H_3(s;x)$ are determined  in matching the first intermediate system  \eqref{eq:iter1}-\eqref{eq:iter6} and the second intermediate system \eqref{eq:iter21}-\eqref{eq:iter26} via \eqref{eq:trans2a}-\eqref{eq:trans2b}.

Inserting \eqref{eq:trans2a}-\eqref{eq:trans2b} into \eqref{eq:iter4} along \eqref{eq:iter23}-\eqref{eq:iter24}, and applying \eqref{eq:t2c1}-\eqref{eq:t2c2}, \eqref{eq:barN}, we have
\begin{align}
&{{\tilde w}_t}(x,t) - {q_2}{{\tilde w}_x}(x,t)+ {c_2}\tilde z(x,t)+ {c_2}\tilde w(x,t)\notag\\
& + {h_3(\tilde v(2,t);x)}\notag\\
 &= {q_1}\psi (x,1)\tilde\alpha (1,t) + {h_3(\tilde v(2,t);x)}=0,\label{eq:wtwx}
\end{align}
of which the detailed calculation is shown in \eqref{eq:appB1}. We thus know the following equation should be satisfied
\begin{align}
{q_1}\psi (x,1)\tilde z (1,t) + {h_3(\tilde v(2,t);x)}=0,\label{eq:f30}
\end{align}
where $\tilde \alpha (1,t)=\tilde z(1,t)$ according to \eqref{eq:trans2b} is used.
Writing \eqref{eq:f30} in the frequency domain and applying \eqref{eq:f3}, \eqref{eq:tran11}, we have
\begin{align}
&\quad {q_1}\psi (x,1)\tilde z (1,s)+ H_3(s;x)\tilde v(2,s)\notag\\
&=\left({q_1}\psi (x,1)+ H_3(s;x)r(s)\right)\tilde z(1,s)=0.\label{eq:redw}
\end{align}
$H_3(s;x)$ should be chosen as
\begin{align}
H_3(s;x)&=-{q_1}\psi (x,1)r(s)^{-1}\notag\\
&=\frac{-{q_1}\psi (x,1)}{C_1e^{-\tau A_1}(sI-{{\bar A}_1})^{-1}{B_1}}\label{eq:H3g3}
\end{align}
to make \eqref{eq:redw} hold. It follows that \eqref{eq:wtwx} holds by rewriting \eqref{eq:redw} in the time domain.
$h_3(\tilde v(2,t);x)$ can then be obtained by \eqref{eq:f3}, \eqref{eq:H3g3}.

Inserting \eqref{eq:trans2a}-\eqref{eq:trans2b} into \eqref{eq:iter3} along \eqref{eq:iter23}-\eqref{eq:iter24}, applying \eqref{eq:t2c3}, \eqref{eq:barM}, we have
\begin{align}
&{{\tilde z}_t}(x,t) + {q_1}{{\tilde z}_x}(x,t)  + {c_1}\tilde w(x,t)+ {c_1}\tilde z(x,t)+ h_2(\tilde v(2,t);x)\notag\\
 &={q_1}\phi (x,1)\tilde \alpha (1,t)+ h_2(\tilde v(2,t);x)=0\label{eq:ztzx}
\end{align}
of which the detailed calculation is shown in \eqref{eq:appB2}. Therefore, $h_2(\tilde v(2,t);x)$ should satisfy
\begin{align}
{q_1}\phi (x,1)\tilde z (1,t)+ h_2(\tilde v(2,t);x)=0\label{eq:f20}
\end{align}
where $\tilde \alpha (1,t)=\tilde z(1,t)$ according to \eqref{eq:trans2b} is used. Taking Laplace transform of \eqref{eq:f20}, and recalling \eqref{eq:f2}, \eqref{eq:tran11}, we have
\begin{align}
&\quad{q_1}\phi (x,1)\tilde z (1,s)+ H_2(s;x)\tilde v(2,s)\notag\\
&=\left({q_1}\phi (x,1)+ H_2(s;x)r(s)\right)\tilde z (1,s)=0.\label{eq:redz}
\end{align}
$H_2(s;x)$ is obtained as
\begin{align}
H_2(s;x)=&-{q_1}\phi (x,1)r(s)^{-1}\notag\\
=&\frac{-{q_1}\phi (x,1)}{C_1e^{-\tau A_1}(sI-{{\bar A}_1})^{-1}{B_1}}\label{eq:Hg2}
\end{align}
to ensure \eqref{eq:ztzx} holds. $h_2(\tilde v(2,t);x)$ can thus be defined by \eqref{eq:Hg2}, \eqref{eq:f2}.

Boundary conditions \eqref{eq:iter2}, \eqref{eq:iter5} follow directly from inserting $x=0$, $x=1$ into \eqref{eq:trans2a}-\eqref{eq:trans2b} and applying \eqref{eq:t2c1a}, \eqref{eq:iter22}, \eqref{eq:iter25}. ODEs \eqref{eq:iter1}, \eqref{eq:iter6} are obtained directly from \eqref{eq:iter21}, \eqref{eq:iter26} via \eqref{eq:trans2a}, \eqref{eq:trans2b} respectively.

The second conversion is thus completed and two PDEs \eqref{eq:iter3}-\eqref{eq:iter4} are decoupled now, which can be seen in \eqref{eq:iter23}-\eqref{eq:iter24}.
\subsection{Third transformation }
In order to decouple the ODE \eqref{eq:iter21} with the PDEs and rebuild the ODE in a stable form, we intend to convert the second intermediate system \eqref{eq:iter21}-\eqref{eq:iter26} to the following target observer error system
\begin{align}
&\dot {\tilde Z}(t) = {{\bar A}_0}\tilde Z(t),\label{eq:targ31}\\
&\tilde \alpha (0,t) = C_0\tilde Z(t),\label{eq:targ32}\\
&{{\tilde \alpha }_t}(x,t) =  - {q_1}{{\tilde \alpha }_x}(x,t) - {c_1}\tilde \alpha (x,t),\label{eq:targ33}\\
&\dot {\tilde Y}(t) = {\bar A_1}\tilde Y(t) + {B_1}\tilde \alpha(1,t),\label{eq:targ36}
\end{align}
for $t\ge t_0=\tau+\frac{1}{q_2}$, where $\bar A_0$ is a Hurwitz matrix defined in \eqref{eq:barA0}.
Please note that $\tilde \beta (x,t)\equiv 0$ after $t_0=\tau+\frac{1}{{q_2}}$ recalling \eqref{eq:iter24}-\eqref{eq:iter25}, and then $\tilde \beta (x,t)$ can be removed for brevity. \eqref{eq:iter21}-\eqref{eq:iter26} can thus be rewritten as
\begin{align}
\dot {\tilde X}(t) &= {{A}_0}\tilde X(t) - {E_0}\int_0^1 {\psi (0,y)\tilde \alpha (y,t)dy} \notag\\
&\quad + {h_1}(\tilde v(2,t)),\label{eq:iter21b}\\
\tilde \alpha (0,t) &= C_0\tilde X(t)- \int_0^1 {C_0K_1(y)\tilde \alpha (y,t)dy}, \label{eq:iter22b}\\
{{\tilde \alpha }_t}(x,t) &=  - {q_1}{{\tilde \alpha }_x}(x,t)- {c_1}\tilde \alpha (x,t),\label{eq:iter25b}\\
\dot {\tilde Y}(t) &= {\bar A_1}\tilde Y(t) + {B_1}\tilde \alpha(1,t),\label{eq:iter26b}
\end{align} for $t\ge t_0$. Note that \eqref{eq:iter25b}-\eqref{eq:iter26b} are the same as \eqref{eq:targ33}-\eqref{eq:targ36}. We thus only need to convert \eqref{eq:iter21b}-\eqref{eq:iter22b} to \eqref{eq:targ31}-\eqref{eq:targ32}.

The following transformation
\begin{align}
\tilde Z(t) = \tilde X(t) - \int_0^1 {{K_1}(y)\tilde \alpha (y,t)dy}\label{eq:trans3}
\end{align}
is applied to complete the conversion, where $K_1(y)$ satisfies
\begin{align}
&{L_0} - {q_1}{K_1}(0)=0,\label{eq:K1c1}\\
&(\bar A_0+c_1){K_1}(y) - {q_1}{K_1}^\prime (y) - {E_0}\psi (0,y) \notag\\
&+ {L_0}C_0K_1(y)=0.\label{eq:K1c2}
\end{align}
Note that  \eqref{eq:t2c1}-\eqref{eq:t2c3} and \eqref{eq:K1c1}-\eqref{eq:K1c2} is a $2\times 2$ hyperbolic PDE-ODE system,
which is a scalar case of the well-posed kernel equations (17)-(23) in \cite{Meglio2017Stabilization} (setting
dimensions in \cite{Meglio2017Stabilization} as 1). Therefore, the conditions of the kernels $\psi (x,y),\phi (x,y)$ in \eqref{eq:trans2a}-\eqref{eq:trans2b} and ${K_1}(y)$ in \eqref{eq:trans3}, i.e., \eqref{eq:t2c1}-\eqref{eq:t2c3}, \eqref{eq:K1c1}-\eqref{eq:K1c2} are well-posed.

In what follows, $H_1(s)$ is determined in matching \eqref{eq:iter21b}-\eqref{eq:iter22b} and \eqref{eq:targ31}-\eqref{eq:targ32} via \eqref{eq:trans3}.
Submitting \eqref{eq:trans3} into \eqref{eq:targ31}, applying \eqref{eq:iter21b}-\eqref{eq:iter25b}, \eqref{eq:K1c1}-\eqref{eq:K1c2}, we have
\begin{align}
&\dot{\tilde Z}(t) - \bar A_0\tilde Z(t) \notag\\
=& {A_0}\tilde X(t) - {E_0}\int_0^1 {\psi (0,y)\tilde \alpha (y,t)dy}  + {h_1}\tilde v(2,t) \notag\\
 &+ {q_1}{K_1}(1)\tilde \alpha (1,t)- {q_1}{K_1}(0)\tilde \alpha (0,t)\notag\\
 & - {q_1}\int_0^1 {{K_1}^\prime (y)\tilde\alpha (y,t)dy}+c_1\int_0^1 {K_1(y)\tilde \alpha (y,t)dy}  \notag \\
 &- A_0\tilde X(t)+ {L_0}C_0\tilde X(t) + \bar A_0\int_0^1 {{K_1}(y)\tilde \alpha (y,t)dy}\notag\\
 =& {h_1}\tilde v(2,t) + {q_1}{K_1}(1)\tilde \alpha (1,t)\notag\\
 &+ [{L_0} - {q_1}{K_1}(0)]\tilde \alpha (0,t)\notag \\
 &+ \int_0^1 \bigg[\bar A_0{K_1}(y)+c_1{K_1}(y) - {q_1}{K_1}^\prime (y) - {E_0}\psi (0,y) \notag\\
 &+ {L_0}C_0K_1(y)\bigg]\tilde \alpha (y,t)dy\notag \\
=& {h_1}(\tilde v(2,t))+ {q_1}{K_1}(1)\tilde \alpha (1,t)=0, ~~t\ge t_0.\label{eq:dZ}
\end{align}
$H_1(s)$ which determines $h_1(\tilde v(2,t))$ by \eqref{eq:f1} can be solved from
\begin{align}
{h_1}(\tilde v(2,t))+ {q_1}{K_1}(1)\tilde z (1,t)=0\label{eq:f10}
\end{align}
where $\tilde \alpha (1,t)=\tilde z(1,t)$ according to \eqref{eq:trans2b} is used.
Writing \eqref{eq:f10} in the frequency domain and applying \eqref{eq:f1}, \eqref{eq:tran11} yield
\begin{align}
&\quad H_1(s)\tilde v(2,s)+ {q_1}{K_1}(1)\tilde z (1,s)\notag\\
&=\left(H_1(s)r(s)+ {q_1}{K_1}(1)\right)\tilde z (1,s)=0.\label{eq:redZ}
\end{align}
$H_1(s)$ is solved as
\begin{align}
H_1(s)&=-{q_1}{K_1}(1)r(s)^{-1}\notag\\
&=\frac{-{q_1}{K_1}(1)}{C_1e^{-\tau A_1}(sI-{{\bar A}_1})^{-1}{B_1}}.\label{eq:H1}
\end{align}
It follows that \eqref{eq:dZ} holds by rewriting \eqref{eq:redZ} in the time domain. $h_1(\tilde v(2,t))$ can then be defined via \eqref{eq:f1}, \eqref{eq:H1}.

Inserting \eqref{eq:trans3} into \eqref{eq:iter22b}, it is straightforward to obtain \eqref{eq:targ32}. Therefore, \eqref{eq:targ31}-\eqref{eq:targ32} is converted from  \eqref{eq:iter21b}-\eqref{eq:iter22b}  through \eqref{eq:trans3} for $t\ge t_0$. The third transformation is completed and the ODE \eqref{eq:targ31} is independent and exponentially stable now.

After the above three transformations, we have converted the original observer error system \eqref{eq:error1}-\eqref{eq:errorv2} to the target observer error system \eqref{eq:targ31}-\eqref{eq:targ36} (for $t\in[t_0,\infty)$, $\tilde \eta(x,t)\equiv 0$ \eqref{eq:iter7}-\eqref{eq:iter8} and $\tilde \beta(x,t)\equiv 0$ \eqref{eq:iter24}-\eqref{eq:iter25} are removed for brevity). Because the original observer error system \eqref{eq:error1}-\eqref{eq:errorv2} is bounded in the finite time $t\in[0,t_0)$, we prove the exponential stability of \eqref{eq:error1}-\eqref{eq:errorv2} for $t\in[t_0,\infty)$ in the next subsection.

Note that dynamics $H_i$ determined above are some dynamic extensions, of which the states are dynamic output injections in the observer \eqref{eq:observer1}-\eqref{eq:observer8}, denoted as follows
\begin{align}
&y_1(t)=h_1(\tilde v(2,t)),~~y_2(x,t)=h_2(\tilde v(2,t);x),\label{eq:y1}\\
&y_3(x,t)=h_3(\tilde v(2,t);x),\\
&y_4(t)=h_4(\tilde v(2,t)),~~y_5(x,t)=h_5(\tilde v(2,t);x),\label{eq:y5}
\end{align}
which are proved as exponentially convergent to zero as well in the next subsection.
\begin{remark}
\emph{$H_i(s)\tilde v(2,s)$ usually generates time-derivatives of $\tilde v(2,t)$ in the time domain, which often appears in using the frequency-domain design approach. In practice, one way to avoid taking the time derivatives which may lead to measurement noise amplification is measuring $n$-order time-derivative states $\partial_t^n v(2,t)$ and calculating $\tilde v(2,t)$ by $n$ times integrations of $\partial_t^n \tilde v(2,t)$, which is actually equal to multiplying  $H_i(s)$ by $\frac{1}{s^n}$ to make $H_i(s)$ proper. As shown in an application case of control of DCV in Section \ref{sim}, payload oscillation acceleration is measured and the velocity is calculated by integrating with the known initial conditions. Measuring acceleration is a prevalent method in many mechanical systems, because the acceleration sensor is cheaper and far easier to manufacture and install \cite{Basturk2014State}.}
\end{remark}
\subsection{Stability analysis of the observer errors}\label{sec:stabilityerror}
\emph{Notation: }Supposing $u(x,t)$ is on a spatial domain $x\in[d_1,d_2]$, $\|u(\cdot,t)\|=\sqrt{\int_{d_1}^{d_2}u(x,t)^2dx}$ denotes the $L^2$ norm  and $\|u(\cdot,t)\|_{\infty}=\sup_{x\in [d_1,d_2] }\{|u(x,t)|\}$ denotes the $\infty$-norm.
$|\cdot|$ denotes the Euclidean norm.
\begin{thme}\label{thm:1}
For any initial data $(\tilde z(x,0),\tilde w(x,0)$, $\tilde v(x,0)$, $\tilde X(0), \tilde Y(0))\in L^2(0,1)\times L^2(0,1)\times L^2(1,2)\times \mathbb{R}^{n}\times \mathbb{R}^{m}$, internal exponential stability of the observer error system \eqref{eq:error1}-\eqref{eq:errorv2} holds in the sense of the norm
\begin{align}
&\|\tilde z(\cdot,t)\|_{\infty}+\|\tilde w(\cdot,t)\|_{\infty}+\|\tilde v(\cdot,t)\|_{\infty}+\left|\tilde X(t)\right|+\left|\tilde Y(t)\right|+|y_1(t)|\notag\\
&+|y_4(t)|+\|y_2(\cdot,t)\|_{\infty}+\|y_3(\cdot,t)\|_{\infty}+\|y_5(\cdot,t)\|_{\infty},\label{eq:errornorm}
\end{align}
with the decay rate being adjustable by $L_0,L_1$.
\end{thme}
\textbf{Proof.}
The stability of the original observer error system can be obtained by analyzing the stability of the target observer error system \eqref{eq:targ31}-\eqref{eq:targ36} and using the invertibility of the transformations. \eqref{eq:targ31}-\eqref{eq:targ36} is a cascade of ${\tilde Z}(t)$ into $\tilde \alpha(\cdot,t)$ into $\tilde Y(t)$.
From \eqref{eq:targ31}, ${\tilde Z}(t)$ is exponentially convergent to zero because ${\bar A}_0$ is Hurwitz. With the method
of characteristics as \cite{Deutscher2018Output} it is easy to show that  $\tilde\alpha(x,t)$ in the PDE subsystem \eqref{eq:targ31}-\eqref{eq:targ32} are exponentially convergent to zero. Because ${\bar A}_1$ is Hurwitz, we have ${\tilde Y}(t)$ is exponentially convergent to zero. The decay rate $\lambda_{e}$ of \eqref{eq:targ31}-\eqref{eq:targ36} depends on the decay rate of the ODEs ${\tilde Z}(t),{\tilde Y}(t)$. In other words, the decay rate $\lambda_{e}$ is adjustable by $L_0,L_1$ according to \eqref{eq:barA0}-\eqref{eq:barA1}. Recalling $\tilde \eta(x,t)\equiv0$ and $\tilde\beta(x,t)\equiv0$ after $t_0=\frac{1}{q_2}+{\tau}$, we obtain $\bar\Omega(t)=\|\tilde\alpha(\cdot,t)\|_{\infty}+\|\tilde\beta(\cdot,t)\|_{\infty}+\|\tilde\eta(\cdot,t)\|_{\infty}+|\tilde Z(t)|+|\tilde Y(t)|$ is bounded by an exponential decay with the decay rate $\lambda_{e}$ for $t\ge t_0$. Note that the transient in the finite time $[0,t_0)$ can be bounded by an arbitrarily fast decay rate considering a trade off between the decay rate and the overshoot coefficient, i.e.,  the higher the decay rate, the higher the overshoot coefficient. Therefore, we conclude the exponential stability in the sense of $\bar\Omega(t)$ being bounded by an exponential decay rate $\lambda_{e}$ with some overshoot coefficients for $t\ge 0$.  Applying the transformation \eqref{eq:tran1}, \eqref{eq:trans3} and \eqref{eq:trans2a}-\eqref{eq:trans2b}, we respectively have
\begin{align*}
\left\|\tilde v(\cdot,t)\right\|_{\infty}&\le \Upsilon_{1a}\left(\|\tilde\eta(\cdot,t)\|_{\infty}+\left|\tilde Y(t)\right|\right),\notag\\
\left|\tilde X(t)\right|&\le \Upsilon_{1b}\left(\|\tilde\alpha(\cdot,t)\|_{\infty}+\left|\tilde Z(t)\right|\right),\notag\\
\|\tilde z(\cdot,t)\|_{\infty}+\|\tilde w(\cdot,t)\|_{\infty}&\le \Upsilon_{1c}\left(\|\tilde\alpha(\cdot,t)\|_{\infty}+\|\tilde \beta(\cdot,t)\|_{\infty}\right),
\end{align*}
for some positive $\Upsilon_{1a},\Upsilon_{1b},\Upsilon_{1c}$.

According to \eqref{eq:f1}-\eqref{eq:f5}, \eqref{eq:H1}, \eqref{eq:Hg2}, \eqref{eq:H3g3}, \eqref{eq:H4s}, \eqref{eq:H5g5a}, we know the output injection states $y_1(t)$, $y_2(x,t)$, $y_3(x,t)$, $y_4(t)$, $y_5(x,t)$ are the output states of the following extended dynamics:
\begin{align}
&H_{1}(s)=\frac{-{q_1}{K_1}(1)}{C_1e^{-\tau A_1}(sI-{{\bar A}_1})^{-1}{B_1}},\label{eq:Tran1}\\
&H_{2}(s;x)=\frac{-{q_1}\phi (x,1)}{C_1e^{-\tau A_1}(sI-{{\bar A}_1})^{-1}{B_1}},\label{eq:Tran2}\\
&H_{3}(s;x)=\frac{-{q_1}\psi (x,1)}{C_1e^{-\tau A_1}(sI-{{\bar A}_1})^{-1}{B_1}},\label{eq:Tran3}\\
&H_{4}(s)=\frac{q+C_1{(sI-{{\bar A}_1})^{-1}{B_1}}}{C_1e^{-\tau A_1}(sI-{{\bar A}_1})^{-1}{B_1}},\label{eq:Tran4}\\
&H_{5}(s;x)=C_1e^ {-\tau A_1 (x-1)}\Gamma_1-\frac{C_1e^ {-\tau A_1 (x-1)}{B_1}}{C_1e^{-\tau A_1}(sI-{{\bar A}_1})^{-1}{B_1}},\label{eq:Tran5}
\end{align}
of which the input signal is $\tilde v(2,t)$ which is exponentially convergent to zero. Recalling Lemma \ref{Lem:r}, we know there is not pole in the closed right-half plane in the transfer function \eqref{eq:Tran1}-\eqref{eq:Tran5}, the exponential convergence of $|y_1(t)|,\|y_2(x,t)\|_{\infty},\|y_3(x,t)\|_{\infty},|y_4(t)|,\|y_5(x,t)\|_{\infty}$ are thus obtained. Note that $x\in[0,1]$ is just a parameter in the transfer functions \eqref{eq:Tran2},\eqref{eq:Tran3},\eqref{eq:Tran5} and the stability result would not be affected. \QEDA
\section{Output-feedback control design}\label{sec:output}
In the last section, we have built the observer which can compensate the time-delay in the output measurement $y_{\rm out}(t)$ of the distal ODE, which is the only one measurement used in the observer, to
track the states of the overall sandwiched PDE system \eqref{eq:plant1}-\eqref{eq:plant7}. In this
section, we design an output-feedback control law $U(t)$ based on the observer \eqref{eq:observer1}-\eqref{eq:observer8} by using backstepping
transformations and frequency-domain designs.

First, two transformations are introduced to transform the observer \eqref{eq:observer1}-\eqref{eq:observer8} to a
target system \eqref{eq:ta1}-\eqref{eq:targ8a}, which is in a ¡°stable-like¡± form except for the proximal
ODE which is influenced by perturbations originating from the
PDEs and distal ODE. Representing this ``target system''
in the frequency domain by using Laplace transform, the algebraic relationships \eqref{eq:r1}-\eqref{eq:rf} between
the states of the proximal ODE and the states of the PDEs and
distal ODE are obtained. Inserting these algebraic relationships to rewrite the perturbations in the proximal ODE, a new ODE \eqref{eq:ZODEs3} without external perturbations can be built in the frequency domain, where the control input to exponentially stabilize this ODE can be designed.
\subsection{First transformation}
The aim of the first transformation is to remove the source terms in the PDE domain $x\in[0,1]$, i.e., couplings in \eqref{eq:observer3}-\eqref{eq:observer4}, and to build the state matrix of the distal ODE \eqref{eq:observer6} as a Hurwitz matrix. A PDE backstepping transformation in the following form \cite{Meglio2017Stabilization}
\begin{align}
\alpha (x,t) =& \hat z(x,t) - \int_x^1 {K_3}(x,y)\hat z(y,t)dy\notag\\
& -  \int_x^1 {J_3}(x,y)\hat w(y,t)dy -  \gamma (x)\hat Y(t),\label{eq:contran1a}\\
\beta (x,t) =& \hat w(x,t) - \int_x^1 {K_2}(x,y)\hat z(y,t)dy \notag\\
&-  \int_x^1 {J_2}(x,y)\hat w(y,t)dy -  \lambda (x)\hat Y(t)\label{eq:contran1b}
\end{align}
is introduced, where the kernels ${K_3}(x,y)$, ${J_3}(x,y)$, $\gamma (x)$, ${K_2}(x,y)$, ${J_2}(x,y)$, $\lambda (x)$ are to be determined later,
to convert \eqref{eq:observer1}-\eqref{eq:observer8} to the following intermediate system:
\begin{align}
\dot {\hat X}(t) =& {A_0}\hat X(t) + {E_0}\beta (0,t) + \int_0^1 {{\bar K}_4}(x)\alpha (x,t)dx \notag\\
&+ \int_0^1 {{{\bar K}_5}(x)\beta (x,t)d} x + {{\bar K}_6}\hat Y(t) \notag\\
&+ {B_0}U(t) + {h_1}(\tilde v(2,t)),\label{eq:targ6}\\
\alpha (0,t) =& p\beta (0,t)+ {C_0}\hat X(t) + \int_0^1 {{\bar K}_1}(x)\alpha (x,t)dx \notag\\
 &+  {\bar K}_3\hat Y(t)+ \int_0^1 {{\bar K}_2}(x)\beta (x,t)dx ,\label{eq:targ5}\\
{\alpha _t}(x,t) =&  - {q_1}{\alpha _x}(x,t) - {c_1}\alpha (x,t)- \gamma (x){\Gamma _1}\tilde v(2,t) \notag\\
&- \int_x^1 {{J_2}(x,y){h_3}(\tilde v(2,t);y)dy}\notag\\
&- \int_x^1 {{K_3}(x,y){h_2}(\tilde v(2,t);y)dy} \notag\\
  &+ h_2(\tilde v(2,t);x) - {q_2}{J_3}(x,1){h_4}(\tilde v(2,t)),\label{eq:targ3}\\
{\beta _t}(x,t) =& {q_2}{\beta _x}(x,t) - {c_2}\beta (x,t)\notag\\
&- \lambda (x){\Gamma _1}\tilde v(2,t) - \int_x^1 {{J_2}(x,y){h_3}(\tilde v(2,t);y)dy}\notag\\
  &- \int_x^1 {{K_2}(x,y){h_2}(\tilde v(2,t);y)dy} \notag\\
  & + {h_3}(\tilde v(2,t);x) - {q_2}{J_2}(x,1){h_4}(\tilde v(2,t)),\label{eq:targ4}\\
  \beta (1,t) =& q\alpha (1,t) + {h_4}(\tilde v(2,t)),\label{eq:targ2}\\
\dot {\hat Y}(t) =& {\hat A_1}\hat Y(t) + {B_1}\alpha (1,t)+ {\Gamma _1}\tilde v(2,t),\label{eq:targ1}\\
\hat v(1,t) =& {C_1}\hat Y(t),\label{eq:targ7}\\
{{\hat v}_t}(x,t) =&  - \frac{1}{\tau}{{\hat v}_x}(x,t) + {h_5}(\tilde v(2,t);x),\label{eq:targ8}
\end{align}
where $\hat A_1$ is a Hurwitz matrix by choosing the control parameter $F_1$ according to Assumption \ref{as:controllable}. ${{\bar K}_1}(x)$, ${{\bar K}_2}(x)$, ${{\bar K}_3}$, ${{\bar K}_4}(x)$, ${{\bar K}_5}(x)$, ${{\bar K}_6}$ satisfy
\begin{align}
{{\bar K}_1}(x)=&p{K_2}(0,x) - {K_3}(0,x) + \int_0^x {{\bar K}_1}(y){K_3}(y,x)dy\notag\\
& + \int_0^x {{{\bar K}_2}(y){K_2}(y,x)dy}  ,\label{eq:barK1}\\
{{\bar K}_2}(x)=&-p{J_2}(0,x) + {J_3}(0,x) + \int_0^x {{\bar K}_1}(y){J_3}(y,x)dy\notag\\
& + \int_0^x {{{\bar K}_2}(y){J_2}(y,x)dy }  ,\label{eq:barK2}\\
{{\bar K}_3}=&\int_0^1 {{{\bar K}_2}(x)\lambda (x)dx}  + \int_0^1 {{\bar K}_1}(x)\gamma (x)dx \notag\\
&+ p\lambda (0) - \gamma (0),\label{eq:barK3}\\
{{\bar K}_4}(x)=&\int_0^x {{\bar K}_4}(y){K_3}(y,x)dy + \int_0^x {{\bar K}_5}(y){K_2}(y,x)dy\notag\\
& - {E_0}{K_2}(0,x) ,  \\
{{\bar K}_5}(x)=&\int_0^x {{\bar K}_4}(y){J_3}(y,x)dy + \int_0^x {{\bar K}_5}(y){J_2}(y,x)dy\notag\\
& + {E_0}{J_2}(0,x),    \\
{{\bar K}_6}=&\int_0^1 {{{\bar K}_5}(x)\lambda (x)dx}  + \int_0^1 {{\bar K}_4}(x)\gamma (x)dx \notag\\
&+ {E_0}\lambda (0),\label{eq:barK6}
\end{align}
which are obtained by matching \eqref{eq:targ6}-\eqref{eq:targ5} and \eqref{eq:observer1}-\eqref{eq:observer2} via \eqref{eq:contran1a}-\eqref{eq:contran1b} (Please see Step 4 in the Appendix-A for the details). The following conditions of the kernels in the transformations \eqref{eq:contran1a}-\eqref{eq:contran1b} are obtained by matching \eqref{eq:targ3}-\eqref{eq:targ1} and \eqref{eq:observer3}-\eqref{eq:observer6} (the detailed process is shown in Steps. 1-3 in the Appendix-A):
\begin{align}
&{q_1}{K_3}(x,1)=  {q_2}{J_3}(x,1)q + \gamma (x){B_1},\label{eq:Kerc1}\\
 &{J_3}(x,x)=\frac {c_1}{{q_2} + {q_1}},\label{eq:q2q1}\\
 & - {q_1}{J_{3x}}(x,y)+ {q_2}{J_{3y}}(x,y)+ ({c_2} - {c_1}){J_3}(x,y) \notag\\
& +{c_1}{K_3}(x,y)=0,\label{eq:c1K3}\\
& - {q_1}{K_{3x}}(x,y) - {q_1}{K_{3y}}(x,y)+{c_2}{J_3}(x,y)=0,\label{eq:c2L3}\\
&\gamma (1) = -F_1,\label{eq:gammaF1}\\
& - {q_1}\gamma'(x) - \gamma (x)({A_1} + {c_1}) - {q_2}{J_3}(x,1){C_1}=0,\label{eq:Kerc1a}\\
& {q_2}q{J_2}(x,1)= {q_1}{K_2}(x,1)- \lambda (x){B_1},\label{eq:Kerc2}\\
&{K_2}(x,x)=\frac {-c_2}{{q_1} + {q_2}},\label{eq:c2q1}\\
&{q_2}{J_{2x}}(x,y) + {q_2}{J_{2y}}(x,y)+{c_1}{K_2}(x,y) =0,\label{eq:c1K2}\\
 & {q_2}{K_{2x}}(x,y) - {q_1}{K_{2y}}(x,y)+ ({c_1}- {c_2}){K_2}(x,y)\notag\\
 &+{c_2}{J_2}(x,y) =0,\label{eq:c2L2}\\
& {q_2}\lambda '(x) - \lambda (x)({A_1} + {c_2}) - {q_2}{J_2}(x,1){C_1}=0,\label{eq:q2lam}\\
& \lambda (1)=q\gamma (1)  + {C_1} .\label{eq:Kerc2a}
\end{align}
The well-posedness of \eqref{eq:Kerc1}-\eqref{eq:Kerc2a} is shown in the following lemma.
\begin{lem}\label{lm:wp}
The kernel equations \eqref{eq:Kerc1}-\eqref{eq:Kerc1a} have a unique solution ${K_3},{J_3}\in C^1(D)$, $\gamma (x)\in C^1([0,1])$ and \eqref{eq:Kerc2}-\eqref{eq:Kerc2a} have a unique solution ${K_2},{J_2}\in C^1(D)$, $\lambda (x)\in C^1([0,1])$, with $D=\{(x,y)|0\le x\le y\le 1\}$.
\end{lem}
\textbf{{Proof.}}
\eqref{eq:Kerc1}-\eqref{eq:Kerc1a} and \eqref{eq:Kerc2}-\eqref{eq:Kerc2a} have the analogous structure with (19)-(24) in \cite{J2017Control}. Following the proof of Lemma 1 in \cite{J2017Control}, we can obtain this lemma.\QEDA

Similarly, the inverse transformation can be obtained as
\begin{align}
\hat z (x,t) =& \alpha(x,t) - \int_x^1 {\mathcal M}(x,y)\alpha(y,t)dy\notag\\
& -  \int_x^1 {\mathcal N}(x,y)\beta(y,t)dy -  \mathcal G (x)\hat Y(t),\label{eq:contran1aI}\\
\hat w (x,t) =& \beta(x,t) - \int_x^1 {\mathcal D}(x,y)\alpha(y,t)dy \notag\\
&-  \int_x^1 {\mathcal T}(x,y)\beta(y,t)dy -  \mathcal P(x)\hat Y(t),\label{eq:contran1bI}
\end{align}
where ${\mathcal M}(x,y), {\mathcal N}(x,y), \mathcal G (x), {\mathcal D}(x,y), {\mathcal T}(x,y), \mathcal P(x)$ are kernels which can be determined through a similar process in the Appendix-A. The first transformation in the control design is completed.
\subsection{Second transformation}
In order to remove the last three terms in the boundary condition \eqref{eq:targ5} and form a Hurwitz matrix of the proximal ODE \eqref{eq:targ6}, we introduce the second transformation
\begin{align}
&\hat Z(t)=\hat X(t) + {C_0}^ + \int_0^1 {{\bar K}_1}(x)\alpha(x,t)dx \notag\\
&+  {C_0}^+ \int_0^1 {{\bar K}_2}(x)\beta(x,t)dx + {C_0}^ +  {{\bar K}_3}\hat Y(t),\label{eq:contran2}
\end{align}
where ${C_0}^+$ denotes the Moore-Penrose right inverse
of ${C_0}$. Note that because $C_0$ is full-row rank (with rank
equal to 1), a right inverse exists for $C_0$. i.e., ${C_0}{C_0}^+=I$. A choice of ${C_0}^+$  is  ${C_0}^+=C_0^T(C_0C_0^T)^{-1}$.

Using \eqref{eq:contran2}, then \eqref{eq:targ6}-\eqref{eq:targ5} is converted to
\begin{align}
&\dot {\hat Z}(t) = {{\hat A}_0}\hat Z(t) + {q_1}{C_0}^ + {{\bar K}_1}(0){C_0}\hat Z(t)+ {B_0}\bar U(t) \notag\\
&+ {{ M}_Y}\hat Y(t) + \int_0^1 {{M_\alpha (x) }} \alpha (x,t)dx+ \int_0^1 {{M_\beta(x) }} \beta (x,t)dx\notag\\
&  + {N_1}\alpha (1,t) + {N_2}\beta (0,t) +\mathcal H[h_1(\tilde v(2,t)),h_2(\tilde v(2,t);x),\notag\\
&h_3(\tilde v(2,t);x),h_4(\tilde v(2,t)),h_5(\tilde v(2,t);x),\tilde v(2,t)],\label{eq:ZODE}\\
&\alpha (0,t) = p\beta (0,t) + {C_0}\hat Z(t),\label{eq:alpha0}
\end{align}
where
\begin{align}
&\mathcal H[h_1(\tilde v(2,t)),h_2(\tilde v(2,t);x),h_3(\tilde v(2,t);x),h_4(\tilde v(2,t)),\notag\\
& h_5(\tilde v(2,t);x),\tilde v(2,t)]\notag\\
= &{h_1}(\tilde v(2,t))+ {C_0}^ + {{\bar K}_3}{\Gamma _1}\tilde v(2,t)+ {q_2}{C_0}^ + {{\bar K}_2}(1)q{h_4}(\tilde v(2,t))\notag\\
  &+ {C_0}^ + \int_0^1 {{{\bar K}_1}(x)\int_x^1 {{J_2}(x,y){h_3}(\tilde v(2,t);y)} dydx}\notag \\
 &- {C_0}^ + \int_0^1 {{{\bar K}_1}(x)\int_x^1 {{K_3}(x,y){h_2}(\tilde v(2,t);y)} dydx} \notag\\
  &+ {C_0}^ + \int_0^1 {{{\bar K}_1}(x)h_2(\tilde v(2,t);x)dx}\notag\\
  &  - {C_0}^ + \int_0^1 {{{\bar K}_1}(x){q_2}{J_3}(x,1)dx{h_4(\tilde v(2,t))}}\notag \\
 &- {C_0}^ + \int_0^1 {{{\bar K}_2}(x)\lambda (x)dx{\Gamma _1}\tilde v(2,t)} \notag\\
  &- {C_0}^ + \int_0^1 {{{\bar K}_2}(x)\int_x^1 {{J_2}(x,y){h_3}(\tilde v(2,t);y)} dydx} \notag\\
 &- {C_0}^ + \int_0^1 {{{\bar K}_2}(x)\int_x^1 {{K_2}(x,y){h_2}(\tilde v(2,t);y)} dydx} \notag\\
  &+ {C_0}^ + \int_0^1 {{{\bar K}_2}(x){h_3}(\tilde v(2,t);x)dx} \notag\\
  & - {C_0}^ + \int_0^1 {{{\bar K}_2}(x){q_2}{J_2}(x,1)dx{h_4(\tilde v(2,t))}},\label{eq:hm}
\end{align}
and
\begin{align}
\bar U(t)=U(t)-{F_0}\hat Z(t).\label{eq:bU}
\end{align}
$\hat A_0$ is Hurwitz by choosing the control parameter $F_0$ considering Assumption \ref{as:controllable}.
${N_1}$,${N_2}$,${M_\alpha }(x)$,${M_\beta }(x)$,${M_Y}$ in \eqref{eq:ZODE} are
\begin{align}
&{N_1} = {C_0}^ + {{\bar K}_3}{B_1} - {q_1}{C_0}^ + {{\bar K}_1}(1) + {q_2}{C_0}^ + {{\bar K}_2}(1)q,\label{eq:N1}\\
&{N_2} = {E_0} - {q_2}{C_0}^ + {{\bar K}_2}(0) + {q_1}{C_0}^ + {{\bar K}_1}(0)p,\\
&{M_\alpha }(x) = {{\bar K}_4}(x)  + {q_1}{C_0}^ + {{\bar K}_1}^\prime (x)- ({{\hat A}_0}+{c_1}){C_0}^ + {{\bar K}_1}(x),\\
&{M_\beta }(x) = {{\bar K}_5}(x) - {q_2}{C_0}^ + {{\bar K}_2}^\prime (x) - ({{\hat A}_0}+{c_2}){C_0}^ + {{\bar K}_2}(x),\\
&{M_Y} = {C_0}^ + {{\bar K}_3}{{\hat A}_1} + {{\bar K}_6} - {{\hat A}_0}{C_0}^ + {{\bar K}_3}.\label{eq:MY}
\end{align}
We thus arrive at the target system consisting of \eqref{eq:targ3}-\eqref{eq:targ8}, \eqref{eq:ZODE}-\eqref{eq:alpha0}, which includes dynamic output injections in \eqref{eq:hm}.
Considering Theorem \ref{thm:1} and \eqref{eq:y1}-\eqref{eq:y5}, we know $h_1(\tilde v(2,t))$, $h_2(\tilde v(2,t);x)$, $h_3(\tilde v(2,t);x)$, $h_4(\tilde v(2,t))$, $h_5(\tilde v(2,t);x)$ and $\Gamma_1\tilde v(2,t)$ in the target system \eqref{eq:targ3}-\eqref{eq:targ8}, \eqref{eq:ZODE}-\eqref{eq:alpha0} can be regarded as zero, i.e., $\mathcal H=0$, for brevity. Therefore, the target system \eqref{eq:targ3}-\eqref{eq:targ8}, \eqref{eq:ZODE}-\eqref{eq:alpha0} can be rewritten as
\begin{align}
&\dot {\hat Z}(t) = {{\hat A}_0}\hat Z(t) + {q_1}{C_0}^ + {{\bar K}_1}(0){C_0}\hat Z(t) \notag\\
&+ {{ M}_Y}\hat Y(t) + \int_0^1 {{M_\alpha(x) }} \alpha (x,t)dx+ \int_0^1 {{M_\beta(x) }} \beta (x,t)dx\notag\\
&  + {N_1}\alpha (1,t) + {N_2}\beta (0,t) + {B_0}\bar U(t),\label{eq:ta1}\\
&\alpha (0,t) = p\beta (0,t) + {C_0}\hat Z(t),\\
&{\alpha _t}(x,t) =  - {q_1}{\alpha _x}(x,t) - {c_1}\alpha (x,t), ~~x\in[0,1]\\
&{\beta _t}(x,t) = {q_2}{\beta _x}(x,t) - {c_2}\beta (x,t),~~x\in[0,1]\\
  &\beta (1,t) = q\alpha (1,t),\\
&\dot {\hat Y}(t) = {\hat A_1}\hat Y(t) + {B_1}\alpha (1,t),\label{eq:ta6}\\
&\hat v(1,t) = {C_1}\hat Y(t),\label{eq:targ7a}\\
&{{\hat v}_t}(x,t) =  - \frac{1}{\tau}{{\hat v}_x}(x,t) ,~~x\in[1,2].\label{eq:targ8a}
\end{align}
\subsection{Control design in frequency domain}\label{eq:cdesign}
In the last two subsections, the system \eqref{eq:observer1}-\eqref{eq:observer8} is converted to the target system \eqref{eq:ta1}-\eqref{eq:targ8a}, through the two transformations \eqref{eq:contran1a}-\eqref{eq:contran1b} and \eqref{eq:contran2}. In this section, the control $\bar U(t)$ in \eqref{eq:ta1} of the target system \eqref{eq:ta1}-\eqref{eq:targ8a}
will be designed in the frequency domain by using Laplace transform.

Taking Laplace transform of \eqref{eq:ta1}-\eqref{eq:targ8a}, we have
\begin{align}
&(sI - {{\hat A}_0})\hat Z(s) = {q_1}{C_0}^ + {{\bar K}_1}(0){C_0}\hat Z(s) + {{M}_Y}\hat Y(s)\notag\\
&+ \int_0^1 {{M_\alpha }(x)} \alpha (x,s)dx + \int_0^1 {{M_\beta }} (x)\beta (x,s)dx\notag\\
&+ {N_1}\alpha (1,s) + {N_2}\beta (0,s) + {B_0}\bar U(s),\label{eq:ZODEs}\\
&\alpha (0,s) = p\beta (0,s) + {C_0}\hat Z(s),\label{eq:alpha0s}\\
&s{\alpha}(x,s) =  - {q_1}{\alpha _x}(x,s) - {c_1}\alpha (x,s),\label{eq:targ3s}\\
&s{\beta}(x,s) = {q_2}{\beta _x}(x,s) - {c_2}\beta (x,s),\label{eq:targ4s}\\
&\beta (1,s) = q\alpha (1,s),\label{eq:targ2s}\\
&(sI-{\hat A_1}){\hat Y}(s) = {B_1}\alpha (1,s),\label{eq:targ1s}\\
&\hat v(1,s) = {C_1}\hat Y(s),\label{eq:targ7as}\\
&s{{\hat v}}(x,s) =  - \frac{1}{\tau}{{\hat v}_x}(x,s).\label{eq:targ8as}
\end{align}
Note: for brevity, we consider all zero initial conditions
while taking Laplace transform (arbitrary initial condition-
s could be incorporated into the stability statement through
an expanded analysis which is routine).

Defining
\begin{align}
h(s)=1 - pq{e^{-(\frac{{{c_2}}}{{{q_2}}} + \frac{{{c_1}}}{{{q_1}}})}}{e^{ - (\frac{1}{{{q_2}}} + \frac{1}{{{q_1}}})s}},
\end{align}
according to \eqref{eq:alpha0s}-\eqref{eq:targ8as} and Sec. 3.2 in \cite{Meglio2019Robust}, we obtain the following algebraic relationships between $C_0\hat Z(s)$ and other states in \eqref{eq:alpha0s}-\eqref{eq:targ8as} as
\begin{align}
&h(s)\alpha(x,s)=e^{\frac{-(c_1+s)}{q_1}x}C_0\hat Z(s),\label{eq:afZ}\\
&h(s)\beta(x,s)=q{e^{\frac{{-(c_2 + s)}}{{{q_2}}}(1 - x) - \frac{{({c_1} + s)}}{{{q_1}}}}}C_0\hat Z(s),\label{eq:BeZ}\\
&h(s)\hat v(x,s)=C_1{(sI - {{\hat A}_1})^{ - 1}}B_1{e^{\frac{{-(c_1 + s)}}{{{q_1}}}-\tau (x-1)s}}{C_0}\hat Z(s),\label{eq:Bevx}\\
&h(s)\hat v(1,s)=C_1{(sI - {{\hat A}_1})^{ - 1}}B_1{e^{\frac{{-(c_1 + s)}}{{{q_1}}}}}{C_0}\hat Z(s),\\
&h(s)\alpha (0,s) = {C_0}\hat Z(s),\label{eq:r1}\\
&h(s)\beta (1,s) = q{e^{\frac{{-(c_1 + s)}}{{{q_1}}}}}{C_0}\hat Z(s),\\
&h(s)\beta (0,s) = q{e^{\frac{{-(c_2 + s)}}{{{q_2}}} - \frac{{({c_1} + s)}}{{{q_1}}}}}{C_0}\hat Z(s),\\
&h(s)\alpha (1,s) = {e^{\frac{{-(c_1 + s)}}{{{q_1}}}}}{C_0}\hat Z(s),\\
&h(s)\hat Y(s)= {(sI - {{\hat A}_1})^{ - 1}}{B_1}{e^{\frac{{-(c_1 + s)}}{{{q_1}}}}}{C_0}\hat Z(s),\label{eq:rf5}\\
&h(s)\int_0^1 {{M_\beta }} (y)\beta (y,s)dy \notag\\
&= \int_0^1 {{M_\beta }} (y)q{e^{\frac{{-(c_2 + s)}}{{{q_2}}}(1 - y) - \frac{{({c_1} + s)}}{{{q_1}}}}}dy{C_0}\hat Z(s),\\
&h(s)\int_0^1 {{M_\alpha }} (y)\alpha (y,s)dy\notag\\
&= \int_0^1 {{M_\alpha }} (y){e^{\frac{{-(c_1 + s)}}{{{q_1}}}y}}dy{C_0}\hat Z(s).\label{eq:rf}
\end{align}
Multiplying both sides of \eqref{eq:ZODEs} by scalar
$h(s)$, and substituting \eqref{eq:r1}-\eqref{eq:rf} therein yields
\begin{align}
&h(s)(sI - {{\hat A}_0})\hat Z(s)\notag\\
=& h(s){q_1}{C_0}^ + {{\bar K}_1}(0){C_0}\hat Z(s) \notag\\
&+ {{M}_Y}{(sI - {{\hat A}_1})^{ - 1}}{B_1}{e^{\frac{{-(c_1 + s)}}{{{q_1}}}}}{C_0}\hat Z(s)\notag\\
 &+ \int_0^1 {{M_\alpha }} (y){e^{\frac{{-(c_1 + s)}}{{{q_1}}}y}}dy{C_0}\hat Z(s) \notag\\
 &+ \int_0^1 {{M_\beta }} (y)q{e^{\frac{{-(c_2 + s)}}{{{q_2}}}(1 - y) - \frac{{({c_1} + s)}}{{{q_1}}}}}dy{C_0}\hat Z(s)\notag\\
 &+ {N_1}{e^{\frac{{-(c_1 + s)}}{{{q_1}}}}}{C_0}\hat Z(s) + {N_2}q{e^{\frac{{-(c_2 + s)}}{{{q_2}}} - \frac{{({c_1} + s)}}{{{q_1}}}}}{C_0}\hat Z(s)\notag\\
 & + h(s){B_0}\bar U(s).\label{eq:ZODEs1}
\end{align}
Recalling Assumption \ref{as:pq}, we know $h(s)$ is nonzero for any $s\in\mathbb{C}$, ${\Re}(s)\ge 0$ and then $h(s)$ has an inverse $h(s)^{-1}$. Multiplying both sides of \eqref{eq:ZODEs1} by $h(s)^{-1}$  and defining
\begin{align}
\hat\xi(t)=C_0\hat Z(t),
\end{align}
\eqref{eq:ZODEs1} is then rewritten as
\begin{align*}
&(sI - {{\hat A}_0})\hat Z(s)\notag\\
 =& {q_1}{C_0}^ + {{\bar K}_1}(0)\hat \xi(s) \\
 &+ h(s)^{-1}{{M}_Y}{(sI - {{\hat A}_1})^{ - 1}}{B_1}{e^{\frac{{-(c_1 + s)}}{{{q_1}}}}}\hat \xi(s) \\
& + h(s)^{-1}\int_0^1 {{M_\alpha }} (y){e^{\frac{{-(c_1 + s)}}{{{q_1}}}y}}dy\hat \xi(s) \\
 &+ h(s)^{-1}\int_0^1 {{M_\beta }} (y)q{e^{\frac{{-(c_2 + s)}}{{{q_2}}}(1 - y) - \frac{{({c_1} + s)}}{{{q_1}}}}}dy\hat \xi(s) \\
 &+ h(s)^{-1}{N_1}{e^{\frac{{-(c_1 + s)}}{{{q_1}}}}}\hat \xi (s)\\
 &+ h(s)^{-1}{N_2}q{e^{\frac{{-(c_2 + s)}}{{{q_2}}} - \frac{{({c_1} + s)}}{{{q_1}}}}}\hat \xi(s)+ {B_0}\bar U(s)
\end{align*}
for any $s\in\mathbb{C}$, ${\Re}(s)\ge 0$. Defining
\begin{align}
&G(s)={q_1}{C_0}^ + {{\bar K}_1}(0)+ h(s)^{-1}\bigg[{{M}_Y}{(sI - {{\hat A}_1})^{ - 1}}{B_1}{e^{\frac{{-(c_1 + s)}}{{{q_1}}}}}\notag\\
 &+ \int_0^1 {{M_\alpha }} (y){e^{\frac{{-(c_1 + s)}}{{{q_1}}}y}}dy + \int_0^1 {{M_\beta }} (y)q{e^{\frac{{-(c_2 + s)}}{{{q_2}}}(1 - y) - \frac{{({c_1} + s)}}{{{q_1}}}}}dy \notag\\
 &+ {N_1}{e^{\frac{{-(c_1 + s)}}{{{q_1}}}}} + {N_2}q{e^{\frac{{-(c_2 + s)}}{{{q_2}}} - \frac{{({c_1} + s)}}{{{q_1}}}}}\bigg]\label{eq:Gs}
\end{align}
which is a stable, proper transfer matrix, we have
\begin{align}
(sI - {{\hat A}_0})\hat Z(s) = G(s)\hat \xi(s)+ {B_0}\bar U(s).\label{eq:ZODEs2}
\end{align}
Recalling $\hat A_0$ being Hurwitz, ${{\rm det}}(sI - {{\hat A}_0})$ does not have any zeros in the closed
right-half plane. Then the matrix $(sI - {{\hat A}_0})$ is invertible for any $s\in\mathbb{C}$, ${\Re}(s)\ge 0$. Multiplying both sides of \eqref{eq:ZODEs2} $C_0(sI - {{\hat A}_0})^{-1}$, we obtain
\begin{align}
C_0\hat Z(s) =& C_0(sI - {{\hat A}_0})^{-1}G(s)\hat \xi(s)\notag\\
&+ C_0(sI - {{\hat A}_0})^{-1}{B_0}\bar U(s).
\end{align}
That is
\begin{align}
\hat \xi(s) = C_0(sI - {{\hat A}_0})^{-1}G(s)\hat \xi(s)+ {W_0}\bar U(s),\label{eq:ZODEs3}
\end{align}
where ${W_0}(s)=C_0(sI - {{\hat A}_0})^{-1}{B_0}$.

Recalling Assumption \ref{ABC0} which is equivalent to the existence of a right inverse for $W_0$. A
possible choice is given by the Moore-Penrose right inverse $W_0^+(s)=W_0^T(s)(W_0(s)W_0^T(s))^{-1}$ \cite{Saba2019}.

Choose $\bar U(s)$ in \eqref{eq:ZODEs3} as
\begin{align}
\bar U(s)&=-W_0^+(s)\Omega(s)C_0(sI - {{\hat A}_0})^{-1}G(s)\hat \xi(s)\notag\\
&=F(s)\hat \xi(s)\label{eq:Us}
\end{align}
where a SISO low-pass filter $\Omega(s)$ satisfying
\begin{align}
|1-\Omega(j\omega )|<\frac{1}{\sup_{\omega\in R}\bar\sigma(G(j\omega))\bar\sigma(C_0(j\omega I- {{\hat A}_0})^{-1})}, \forall \omega\in R\label{eq:w}
\end{align}
is adopted to make sure $F(s)$ strictly proper. Note that because $G(s)$ is uniformly bounded in the closed
right-half plane, $\sup_{\omega\in R}\bar\sigma(G(j\omega))$ is bounded where $\bar\sigma$ stands for the largest singular value. A low-pass filter $\Omega(s)$ always can be chosen to ensure $F(s)$ strictly proper and satisfy \eqref{eq:w} concurrently, because there exists a $\omega_1$  to make the right hand side of \eqref{eq:w} larger than 1 at $\omega\ge\omega_1$ ($\sup_{\omega\in R}\bar\sigma(G(j\omega))$ is bounded and $\bar\sigma(C_0(j\omega - {{\hat A}_0})^{-1})$ can be small enough at sufficiently high frequencies), and thus \eqref{eq:w} still holds  even if the gain $|\Omega(j\omega )|$ of the low-pass filter is close to zero at $\omega\ge\omega_1$. It means that a choice of the cut-off frequency of the low-pass filter $\Omega(s)$ is $\omega_1$.

Note that $\bar U$ has been chosen as strictly
proper by introducing the low-pass filter $\Omega(s)$, which means that the controller is robust
to small input delays \cite{Saba2019}.

Substituting \eqref{eq:Us} into \eqref{eq:ZODEs3}, we have
\begin{align}
\hat \xi(s) &= (1-\Omega(s))C_0(sI - {{\hat A}_0})^{-1}G(s)\hat \xi(s)\notag\\
&=\Phi(s)\hat \xi(s).
\end{align}
That is
\begin{align}
(1-\Phi(s))\hat \xi(s) = 0,\label{eq:hxi1}
\end{align}
where
\begin{align}
\bar\sigma\left(\Phi(j\omega ) \right)&\le |1-\Omega(j\omega )|\bar\sigma(C_0(j\omega I - {{\hat A}_0})^{-1})\sup_{\omega\in R}\bar\sigma(G(j\omega ))\notag\\
&<1\label{eq:hxi2}
\end{align}
by recalling \eqref{eq:w}, which is a sufficient condition for exponential convergence to zero of $\hat \xi$.
Considering \eqref{eq:Us}, \eqref{eq:bU}, $U(s)$ can be written as
\begin{align}
U(s)&=\bar U(s)+F_0\hat Z(s)\notag\\
&=[{F_0}-W_0^+(s)\Omega(s)C_0(sI - {{\hat A}_0})^{-1}G(s)C_0]\hat Z(s),\label{eq:U}
\end{align}
where inverse Laplace transform is required to represent $U(s)$ in the time domain considering implementation of the controller and  $\hat Z$ can be replaced as the observer states by \eqref{eq:contran2}, \eqref{eq:contran1a}-\eqref{eq:contran1b}.
\section{Stability analysis of the closed-loop system}\label{sec:stability}
\begin{figure}
\centering
\includegraphics[width=8.5cm]{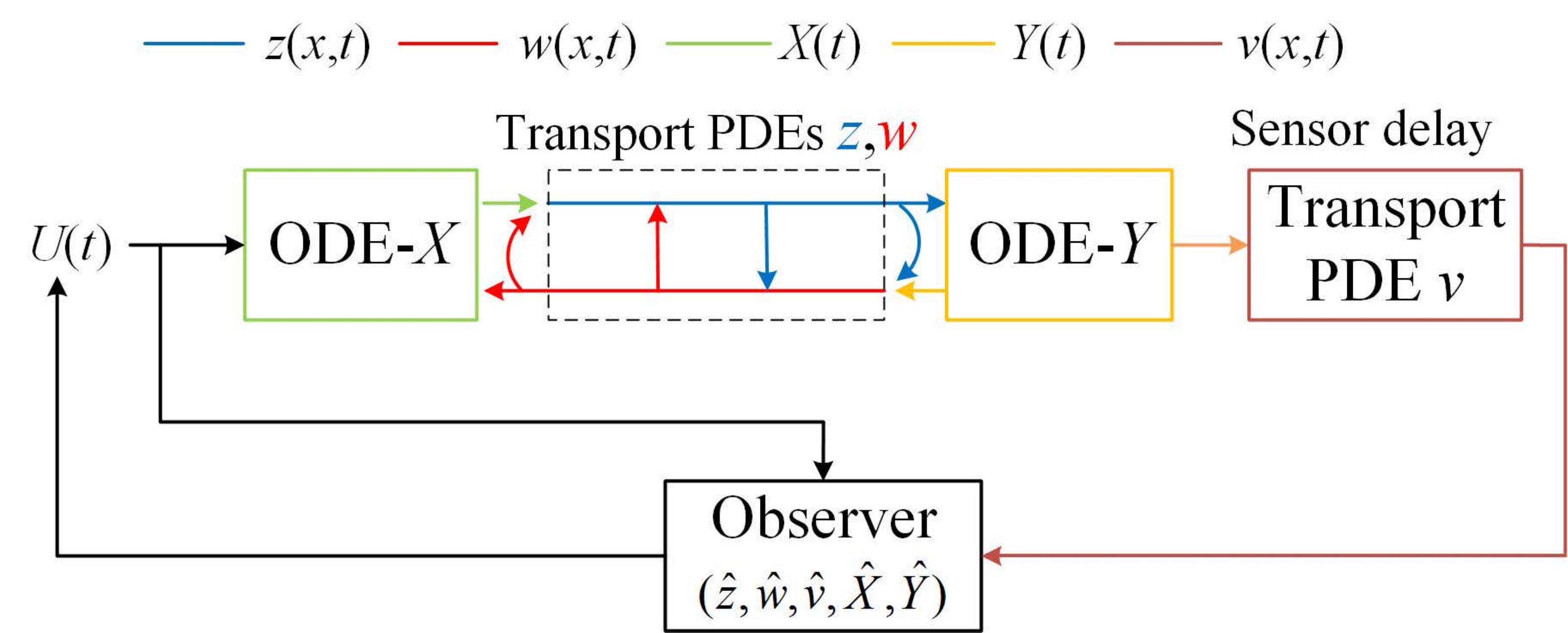}
\caption{Diagram of the closed-loop system.}
\label{fig:1}
\end{figure}

The closed-loop system includes the plant \eqref{eq:plant1}-\eqref{eq:plant7}, the observer \eqref{eq:observer1}-\eqref{eq:observer8} and the controller \eqref{eq:U}. The block diagram of the closed-loop system is shown in Fig. \ref{fig:1}.
We have given Theorem \ref{thm:1} showing the observer error states between the plant and the observer are exponentially convergent to zero in the sense of the norm \eqref{eq:errornorm} in Section \ref{sec:stabilityerror}. Considering \eqref{eq:errorstates}, in order to prove the exponential stability result of the closed-loop system, we present the next lemma to show the exponential stability of the system-$(\hat z(x,t), \hat w(x,t), \hat v(x,t), \hat X(t), \hat Y(t))$ \eqref{eq:observer1}-\eqref{eq:observer8} under the controller \eqref{eq:U}.
\begin{lem}\label{lem:1}
For any initial data $(\hat z(x,0),\hat w(x,0),\hat v(x,0), \hat X(0)$, $\hat Y(0))\in L^2(0,1)\times L^2(0,1)\times L^2(1,2)\times \mathbb{R}^{n}\times \mathbb{R}^{m}$, exponential stability of the system \eqref{eq:observer1}-\eqref{eq:observer8} under the controller \eqref{eq:U} holds in the sense of the norm
\begin{align*}
\|\hat z(\cdot,t)\|_{\infty}+\|\hat w(\cdot,t)\|_{\infty}+\|\hat v(\cdot,t)\|_{\infty}+\left|\hat X(t)\right|+\left|\hat Y(t)\right|.
\end{align*}
with the convergence rate being adjustable by $F_0,F_1$.
\end{lem}
\textbf{{Proof.}}
We would prove the exponential convergence of the states in the overall system based on the exponential convergence of $\hat\xi(t)$ by applying their algebraic relationships obtained in Section \ref{eq:cdesign}.

According to exponential convergence to zero of $\hat\xi(t)=C_0\hat Z(t)$, which is obtained from \eqref{eq:hxi1}-\eqref{eq:hxi2}, recalling \eqref{eq:afZ}-\eqref{eq:Bevx} and \eqref{eq:rf5}-\eqref{eq:rf}, we have $\alpha(x,t)$, $\beta(x,t)$, $\hat v(x,t)$, $\|\alpha(\cdot,t)\|$, $\|\beta(\cdot,t)\|$, $|\hat Y(t)|$ are exponentially convergent to zero, where the convergence rate is adjustable by $F_0,F_1$ considering \eqref{eq:hatA0}-\eqref{eq:bA1}.

Substituting \eqref{eq:Us} into \eqref{eq:ZODEs2},
\begin{align}
&\quad\hat Z(s)\notag\\
&= (sI - {{\hat A}_0})^{-1}[G(s)-{B_0}W_0^+\Omega(s)C_0(sI - {{\hat A}_0})^{-1}G(s)]\hat \xi(s).\label{eq:ZODEs0}
\end{align}
Because $(sI - {{\hat A}_0})^{-1}[G(s)-{B_0}W_0^+\Omega(s)C_0(sI - {{\hat A}_0})^{-1}G(s)]$ is a (stable) proper transfer matrix, using the exponential convergence result of $\hat \xi$, we also obtain exponential convergence to zero of $\hat Z$ via \eqref{eq:ZODEs0}.

Applying Cauchy-Schwarz inequality into the inverse transformations \eqref{eq:contran1aI}-\eqref{eq:contran1bI}, and transformations \eqref{eq:contran2}, we obtain
\begin{align}
&\left|\hat z (x,t)\right|+\left|\hat w (x,t)\right|\le \Upsilon_{2a}\bigg( |\alpha(x,t)| +|\beta(x,t)|+\|\alpha(\cdot,t)\| \notag\\
&\quad\quad\quad\quad\quad\quad\quad\quad+\|\beta(\cdot,t)\|+ \left|\hat Y(t)\right|\bigg),\notag\\
&\left|\hat X(t)\right|\le \Upsilon_{2b}\bigg( \|\alpha(\cdot,t)\| +\|\beta(\cdot,t)\|+\left|\hat Y(t)\right|+ \left|\hat Z(t)\right|\bigg)
\end{align}
for some positive  $\Upsilon_{2a},\Upsilon_{2b}$. Recalling the obtained exponential convergence of $\alpha(x,t)$, $\beta(x,t)$, $\|\alpha(\cdot,t)\|$, $\|\beta(\cdot,t)\|$, $|\hat Y(t)|$, $|\hat Z(t)|$, we thus obtain the exponential convergence to zero of $\hat z (x,t)+\hat w (x,t)+|\hat X(t)|$. Recalling the exponential convergence to zero of $|\hat Y(t)|$ and $\hat v(x,t)$, we obtain Lemma \ref{lem:1}.\QEDA
\begin{thme}\label{thm:main}
For any initial data $(z(x,0),w(x,0)$, $v(x,0)$, $X(0), Y(0))\in L^2(0,1)\times L^2(0,1)\times L^2(1,2)\times \mathbb{R}^{n}\times \mathbb{R}^{m}$, considering the closed-loop system including the plant \eqref{eq:plant1}-\eqref{eq:plant7}, the observer \eqref{eq:observer1}-\eqref{eq:observer8} and the controller \eqref{eq:U},

1) The internal exponential stability holds in the sense of the norm
\begin{align*}
&\|z(\cdot,t)\|_{\infty}+\|w(\cdot,t)\|_{\infty}+\|v(\cdot,t)\|_{\infty}+|X(t)|+|Y(t)|\notag\\
&\|\hat z(\cdot,t)\|_{\infty}+\|\hat w(\cdot,t)\|_{\infty}+\|\hat v(\cdot,t)\|_{\infty}+\left|\hat X(t)\right|+\left|\hat Y(t)\right|\notag\\
&+|y_1(t)|+|y_4(t)|+\|y_2(\cdot,t)\|_{\infty}+\|y_3(\cdot,t)\|_{\infty}+\|y_5(\cdot,t)\|_{\infty}
\end{align*}
with the convergence rate being adjustable by $L_0,L_1,F_0,F_1$.

2) There exist positive constants $\Gamma_{c}$ and $\lambda_c$ making the dynamic feedback control $U(t)$ bounded and exponentially convergent to zero in the sense of
\begin{align*}
|U(t)|\le \Gamma_{c}e^{-\lambda_ct}.
\end{align*}
\end{thme}
\textbf{{Proof.}}
1) Applying \eqref{eq:errorstates} and Cauchy-Schwarz inequality, recalling Theorem \ref{thm:1} and Lemma \ref{lem:1}, we straightforwardly obtain 1) in Theorem \ref{thm:main}.

2) According to the control design in Section \ref{eq:cdesign}, we know $F(s)=W_0^+\Omega(s)C_0(sI - {{\hat A}_0})^{-1}G(s)C_0$ in \eqref{eq:Us} is strictly proper. It follows that ${F_0}-W_0^+\Omega(s)C_0(sI - {{\hat A}_0})^{-1}G(s)C_0$ in \eqref{eq:U} is a (stable) proper transfer function because $F_0$ is a constant matrix. Recalling \eqref{eq:U} and the exponential convergence of $\hat Z$ proved in Lemma \ref{lem:1}, we obtain the exponential convergence to zero of the dynamic feedback control $U(t)$, which is a dynamic extension generated by utilizing the frequency-domain design approach.

The proof of Theorem \ref{thm:main} is completed.\QEDA
\section{Application in control of a deepwater construction vessel}\label{sim}
\begin{figure}
\centering
\includegraphics[width=5cm]{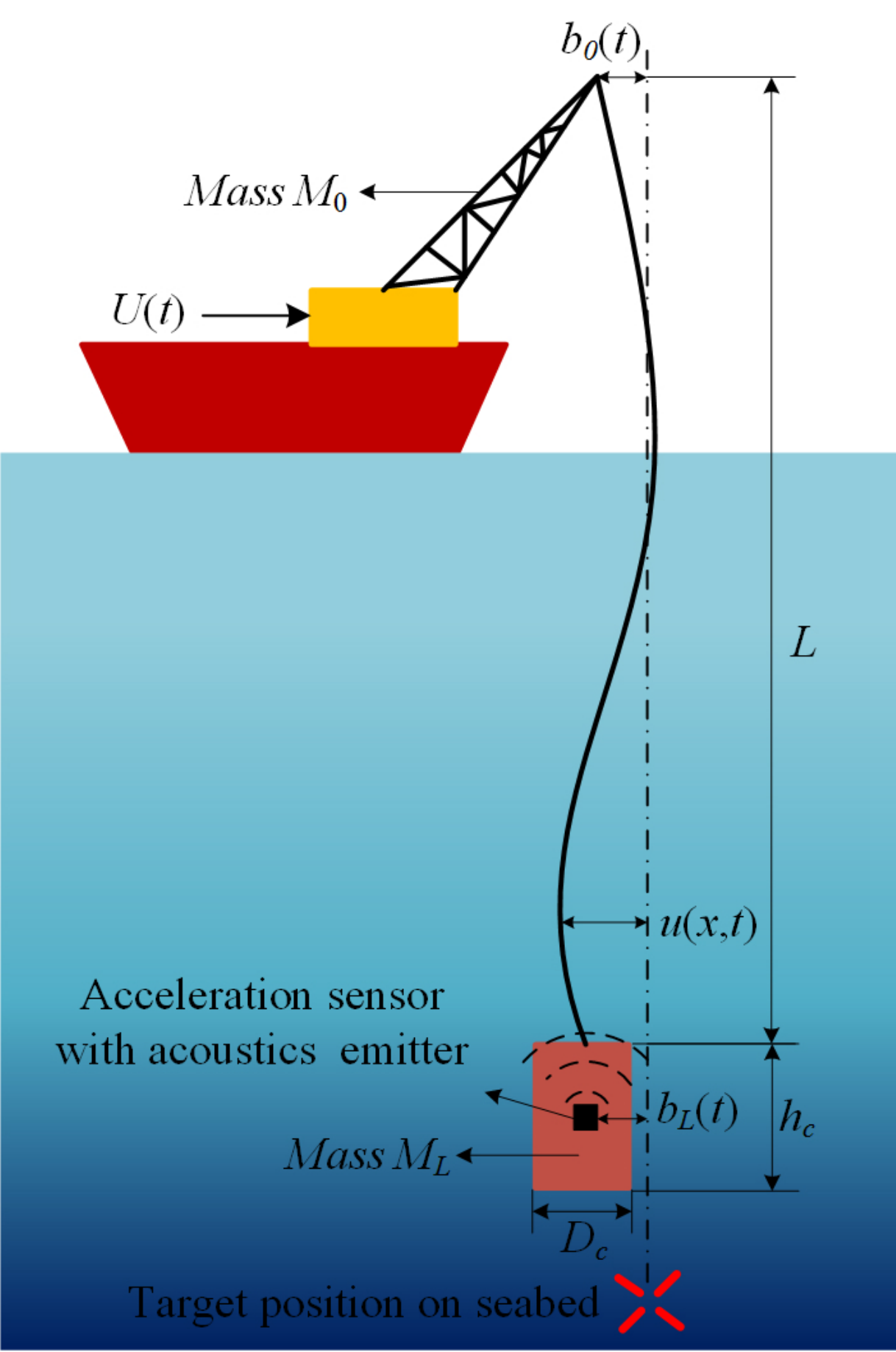}
\caption{Schematic of a DCV used in seafloor installation.}
\label{fig:DCV}
\end{figure}

A DCV is used to place equipment to be installed at the predetermined location on the seafloor for off-shore oil drilling, which is shown in Fig. \ref{fig:DCV} and described in the first paragraph in Section \ref{sec:Intro}. The equipment, referred to as payload, have to be installed accurately at the predetermined location with a tight tolerance, such as the permissible maximum tolerance for a typical
subsea installation in  \cite{How2011} is 2.5 m. In this section, we design an output-feedback control force at the crane to reduce oscillations
of the long cable and position the payload in the target area with compensating the sensor delay, where the details of applying the above theoretical results in observer and controller designs would be presented. Note that we only consider one-dimensional oscillations of DCV and the end phase of the descending process, i.e., the cable length being constant. Control problems of two-dimensional coupled oscillations of DCV in the whole descending/ascending process with a time-varying-length cable are considered in \cite{J2020Arxiv} which, however, is not a sandwiched system by neglecting the crane dynamics, and does not include delay compensation.
\subsection{Modeling}
A nonlinear PDE model of the DCV consisting of vessel, crane, cable,
and payload has been presented in \cite{How2011}.
\subsubsection{Simplifies from a nonlinear model of the vessel crane}
Impose the following simplifies on the nonlinear model of the DCV, Eqs. (1), (2), (6)-(8)  in \cite{How2011}:
\begin{itemize}
\item Nonuniform distributed tension $T(z,t)$ in Eq. (2) in \cite{How2011}, which introduces nonlinearity into the model, is simplified as an uniform distributed tension $T_0$ which is defined by an equation of static equilibrium.
    \item The dynamics of the ocean surface vessel $y_s(t)$  (Eq. (1) in \cite{How2011}) is neglected and only consider regulating dynamics of crane-cable-payload, because the vessel can be kept at the desired position by the ship dynamic positioning system.
        \item The time-varying damping coefficients $d_0(t)$, $d_L(t)$ of the onboard crane and the payload in \cite{How2011} are considered as constants.
\end{itemize}
Thus the  nonlinear model of the DCV (Eqs. (1), (2), (6)-(8)  in \cite{How2011}) can be simplified as the following linear model (note that $u_0(t),u_L(t)$ in (7)-(8) in \cite{How2011} are two control inputs which are not required in our design):
\begin{align}
&M_0\ddot b_0(t) = -d_0\dot b_0 (t)+ T_0u_x(0,t)+U(t),\label{eq:5.1}\\
&u(0,t) = b_0(t),\label{eq:5.2}\\
&\rho u_{tt}(x,t) = T_0u_{xx}(x,t) - d_c{u_t}(x,t)+f(x,t), x\in[0,L]\label{eq:5.3}\\
&u(L,t) = b_L(t),\label{eq:5.4}\\
&M_L\ddot b_L(t) = -d_L\dot b_L(t) + T_0u_x(L,t)+f_L(t),\label{eq:5.5}
\end{align}
Therein, $T_0$ is static tension defined as $T_0=M_Lg-F_{buoyant}$ with $F_{buoyant}=\frac{1}{4}\pi D_c^2 h_c\rho_{s}g$. $u(x,t)$ denote distributed transverse displacements along the cable. $b_0(t)$ and $b_L(t)$ represent transverse displacements of the onboard crane and the payload. $f(z,t),f_L(t)$ are ocean current disturbances , i.e., external drag forces at the cable and payload.  The physical parameters of the DCV in simulation are from \cite{How2011} and shown in Tab. \ref{table1}.

Note that even though the DCV model here is a linear model while that in \cite{How2011} is a more complicated nonlinear model, some boundedness assumptions in \cite{How2011} are not required here. Moreover, only one control input at the onboard crane is required here while one more control input for the payload is required in the control system in \cite{How2011}.
\begin{table}
\centering
\caption{Physical parameters of the DCV.}
\begin{tabular}{lccc}
\hline
Parameters (units)&values\\ \hline
Cable length ${L}$ (m) &1000\\
Cable diameter ${R_D}$ (m) &0.2\\
Cable effective Young¡¯s Modulus $E$ (N/m$^{2}$) &4.0$\times$${10^{9}}$\\
Cable linear density ${\rho}$ (kg/m) &8.02\\
Crane mass ${M_0}$ (kg) &1.0$\times$ $10^{6}$\\
Payload mass ${M_L}$ (kg) &4.0$\times$ $10^5$\\
Gravitational acceleration $g$ (m/s$^{2}$) &9.8\\
Cable material damping coefficient $d_c$ (N$\cdot$s/m) & 0.5\\
Height of payload modeled as a cylinder $h_c$ (m)& 10\\
Diameter of payload modeled as a cylinder $D_c$ (m)& 5\\
Damping coefficient at payload $d_L$ (N$\cdot$s/m) & 2.0$\times$${10^{5}}$\\
Damping coefficient at crane $d_0$ (N$\cdot$m$\cdot$s/rad) & 8.0$\times$ ${10^{5}}$\\
Seawater density $\rho_s$ (kgm$^{-3}$)& $1024$  \\\hline
\end{tabular}
\label{table1}
\end{table}
\subsubsection{Reformulation of the linear model}
Apply Riemann transformations
\begin{align}
{z}(x,t) &= {u_t}(x,t) - \sqrt {\frac{T_0}{\rho}} {u_x}(x,t),\label{eq:t11}\\
{w}(x,t) &= {u_t}(x,t) + \sqrt {\frac{T_0}{\rho}} {u_x}(x,t)\label{eq:t12}
\end{align}
and define new variables $X(t)=\dot b_0(t)$, $Y(t)=\dot b_L(t)$, \eqref{eq:5.1}-\eqref{eq:5.5} can be rewritten  as
\begin{align}
&\dot X(t) =A_0X(t) + E_0{w}(0,t)+B_0U(t) ,\label{eq:x1}\\
&z(0,t)= pw(0,t)+C_0X(t),\label{eq:x4}\\
&{z_{t}}(x,t) =  - q_1 {z_{x}}(x,t)- c_1({z}(x,t) + {w}(x,t))+f(x,t),\label{eq:x2}\\
&{w_{t}}(x,t) = q_2{w_{x}}(x,t) - c_2({z}(x,t) + {w}(x,t))+f(x,t),\label{eq:x3}\\
&w(L,t) = qz(L,t)+C_1Y(t),\label{eq:x5}\\
&\dot Y(t) = A_1Y(t) + B_1z(1,t)+f_L(t),\label{eq:5.4b}\\
&y_{\rm out}(t) = C_1Y(t-\tau),\label{eq:5.6b}
\end{align}
where $y_{\rm out}(t)$ is the delayed measurement output and $y_{\rm out}(t)=0,t\in[0,\tau)$ because the sensing signal has not been received. The observer and controller design in the next subsection is based on \eqref{eq:x1}-\eqref{eq:5.6b} except for the disturbances $f(x,t),f_L(t)$, which are regarded as model uncertainties in the simulation to test the robustness of the controller. The sensor delay $\tau$ is considered as $0.1 s$. Note that $q_1=q_2=\sqrt {\frac{T_0}{\rho}}$, $c_1=c_2=\frac{{{d_c}}}{2\rho}$. $p,q$ satisfy Assumption \ref{as:pq} ($|pq|=1<e^{\frac{{{c_2}}}{{{q_2}}} + \frac{{{c_1}}}{{{q_1}}}}=1.0014$), and
\begin{align}
&A_0=\frac{-d_0}{M_0}-\frac{\sqrt{T_0\rho}}{M_0},~E_0=\frac{\sqrt {T_0\rho}}{M_0},~B_0=\frac{1}{M_0},C_0=2,\label{eq:ABC0}\\ &A_1=\frac{-d_L}{M_L}+\frac{\sqrt{T_0\rho}}{M_L},~B_1=-\frac{\sqrt {{T_0}\rho}}{M_L},~C_1=2 \label{eq:ABC1}
\end{align}
satisfy Assumptions \ref{as:controllable}-\ref{ABC1a}.
\subsubsection{Ocean current disturbances in the model}
The time-varying ocean surface current velocity can be modeled by a first-order Gauss-Markov process \cite{Fossen2002}
\begin{align}
\dot P(t)+\mu P(t)=\mathcal G(t),~~ P_{\min}\le P(t)\le P_{\max},
\end{align}
where $\mathcal G(t)$ is Gaussian white noise. Constants $P_{\min},P_{\max}$ and $\mu$ are chosen as $1.6 ms^{-1}$, $2.4 ms^{-1}$ and 0 \cite{How2011}. The full current load $P(t)$ is applied at the cable from $x=0$ to $x=300m$  and
thereafter linearly decline to $0.1P(t)$ at the bottom of the cable, i.e., $x=1000$ m \cite{How2011}. The depth dependent ocean current profile $P(x,t)$ is thus obtained as \begin{equation}
P(x,t)=\left\{
\begin{aligned}
& P(t),~~0\le x\le 300\\
& \frac{970-0.9x}{700}P(t),~~300\le x\le L
\end{aligned}
\right.
\end{equation}
which determines the ocean current disturbances $f(x,t)$ $f_L(t)$ as following. $f(x,t)$ can be modeled as an oscillating drag force \cite{How2011}:
\begin{align}
f(x,t)=\frac{1}{2}\rho_{s}C_dP(x,t)^2R_DA_D\cos\left(4\pi\frac{S_tP(x,t)}{R_D}t+\varsigma\right)\label{eq:fxt}
\end{align}
where $C_d=1$ denoting the drag coefficient, $\varsigma=\pi$ being the phase angle, $A_D=400$ denoting the amplitude of the oscillating drag force, $S_t=0.2$ being the Strouhal number \cite{Faltinsen1993}.
The drag force $f_L(t)$ at the payload which is considered as a cylinder is derived from Morison's equation \cite{How2011}:
\begin{align}
f_L(t)=\frac{1}{2}C_d\rho_{s}h_cD_c\left|P(L,t)\right|P(L,t).\label{eq:fL}
\end{align}
\subsection{Observer and controller}
The observer and controller would be defined by applying the theoretical results obtained in the previous sections.
\subsubsection{Observer}
Defining control parameters $L_0$ and $L_1$, $\bar A_0=A_0-L_0C_0$ and $\bar A_1=A_1-L_1C_1$ are obtained as
\begin{align}
\bar A_0&=\frac{-d_0}{M_0}-\frac{\sqrt{T_0\rho}}{M_0}-2L_0\label{eq:barA0s},\\
\bar A_1&=\frac{-d_L}{M_L}+\frac{\sqrt{T_0\rho}}{M_L}-2L_1\label{eq:barA1s}
\end{align}
according to \eqref{eq:ABC0}-\eqref{eq:ABC1}, where the control parameters $L_0$ and $L_1$ should satisfy
\begin{align}
L_0>\frac{-d_0}{2M_0}-\frac{\sqrt{T_0\rho}}{2M_0},\label{eq:CL0}\\
L_1>\frac{-d_L}{2M_L}+\frac{\sqrt{T_0\rho}}{2M_L}.\label{eq:CL1}
\end{align}
According to the values of the physical parameters in the DCV given in Tab. \ref{table1}, we know $L_0,L_1$ should satisfy $L_0>-0.4,L_1>-0.25$. Considering the robustness to the external disturbances \eqref{eq:fxt}-\eqref{eq:fL}, $L_0$ and $L_1$ are chosen as $0.05$ and $0.1$ in the simulation after adjusting $L_0$ and $L_1$ under the disturbances.

Next, we calculate $r(s)=C_1e^{-\tau A_1}(sI-{{\bar A}_1})^{-1}{B_1}$ in this DCV model. According to $A_1$ \eqref{eq:ABC1}, we have
\begin{align}
e^{-\tau A_1}=e^{\frac{d_L\tau}{M_L}-\frac{\tau\sqrt{T_0\rho}}{M_L}}.\label{eq:etA1}
\end{align}
According to \eqref{eq:ABC1}, \eqref{eq:barA1s}, \eqref{eq:etA1}, $r(s)$ can then be obtained as
\begin{align}
r(s)&=C_1e^{-\tau A_1}(sI-{{\bar A}_1})^{-1}{B_1}\notag\\
&=\frac{-2\sqrt{T_0\rho}e^{\frac{d_L\tau}{M_L}-\frac{\tau\sqrt{T_0\rho}}{M_L}}}{M_Ls+{d_L}-{\sqrt{T_0\rho}}+2L_1{M_L}}.\label{eq:r1si}
\end{align}
From \eqref{eq:r1si}, we know the numerator and denominator of \eqref{eq:r1si} are not zero for $s\ge0$ recalling the choice of $L_1$ \eqref{eq:CL1}. It means $r(s)\neq0$ and $r(s)$ is bounded.
Next, we calculate the form of $H_1(s)$, $H_2(s;x)$, $H_3(s;x)$, $H_4(s)$, $H_5(s;x)$.

Recalling \eqref{eq:H1}, \eqref{eq:Hg2}, \eqref{eq:H3g3}, \eqref{eq:H4s} and \eqref{eq:H5g5a}, $H_1(s)$, $H_2(s;x)$, $H_3(s;x)$, $H_4(s)$, $H_5(s;x)$ in the observer are obtained
\begin{align}
&H_1(s)=\frac{-{q_1}{K_1}(1)}{C_1e^{-\tau A_1}(sI-{{\bar A}_1})^{-1}{B_1}}\notag\\
&={-{q_1}{K_1}(1)}\frac{M_Ls+{d_L}-{\sqrt{T_0\rho}}+2L_1{M_L}}{-2\sqrt{T_0\rho}e^{\frac{d_L\tau}{M_L}-\frac{\tau\sqrt{T_0\rho}}{M_L}}},\label{eq:H1ss}\\
&H_2(s;x)=\frac{-{q_1}\phi (x,1)}{C_1e^{-\tau A_1}(sI-{{\bar A}_1})^{-1}{B_1}}\notag\\
&={-{q_1}\phi (x,1)}\frac{M_Ls+{d_L}-{\sqrt{T_0\rho}}+2L_1{M_L}}{-2\sqrt{T_0\rho}e^{\frac{d_L\tau}{M_L}-\frac{\tau\sqrt{T_0\rho}}{M_L}}},\\
&H_3(s;x)=\frac{-{q_1}\psi (x,1)}{C_1e^{-\tau A_1}(sI-{{\bar A}_1})^{-1}{B_1}}\notag\\
&={-{q_1}\psi (x,1)}\frac{M_Ls+{d_L}-{\sqrt{T_0\rho}}+2L_1{M_L}}{-2\sqrt{T_0\rho}e^{\frac{d_L\tau}{M_L}-\frac{\tau\sqrt{T_0\rho}}{M_L}}},\\
&H_4(s)=\frac{q+C_1{(sI-{{\bar A}_1})^{-1}{B_1}}}{C_1e^{-\tau A_1}(sI-{{\bar A}_1})^{-1}{B_1}}\notag\\
&=\left(q+\frac{-2\sqrt{T_0\rho}}{M_Ls+{d_L}-{\sqrt{T_0\rho}}+2L_1{M_L}}\right)\notag\\
&\quad\times\frac{M_Ls+{d_L}-{\sqrt{T_0\rho}}+2L_1{M_L}}{-2\sqrt{T_0\rho}e^{\frac{d_L\tau}{M_L}-\frac{\tau\sqrt{T_0\rho}}{M_L}}}\notag\\
&=q\frac{M_Ls+{d_L}-{\sqrt{T_0\rho}}+2L_1{M_L}}{-2\sqrt{T_0\rho}e^{\frac{d_L\tau}{M_L}-\frac{\tau\sqrt{T_0\rho}}{M_L}}}+\frac{1}{e^{\frac{d_L\tau}{M_L}-\frac{\tau\sqrt{T_0\rho}}{M_L}}},\label{eq:H4ss}\\
&H_5(s;x)=C_1e^ {-\tau A_1 x}\Gamma_1-\frac{C_1e^ {-\tau A_1 x}{B_1}}{C_1e^{-\tau A_1}(sI-{{\bar A}_1})^{-1}{B_1}}\notag\\
&=2e^ {(\frac{-d_L\tau}{M_L}+\frac{\tau\sqrt{T_0\rho}}{M_L})(1-x)}L_1\notag\\
&+\frac{-2e^{(\frac{d_L\tau }{M_L}-\frac{\tau\sqrt{T_0\rho}}{M_L})x}\sqrt{T_o\rho}}{M_L}\frac{M_Ls+{d_L}-{\sqrt{T_0\rho}}+2L_1{M_L}}{-2\sqrt{T_0\rho}e^{\frac{d_L\tau}{M_L}-\frac{\tau\sqrt{T_0\rho}}{M_L}}}\notag\\
&=2e^ {(\frac{-d_L\tau}{M_L}+\frac{\tau\sqrt{T_0\rho}}{M_L})(1-x)}L_1\notag\\
&+\frac{e^{(\frac{d_L\tau }{M_L}-\frac{\tau\sqrt{T_0\rho}}{M_L})(x-1)}}{M_L}({M_Ls+{d_L}-{\sqrt{T_0\rho}}+2L_1{M_L}}),\label{eq:H5ss}
\end{align}
where
\begin{align}
C_1(sI-{{\bar A}_1})^{-1}{B_1}=\frac{-2\sqrt{T_0\rho}}{M_Ls+{d_L}-{\sqrt{T_0\rho}}+2L_1{M_L}}\label{eq:CeAB}
\end{align}
is used in calculating \eqref{eq:H4ss}, and
\begin{align*}
&\quad {C_1e^ {-\tau A_1 x}{B_1}}=\frac{-2e^{(\frac{d_L\tau }{M_L}-\frac{\tau\sqrt{T_0\rho}}{M_L})x}\sqrt{T_o\rho}}{M_L},\\
&\quad C_1e^ {-\tau A_1 x}\Gamma_1=C_1e^ {\tau A_1 (1-x)}L_1=2e^ {(\frac{-d_L\tau}{M_L}+\frac{\tau\sqrt{T_0\rho}}{M_L})(1-x)}L_1
\end{align*}
are used in calculating \eqref{eq:H5ss}.

According to \eqref{eq:H1ss}-\eqref{eq:H5ss}, the inverse Laplace transform  \eqref{eq:f1}-\eqref{eq:f5} can then be obtained as the following form
\begin{align}
&h_1({\tilde y}_{out}(t))=\frac{{{q_1}{K_1}(1)}M_L}{2\sqrt{T_0\rho}e^{\frac{d_L\tau}{M_L}-\frac{\tau\sqrt{T_0\rho}}{M_L}}}\dot{\tilde y}_{out}(t)\notag\\
&+\frac{{{q_1}{K_1}(1)}({d_L}-{\sqrt{T_0\rho}}+2L_1{M_L})}{2\sqrt{T_0\rho}e^{\frac{d_L\tau}{M_L}-\frac{\tau\sqrt{T_0\rho}}{M_L}}}{\tilde y}_{out}(t),\label{eq:f1sim}\\
&h_2(\tilde y_{\rm out}(t);x)=\frac{{{q_1}\phi (x,1)}M_L}{2\sqrt{T_0\rho}e^{\frac{d_L\tau}{M_L}-\frac{\tau\sqrt{T_0\rho}}{M_L}}}{\dot{\tilde y}}_{out}(t)\notag\\
&+\frac{{{q_1}\phi (x,1)}({d_L}-{\sqrt{T_0\rho}}+2L_1{M_L})}{2\sqrt{T_0\rho}e^{\frac{d_L\tau}{M_L}-\frac{\tau\sqrt{T_0\rho}}{M_L}}}{{\tilde y}}_{out}(t),\label{eq:f2sim}\\
&h_3(\tilde y_{\rm out}(t);x)=\frac{{{q_1}\psi (x,1)}M_L}{2\sqrt{T_0\rho}e^{\frac{d_L\tau}{M_L}-\frac{\tau\sqrt{T_0\rho}}{M_L}}}{\dot{\tilde y}}_{out}(t)\notag\\
&+\frac{{{q_1}\psi (x,1)}({d_L}-{\sqrt{T_0\rho}}+2L_1{M_L})}{2\sqrt{T_0\rho}e^{\frac{d_L\tau}{M_L}-\frac{\tau\sqrt{T_0\rho}}{M_L}}}{{\tilde y}}_{out}(t),\label{eq:f3sim}\\
&h_4({\tilde y}_{out}(t))=\frac{qM_L}{-2\sqrt{T_0\rho}e^{\frac{d_L\tau}{M_L}-\frac{\tau\sqrt{T_0\rho}}{M_L}}}{\dot{\tilde y}}_{out}(t)\notag\\
&+\frac{q({d_L}-{\sqrt{T_0\rho}}+2L_1{M_L})-2\sqrt{T_0\rho}}{-2\sqrt{T_0\rho}e^{\frac{d_L\tau}{M_L}-\frac{\tau\sqrt{T_0\rho}}{M_L}}}{\tilde y}_{out}(t),\label{eq:f4sim}\\
&h_5(\tilde y_{\rm out}(t);x)={e^{(\frac{d_L\tau }{M_L}-\frac{\tau\sqrt{T_0\rho}}{M_L})(x-1)}}\dot{\tilde y}_{out}(t)\notag\\
&+\bigg[\frac{e^{(\frac{d_L\tau }{M_L}-\frac{\tau\sqrt{T_0\rho}}{M_L})(x-1)}({d_L}-{\sqrt{T_0\rho}}+2L_1{M_L})}{M_L}\notag\\
&+2e^ {(\frac{-d_L\tau}{M_L}+\frac{\tau\sqrt{T_0\rho}}{M_L})(1-x)}L_1\bigg]{\tilde y}_{out}(t)\label{eq:f5sim}
\end{align}
where $\tilde y_{\rm out}(t)=y_{\rm out}(t)-\hat v(L+1,t)$. $\hat v(L+1,t)$ is the boundary state of the observer \eqref{eq:observer1}-\eqref{eq:observer8} for this DCV model, where $\hat z(x,t),\hat w(x,t)$ are in the spatial domain $x\in[0,L]$ and $\hat v(x,t)$ is defined in $x\in [L,L+1]$.

Note that transfer functions \eqref{eq:H1ss}-\eqref{eq:H5ss} are non proper because $s$ in the numerators, corresponding to taking the time derivative of $\tilde y_{\rm out}(t)$, i.e., $\dot{\tilde y}_{out}(t)$ in \eqref{eq:f1sim}-\eqref{eq:f5sim}. In practice, instead of taking time derivative, we place an accelerometer at the bottom of the cable then getting delayed oscillation acceleration $\ddot b_L(t-\tau)$ to obtain $\dot{\tilde y}_{out}(t)=\dot{y}_{out}(t)-\hat v_{t}(L+1,t)=C_1\ddot b_L(t-\tau)-\hat v_{t}(L+1,t), t\ge\tau$, and then calculate ${\tilde y}_{out}(t)$ by integration from $\tau$ to $t$ with the known initial value ${\tilde y}_{out}(\tau)=C_1Y(0)-\hat v(L,0)$. Actually it is equal to multiply \eqref{eq:H1ss}-\eqref{eq:H5ss} by $\frac{1}{s}$ and then the transfer functions are made proper. Thereby, high gain at high frequency which may be caused by the non proper transfer functions is avoided.
\subsubsection{Controller}
Recalling \eqref{eq:Gs}, in the DCV model, $G(s)$ in the controller \eqref{eq:U} is
\begin{align}
&G(s)=\sqrt {\frac{T_0}{2\rho}} {{\bar K}_1}(0)+ \frac{1}{1 - pq{e^{ {{\frac{-d_c}{\sqrt{\rho T_0}}}}- 2\sqrt {\frac{\rho}{T_0}}s}}}\notag\\
&\times\bigg[{\frac{-{{M}_Y}\sqrt{T_0\rho}{e^{\left({\frac{-d_c}{2\sqrt{\rho T_0}}} - \sqrt {\frac{\rho}{T_0}}s\right)}}}{M_Ls+{d_L}-{\sqrt{T_0\rho}}+F_1}}\notag\\
 &+ \int_0^1 {{M_\alpha }} (y)q{e^{\left({\frac{-d_c}{2\sqrt{\rho T_0}}} - \sqrt {\frac{\rho}{T_0}}s\right)y}}dy \notag\\
 &+ \int_0^1 {{M_\beta }} (y)q{e^{\left({\frac{-d_c}{2\sqrt{\rho T_0}}} - \sqrt {\frac{\rho}{T_0}}s\right)(2 - y) }}dy \notag\\
 &+ {N_1}{e^{\left({\frac{-d_c}{2\sqrt{\rho T_0}}} - \sqrt {\frac{\rho}{T_0}}s\right)}} + {N_2}q{e^{{\frac{-d_c}{\sqrt{\rho T_0}}} - 2\sqrt {\frac{\rho}{T_0}}s}}\bigg]\label{eq:Gvessel}
\end{align}
which is a proper transfer function. According to \eqref{eq:hatA0} and \eqref{eq:ABC0}, one obtains
\begin{align}
(sI - {{\hat A}_0})^{-1}=\frac{{M_0}}{{M_0}s+{d_0}+{\sqrt{T_0\rho}}+{F_0}}
\end{align}
where $F_0$ should be chosen as
\begin{align}
{F_0}>{-d_0}-{\sqrt{T_0\rho}}\label{eq:F0sim}
\end{align}
to make sure $\hat A_0<0$.  According to \eqref{eq:bA1}, $F_1$ should be chosen as
\begin{align}
F_1<\frac{d_L}{\sqrt{T_0\rho}}-1\label{eq:F1sim}
\end{align}
to make $\hat A_1<0$. Recalling the parameter values of the DCV given in Tab. \ref{table1}, we know $F_0,F_1$ should be chosen to satisfy $F_0>-8.04\times 10^5,F_1<49.6$ according to \eqref{eq:F0sim}-\eqref{eq:F1sim}. Considering robustness to the external disturbance \eqref{eq:fxt}-\eqref{eq:fL}, $F_0$ and $F_1$ are defined as $8.57\times 10^5$ and $-2.9\times 10^6$ in the simulation after tuning $F_0$ and $F_1$ under the disturbances.

The low-pass filter $\Omega(s)$ in \eqref{eq:Gs} is chosen as a traditional second-order low-pass filter because the non proper transfer function
\begin{align}
W_0^+(s)&=\frac{1}{{{C_0}{(sI - {{\bar A}_0})^{ - 1}}{B_0}}}\notag\\
&=\frac{{M_0}s+{d_0}+{\sqrt{T_0\rho}}+{F_0}}{{2}}
\end{align}
includes first order $s$ in the numerator. The cut-off frequency $\omega_c$ of the low-pass filtering $\Omega(s)$ can be chosen as
 \begin{align}
 \omega_c={2M_{\max}} - \frac{d_0}{M_0}-\frac{\sqrt{T_0\rho}}{M_0}-\frac{F_0}{M_0}, \label{eq:omrgac}
 \end{align}
where constant $M_{\max}=\sup_{\omega\in R}|G(j\omega)|$ recalling \eqref{eq:Gvessel}. \eqref{eq:omrgac} is chosen according to $\frac{\omega_c +\frac{d_0}{M_0}+\frac{\sqrt{T_0\rho}}{M_0}+\frac{F_0}{M_0}}{2M_{\max}}\ge1$ from \eqref{eq:w}, which means when $\omega>\omega_c$ in \eqref{eq:w}, the right hand side is larger than 1 and the gain of the low-pass filter $\Omega(s)$ can be even zero.
\subsection{Simulation results}
We consider the end phase (20 s) of the descending process, i.e., the payload near the seafloor and the cable being the total length $L$, which is the most important and challenging phase because the cable is long and the oscillations would be large.  The simulation is based on \eqref{eq:x1}-\eqref{eq:5.6b} using the finite difference
method with the time step and the space step as $0.001$ s and $0.1$ m respectively. Considering the sensor delay $\tau=0.1s$, the measurement output is the 100-time-steps-earlier one.
\subsubsection{Initial values}
The initial conditions are defined as $z(x,0)=4\sin( \frac{\pi x}{L})$, $w(x,0)=4\cos( \frac{\pi x}{L})$, thereby, $X(0)=2$, $Y(0)=-2$ recalling \eqref{eq:x4}, \eqref{eq:x5} which physically means initial oscillation velocities of the crane and payload. The initial conditions of \eqref{eq:5.1}-\eqref{eq:5.5} are determined based on the initial conditions of \eqref{eq:x1}-\eqref{eq:5.6b}, i.e., $z(x,0),w(x,0)$. The initial oscillation velocity of the cable is $u_t(x,0)=\frac{1}{2}(z(x,0)+w(x,0))=2\sin( \frac{\pi x}{L})+2\cos( \frac{\pi x}{L})$. The initial distributed oscillation displacement of the cable is defined as $u(x,0)=0$, thereby, initial offset of the payload $b_L(0)=0$, and $b_0(0)=0$ recalling \eqref{eq:5.2},\eqref{eq:5.4}. Our task is to reduce the oscillations of cable and place the payload in the target area, namely within the permissible tolerance 2.5 m around the predetermined location \cite{How2011}, by applying the observer-based output-feedback control force at the onboard crane.
\subsubsection{Responses of $z,w,X,Y$}
\begin{figure}
  \centerline{\includegraphics[width=7cm]{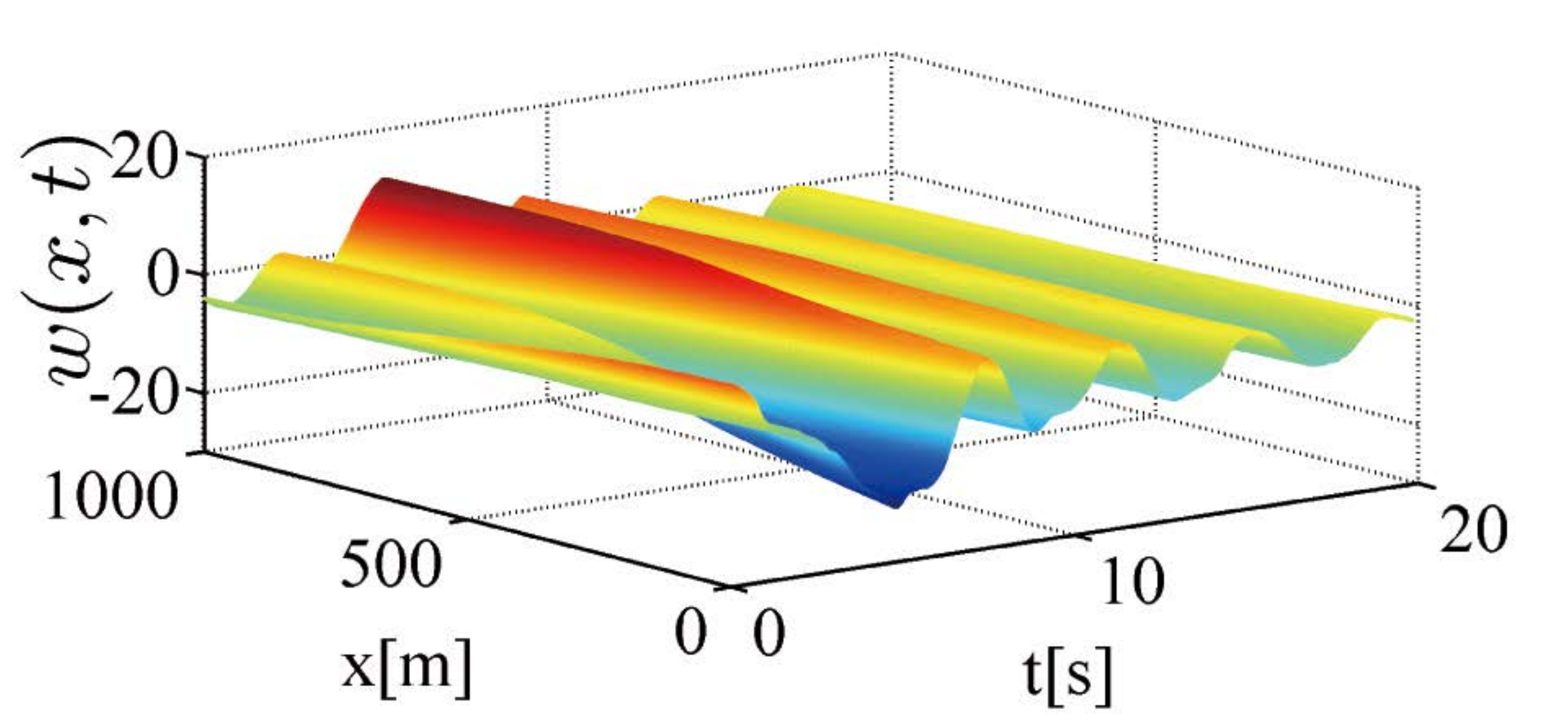}}
  \centerline{(a) $w(x,t)$.}

  \centerline{\includegraphics[width=7cm]{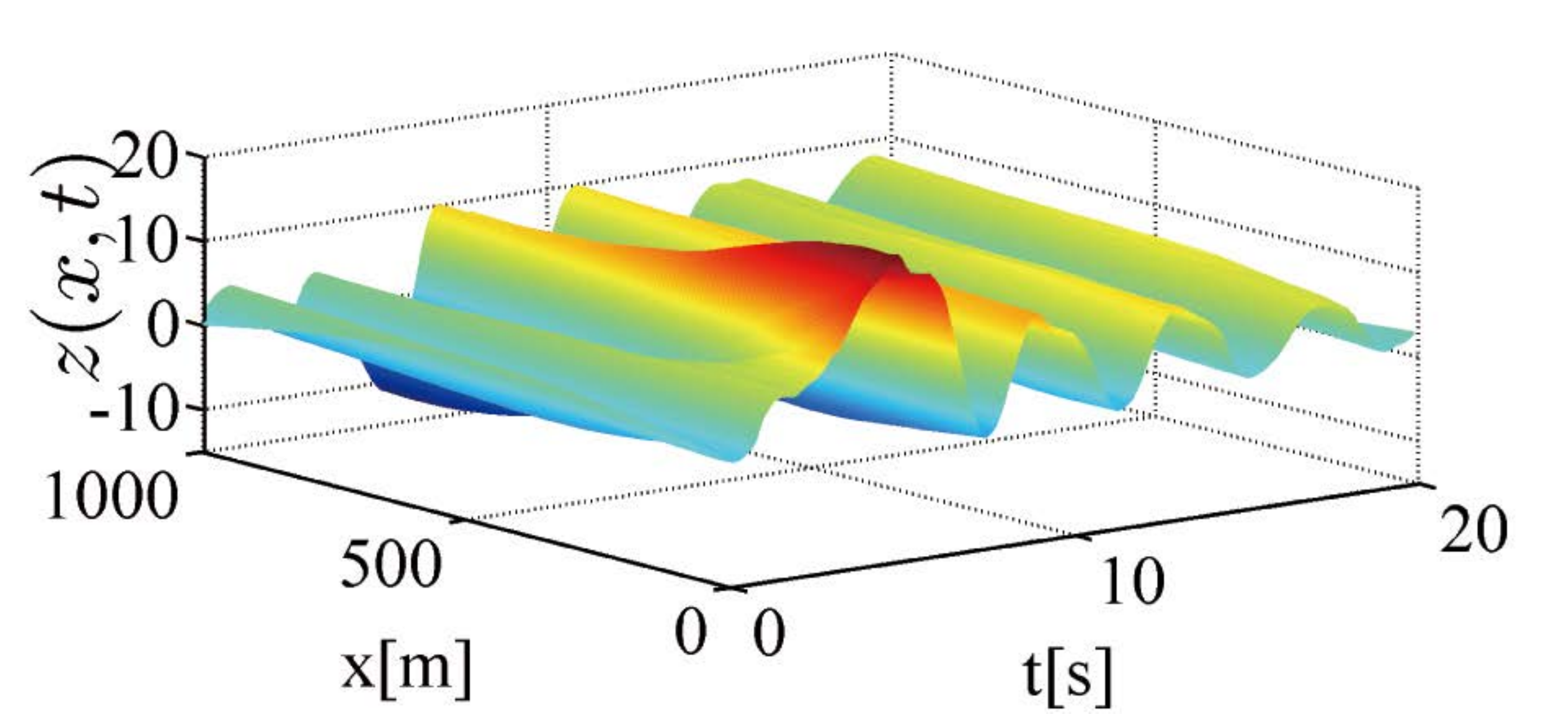}}
  \centerline{(b) $z(x,t)$.}

\caption{Responses of $w(x,t)$, $z(x,t)$ (without control)}
\label{fig:simwzopen}
\end{figure}
\begin{figure}
\centerline{\includegraphics[width=7cm]{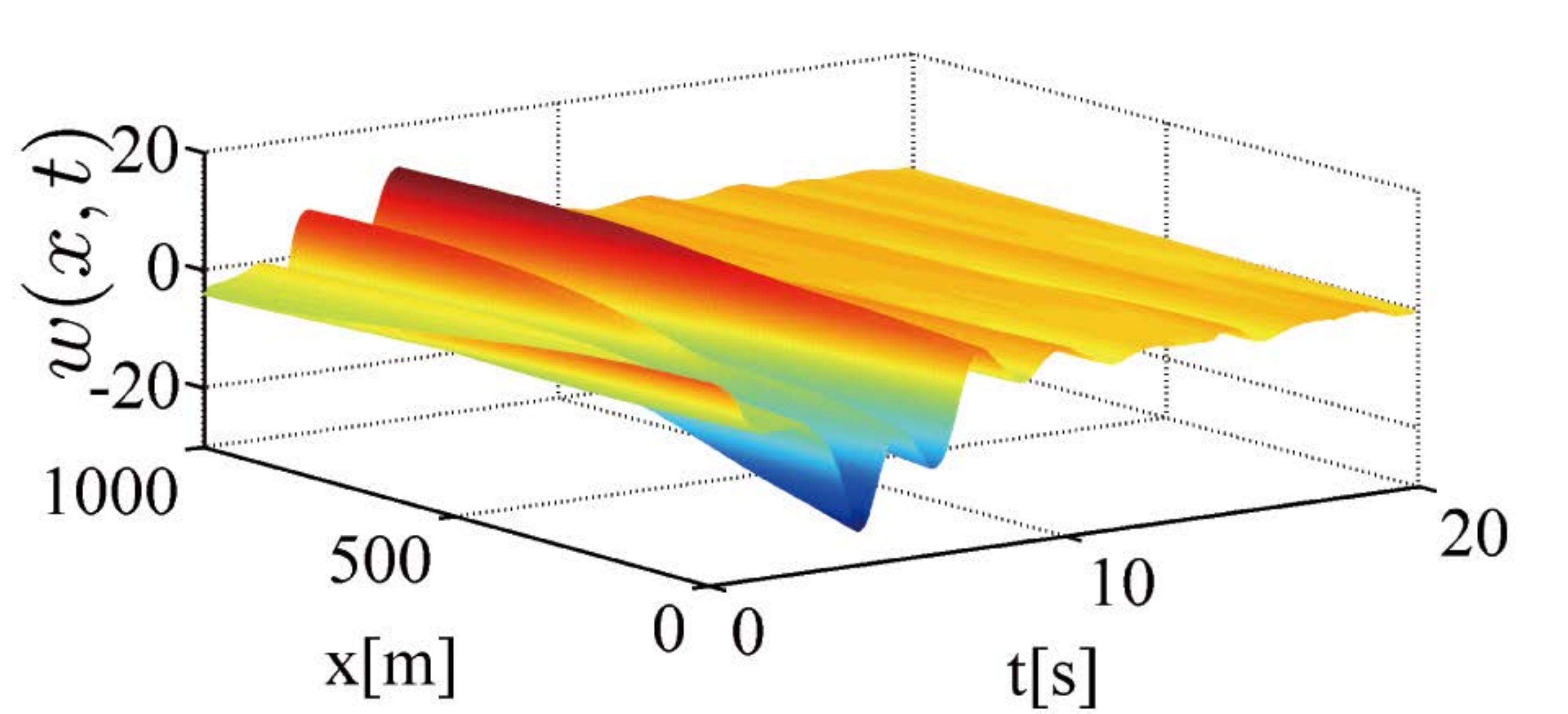}}
  \centerline{(a) $w(x,t)$.}
\centerline{\includegraphics[width=7cm]{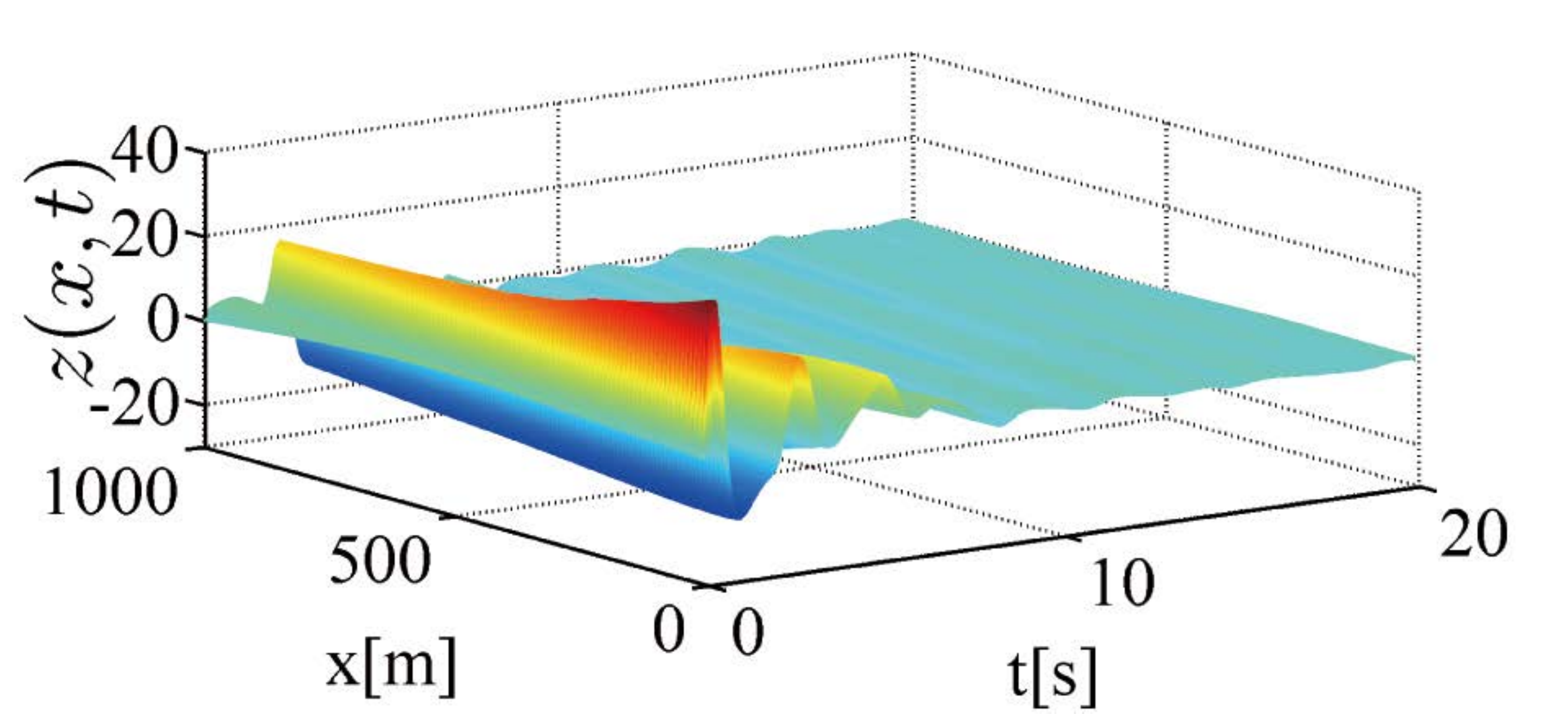}}
  \centerline{(b) $z(x,t)$.}
\caption{Responses of $w(x,t)$, $z(x,t)$ (with control).}
\label{fig:simwzclosed}
\end{figure}
\begin{figure}
\begin{minipage}{0.45\linewidth}
  \centerline{\includegraphics[width=4.5cm]{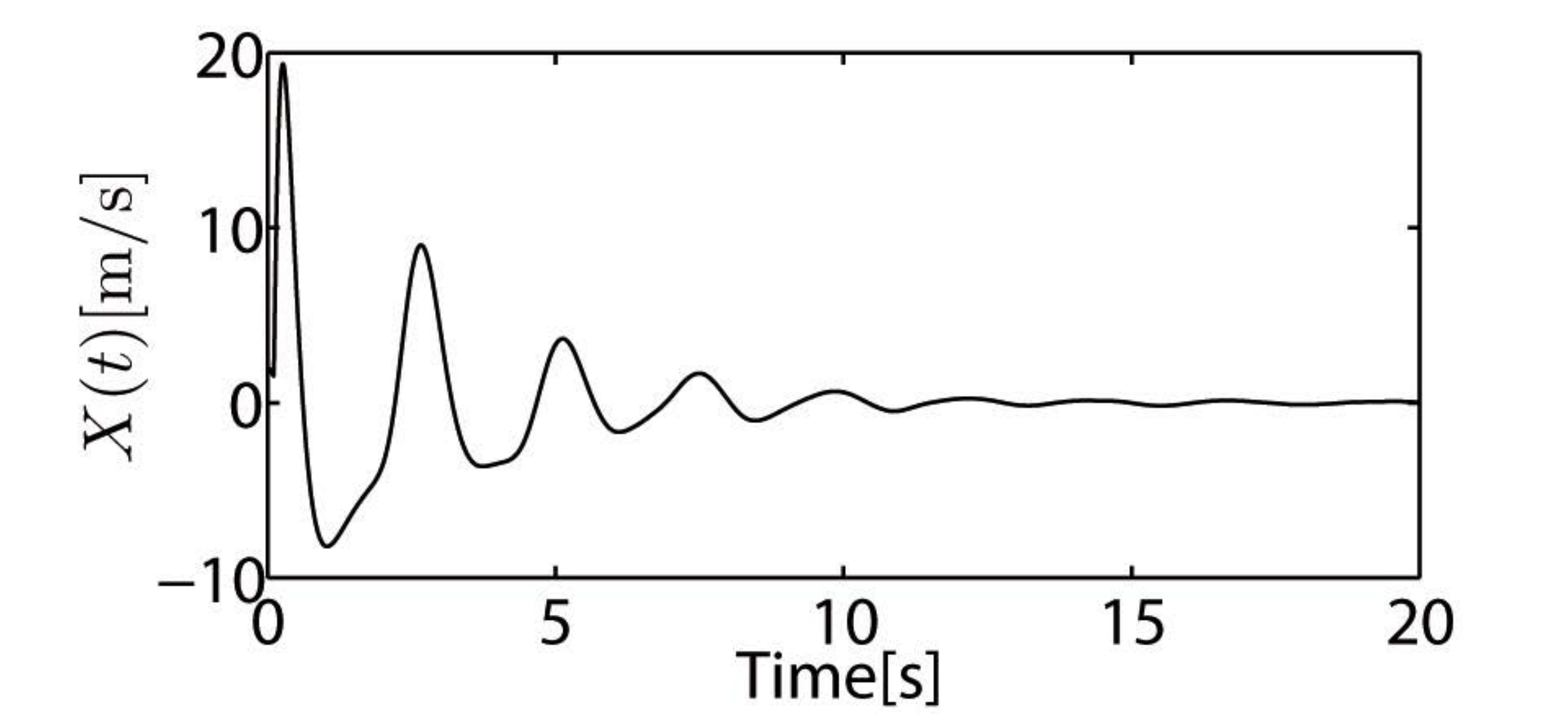}}
  \centerline{(a) $X(t)$.}
\end{minipage}
\begin{minipage}{.55\linewidth}
  \centerline{\includegraphics[width=4.5cm]{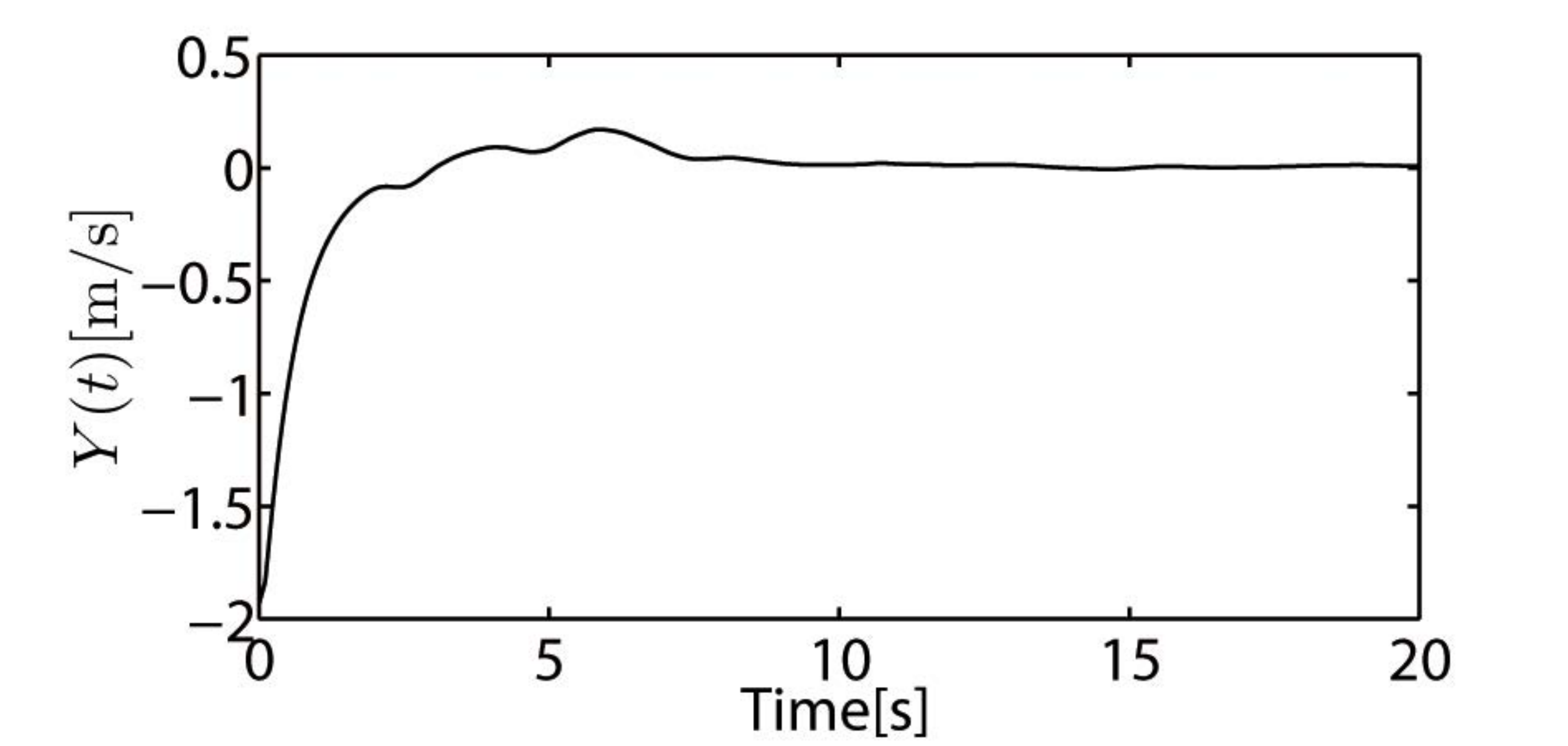}}
  \centerline{(b) $Y(t)$.}
\end{minipage}
\caption{Responses of $X(t)$, $Y(t)$ (with control).}
\label{fig:simXY}
\end{figure}
\begin{figure}
  \centerline{\includegraphics[width=7cm]{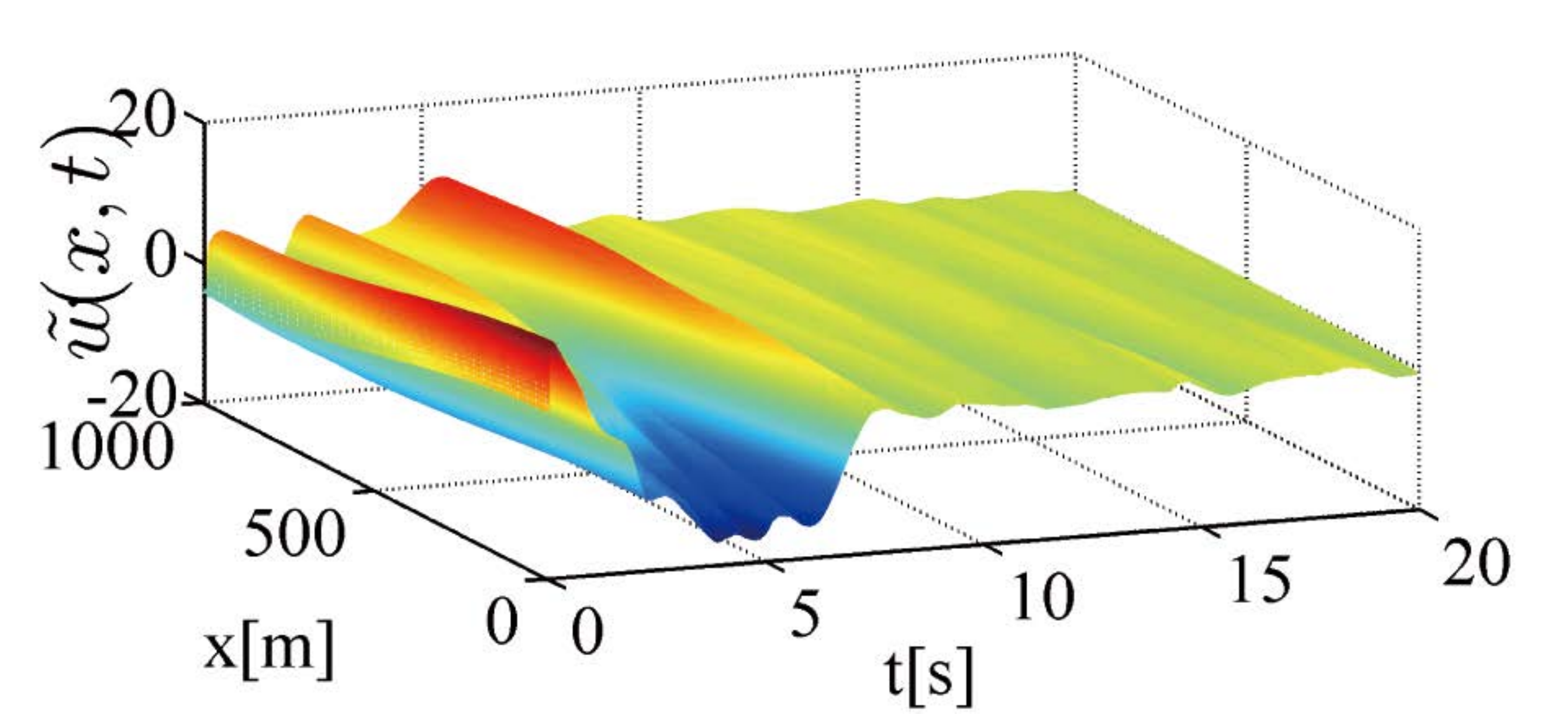}}
  \centerline{(a) $\tilde w(x,t)$.}
  \centerline{\includegraphics[width=7cm]{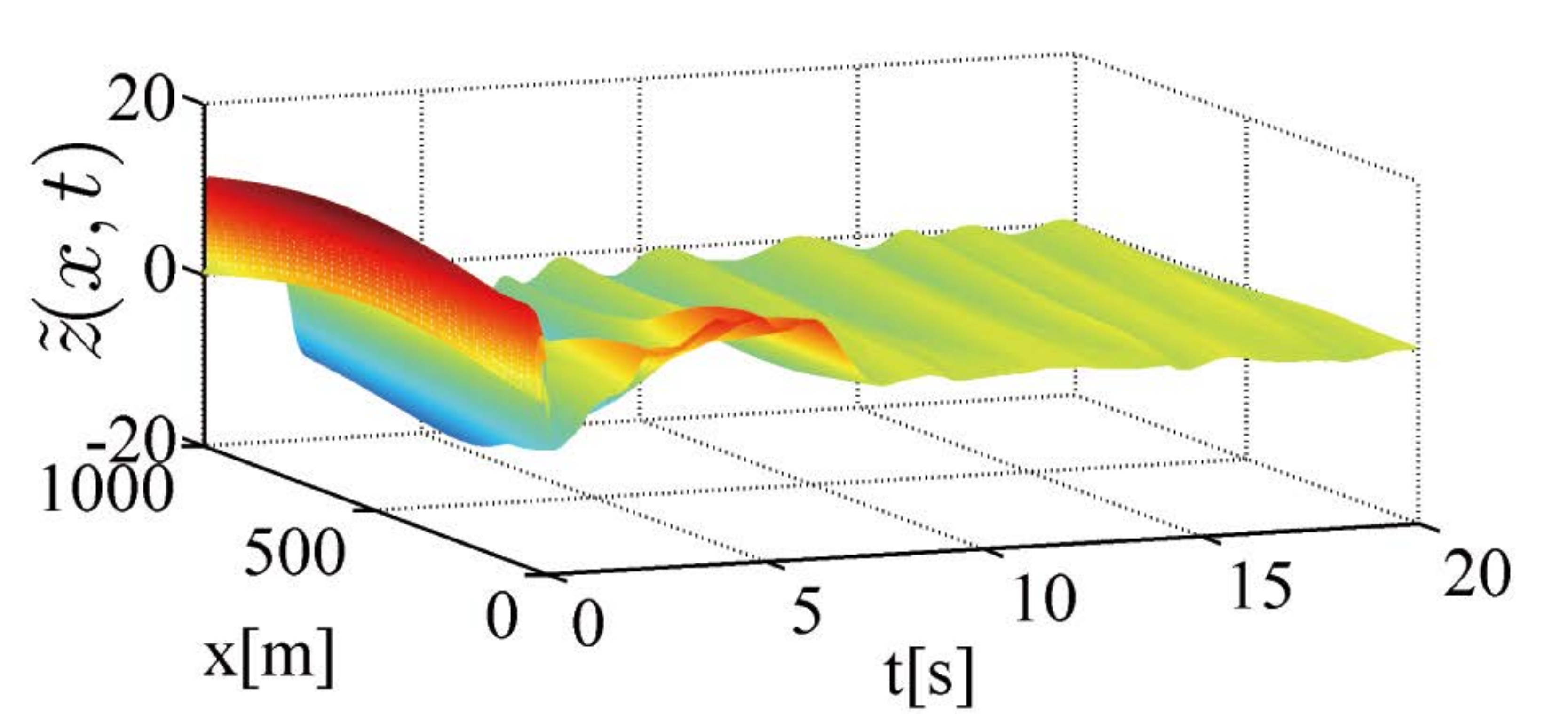}}
  \centerline{(b) $\tilde z(x,t)$.}
\caption{Observer errors $\tilde w(x,t)$, $\tilde z(x,t)$}
\label{fig:simob}
\end{figure}
According to Fig. \ref{fig:simwzopen}, we know the oscillations appear in the responses of $w(x,t)$, $z(x,t)$, which is the result of the property of the long cable and the external disturbances \eqref{eq:fxt}-\eqref{eq:fL}. From Fig. \ref{fig:simwzclosed}, we can observe that the designed control input can effectively reduce the oscillation amplitudes even though the plant is subject to the external disturbance. The moving velocity of the controlled crane and the oscillation velocity of the payload, namely $X(t)$ and $Y(t)$, are shown in Fig. \ref{fig:simXY} from which we know $X(t)$ and $Y(t)$ are convergent to zero. It also can be seen in Fig. \ref{fig:simob} that the observer errors $\tilde w(x,t)$, $\tilde z(x,t)$ converge to a small range around zero under the unknown external disturbances and the sensor delay $\tau$.
\subsubsection{Representing the obtained responses as $u,b_L$ in DSV}
\begin{figure}[!ht]
\hspace{-0.5cm}
\includegraphics[width=8.7cm]{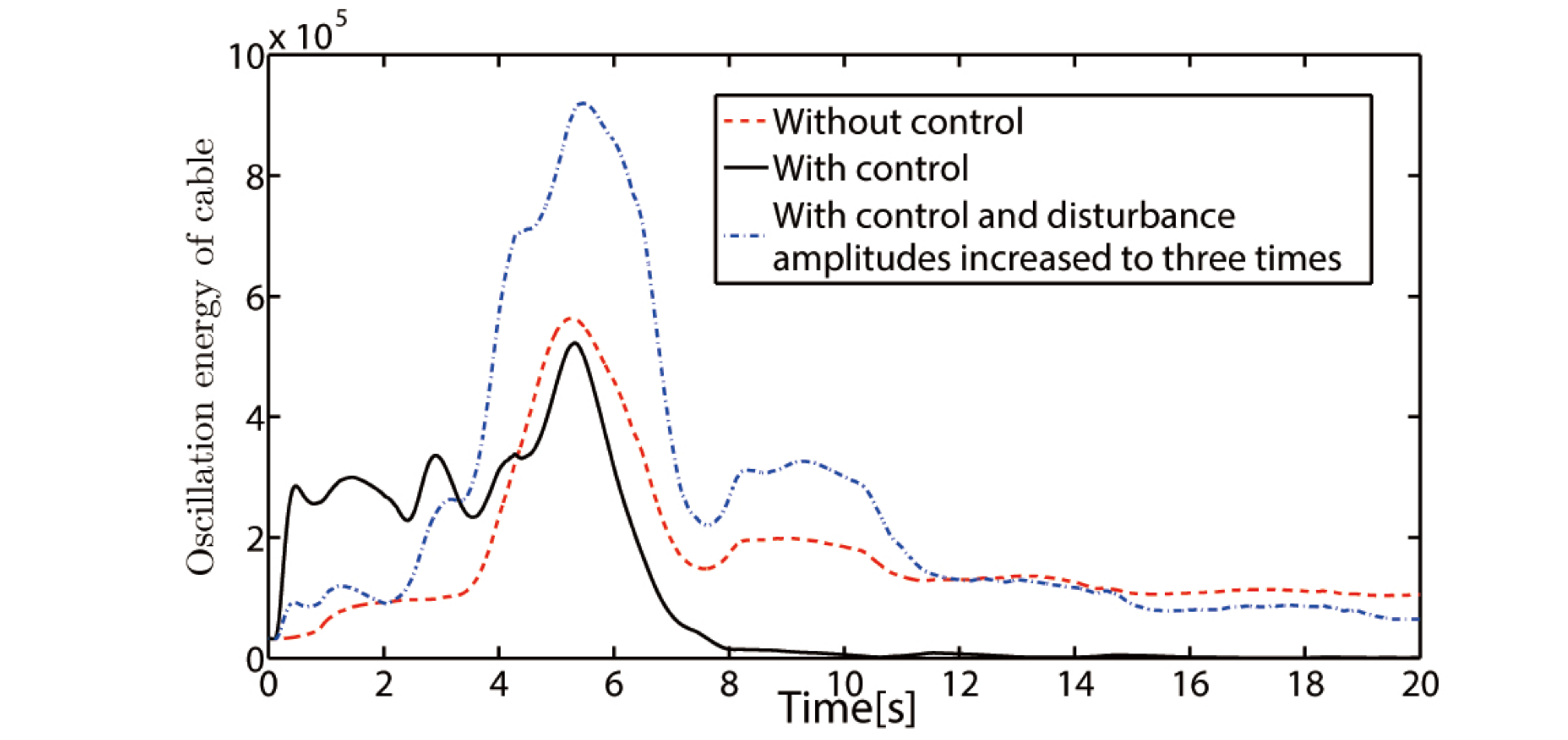}
\caption{Cable transverse oscillation energy:  $\frac{\rho}{2}\|u_t(\cdot,t)\|^2+\frac{T_0}{2}\|u_x(\cdot,t)\|^2$.}
\label{fig:simenergy}
\end{figure}

\begin{figure}[!ht]
\centering
\includegraphics[width=8.5cm]{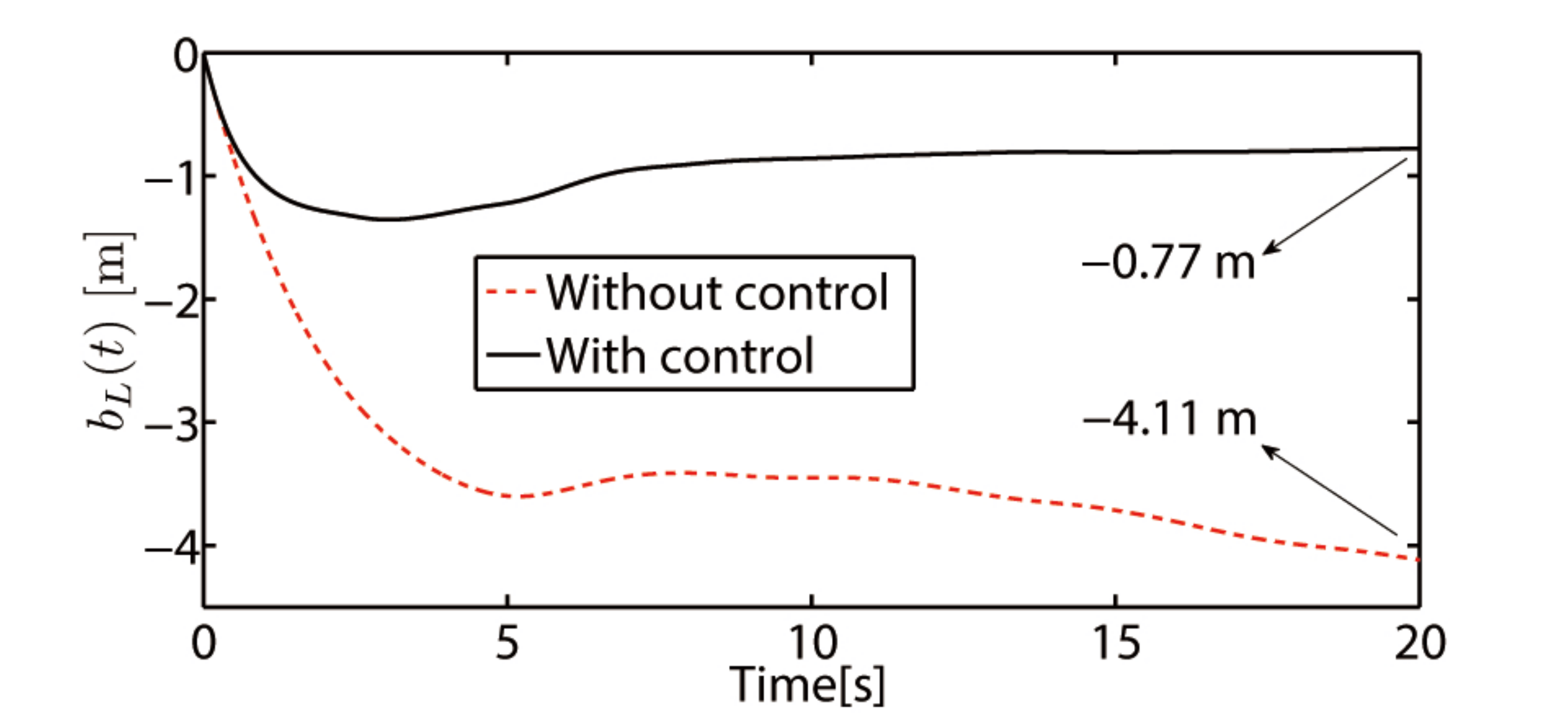}
\caption{Transverse displacement $b_L(t)$ of the payload. The end point at $t=20 s$ means the position error on the sea floor. The permissible tolerance of this typical model is 2.5 m \cite{How2011}.}
\label{fig:simbL}
\end{figure}

\begin{figure}[!ht]
\centering
\includegraphics[width=8.5cm]{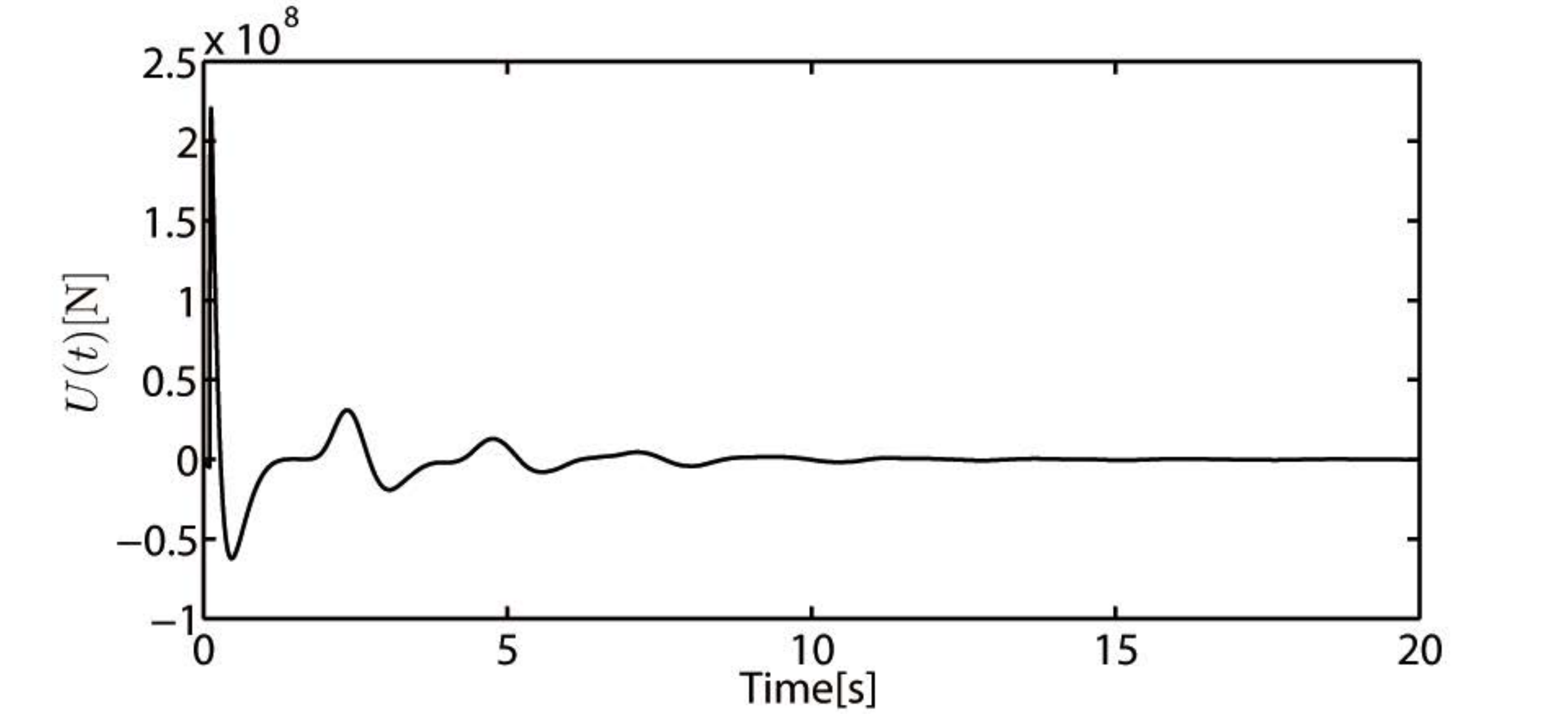}
\caption{Control force of the onboard crane.}
\label{fig:simU}
\end{figure}

The physical meaning of the responses $z,w$ in Figs. \ref{fig:simwzclosed}-\ref{fig:simXY} would be clear after representing them as the responses of the cable oscillation and position error, i.e., $u$ and $b_L$ in \eqref{eq:5.1}-\eqref{eq:5.5}. Through \eqref{eq:t11}-\eqref{eq:t12}, the cable transverse oscillation energy including oscillation kinetic energy $\frac{\rho}{2}\|u_t(\cdot,t)\|^2$ and potential energy $\frac{T_0}{2}\|u_x(\cdot,t)\|^2$ can be represented by $z(x,t),w(x,t)$ as
\begin{align}
&\frac{\rho}{2}\|u_t(\cdot,t)\|^2+\frac{T_0}{2}\|u_x(\cdot,t)\|^2\notag\\
&=\frac{\rho}{8}\|w(\cdot,t)+z(\cdot,t)\|^2+\frac{\rho}{8}\|w(\cdot,t)-z(\cdot,t)\|^2
\end{align}
where $\|u_t(\cdot,t)\|^2$ denotes  $\int_0^Lu_t(\cdot,t)^2dx$. The transverse displacement of the payload $b_L(t)$ can be obtained as
\begin{align}
b_L(t)=u(L,t)=\frac{1}{2}\int_0^{t}(z(L,\delta)+w(L,\delta))d\delta+b_L(0).
\end{align}
As shown in Fig. \ref{fig:simenergy}, the oscillation energy of the cable with the proposed control law is reduced faster and to a level below to the uncontrolled case after  $t=5.5s$, under the external disturbances \eqref{eq:fxt}-\eqref{eq:fL}. This result shows robustness of the proposed control to small disturbances. However, as we continue to increase the amplitude of the disturbance  \eqref{eq:fxt} by gradually raising $A_D$ i.e., the amplitude of the oscillating drag force, from its baseline value $400$, the blue line in Fig. \ref{fig:simenergy} shows that the controller fails to achieve effective vibration suppression once $A_D$ reaches three times the baseline value,  i.e., $A_D=1200$. From Fig. \ref{fig:simbL}, we note that position error of the payload is $-0.77$ m from the desired location on the sea floor, which satisfies the requirement of being within  permissible tolerance of 2.5 m, while the position error is $-4.11$  m in the case without control, which exceeds the tolerance. The control signal shown in Fig. \ref{fig:simU}, is bounded and convergent.
\section{Conclusion and future work}\label{sec:conclusion}
Delay-compensated control of heterodirectional
coupled hyperbolic PDEs sandwiched between two general ODEs is addressed in this paper, where the control input is applied at one ODE and the measurement is placed at another ODE with a sensor delay. Using the delayed measurement output, a full-order state observer which can compensate the sensor delay is designed to estimate the states of the overall system including the plant and sensor delay dynamics. An observer-based output-feedback controller is designed using backstepping transformations and frequency-domain design methods. The exponential stability result of the closed-loop system and the boundedness and exponential convergence of the control input are proved in this paper. The obtained theoretical results are applied to oscillation suppression of a DCV for off-shore oil drilling as a simulation case. The simulation results show the proposed control input applied at the onboard crane can reduce the oscillations of the cable and place the payload in the target area on the sea floor.

In this paper, the observer-based output-feedback controller is designed based on the model
with completely known parameters. In the future work, the model uncertainties, such as unknown plant parameters and external disturbances will be considered in such a sandwiched PDE system, and the adaptive and ADRC (Active Disturbance Rejection Control) technologies would be incorporated into the control design to solve this more practical and complex problem including parameter estimation and disturbance attenuation.
\section*{Appendices}
\subsection*{A. Matching \eqref{eq:observer1}-\eqref{eq:observer6} and \eqref{eq:targ6}-\eqref{eq:targ1}:}
\setcounter{equation}{0}
\renewcommand{\theequation}{A.\arabic{equation}}
\emph{Step 1: }Taking the time and spatial derivative of \eqref{eq:contran1b} along \eqref{eq:observer1}-\eqref{eq:observer6}, substituting the results into \eqref{eq:targ4}, we have
\begin{align}
&{\beta _t}(x,t) - {q_2}{\beta _x}(x,t) + {c_2}\beta (x,t)\notag\\
 & + \lambda (x){\Gamma _1}\tilde v(2,t)+ \int_x^1 {{J_2}(x,y){h_3}(\tilde v(2,t);y)dy}\notag\\
 &+ \int_x^1 {{K_2}(x,y){h_2}(\tilde v(2,t);y)dy}\notag\\
   &  - {h_3}(\tilde v(2,t);x)+ {q_2}{J_2}(x,1){h_4}(\tilde v(2,t))\notag\\
& = {\hat w_t}(x,t) - \int_x^1 {K_2}(x,y){\hat z_t}(y,t)dy \notag\\
&-  \int_x^1 {{J_2}(x,y){\hat w_t}(y,t)dy}  - \lambda (x)\dot {\hat Y}(t)\notag\\
& - {q_2}{\hat w_x}(x,t) + {q_2}\int_x^1{K_{2x}}(x,y)\hat z(y,t)dy\notag\\
& +  {q_2}\int_x^1 {{J_{2x}}(x,y)\hat w(y,t)dy}  + {q_2}\lambda '(x)\hat Y(t)\notag\\
& - {q_2}{K_2}(x,x)\hat z(x,t) - {q_2}{J_2}(x,x)\hat w(x,t)\notag\\
& + {c_2}\hat w(x,t) - {c_2}\int_x^1 {K_2}(x,y)\hat z(y,t)dy \notag\\
&-  {c_2}\int_x^1 {{J_2}(x,y)\hat w(y,t)dy}  - {c_2}\lambda (x)\hat Y(t)\notag\\
 & + \lambda (x){\Gamma _1}\tilde v(2,t)+ \int_x^1 {{J_2}(x,y){h_3}(\tilde v(2,t);y)dy}\notag\\
 &+ \int_x^1 {{K_2}(x,y){h_2}(\tilde v(2,t);y)dy}\notag\\
   &  - {h_3}(\tilde v(2,t);x)+ {q_2}{J_2}(x,1){h_4}(\tilde v(2,t))\notag\\
& =  - {c_2}\hat z(x,t)  + {h_3}(\tilde v(2,t);x) + {q_1}\int_x^1 {K_2}(x,y){\hat z_x}(y,t)dy\notag\\
 &+ \int_x^1 {{c_1}{K_2}(x,y)[{\hat z}(x,t) + {\hat w}(y,t)]dy}\notag \\
 &  - \int_x^1 {{K_2}(x,y){h_2}(\tilde v(2,t);y)dy} \notag\\
 &- {q_2}\int_x^1 {{J_2}(x,y){\hat w_x}(y,t)dy}\notag \\
 & + \int_x^1 {{c_2}{J_2}(x,y)[\hat z(y,t) + \hat w(y,t)]dy}\notag \\
& - \int_x^1 {{J_2}(x,y){h_3}(\tilde v(2,t);y)dy} \notag\\
 &+ {q_2}\int_x^1 {{K_{2x}}(x,y)\hat z(y,t)dy + } {q_2}\int_x^1 {{J_{2x}}(x,y)\hat w(y,t)dy}\notag \\
 &- {q_2}{K_2}(x,x)\hat z(x,t) - {q_2}{J_2}(x,x)\hat w(x,t)\notag\\
 &- \lambda (x){A_1}\hat Y(t) - \lambda (x){B_1}\hat z(1,t) \notag\\
 &- \lambda (x){\Gamma _1}\tilde v(2,t) + {q_2}\lambda '(x)\hat Y(t)\notag\\
 & - {c_2}\int_x^1 {K_2}(x,y)\hat z(y,t)dy \notag\\
 &-  {c_2}\int_x^1 {{J_2}(x,y)\hat w(y,t)dy}  - {c_2}\lambda (x)\hat Y(t)\notag\\
 & + \lambda (x){\Gamma _1}\tilde v(2,t)+ \int_x^1 {{J_2}(x,y){h_3}(\tilde v(2,t);y)dy}\notag\\
 &+ \int_x^1 {{K_2}(x,y){h_2}(\tilde v(2,t);y)dy}\notag\\
   &  - {h_3}(\tilde v(2,t);x)+ {q_2}{J_2}(x,1){h_4}(\tilde v(2,t))\notag\\
& = [ - {c_2} - ({q_1} + {q_2}){K_2}(x,x)]\hat z(x,t)\notag\\
 &+ \int_x^1 \left[{c_1}{K_2}(x,y) + {q_2}{J_{2x}}(x,y) + {q_2}{J_{2y}}(x,y)\right]\hat w(y,t)dy \notag\\
 &+ \int_x^1 \big[{c_2}{J_2}(x,y) + {c_1}{K_2}(x,y) - {c_2}{K_2}(x,y)\notag \\
 &+ {q_2}{K_{2x}}(x,y) - {q_1}{K_{2y}}(x,y)\big]{\hat z}(y,t)dy \notag\\
 &+\left[{q_1}{K_2}(x,1) - {q_2}{J_2}(x,1)q - \lambda (x){B_1}\right]\hat z(1,t)\notag\\
 &+ \left[{q_2}\lambda '(x) - \lambda (x){A_1} - {c_2}\lambda (x) - {q_2}{J_2}(x,1){C_1}\right]\hat Y(t).\label{eq:s1}
\end{align}
For $\eqref{eq:s1}$ to hold, conditions \eqref{eq:Kerc2}-\eqref{eq:q2lam} should be satisfied.

\emph{Step 2:} Taking the time and spatial derivative of \eqref{eq:contran1a} along \eqref{eq:observer1}-\eqref{eq:observer6}, substituting the results into \eqref{eq:targ3}, we have
\begin{align}
&{\alpha _t}(x,t) + {q_1}{\alpha _x}(x,t) + {c_1}\alpha (x,t)\notag\\
&+ \gamma (x){\Gamma _1}\tilde v(2,t) + \int_x^1 {{J_3}(x,y){h_3}(\tilde v(2,t);y)dy}\notag \\
 &+ \int_x^1 {{K_3}(x,y){h_2}(\tilde v(2,t);y)dy}  - h_2(\tilde v(2,t);x)\notag\\
 & + {q_2}{J_3}(x,1){h_4}(\tilde v(2,t))\notag\\
& = {\hat z_t}(x,t) - \int_x^1 {K_3}(x,y){\hat z_t}(y,t)dy -  \int_x^1 {{J_3}(x,y){\hat w_t}(y,t)dy}\notag\\
&   - \gamma (x)\dot {\hat Y}(t)+ {q_1}{\hat z_x}(x,t) - {q_1}\int_x^1 {K_{3x}}(x,y)\hat z(y,t)dy \notag\\
&-  {q_1}\int_x^1 {{J_{3x}}(x,y)\hat w(y,t)dy}  - {q_1}\gamma '(x)\hat Y(t)\notag\\
& + {q_1}{K_3}(x,x)\hat z(x,t) + {q_1}{J_3}(x,x)\hat w(x,t)\notag\\
& + {c_1}\hat z(x,t) - {c_1}\int_x^1 {K_3}(x,y)\hat z(y,t)dy\notag\\
& -  {c_1}\int_x^1 {{J_3}(x,y)\hat w(y,t)dy - {c_1}} \gamma (x)\hat Y(t)\notag\\
&+ \gamma (x){\Gamma _1}\tilde v(2,t) + \int_x^1 {{J_3}(x,y){h_3}(\tilde v(2,t);y)dy}\notag \\
 &+ \int_x^1 {{K_3}(x,y){h_2}(\tilde v(2,t);y)dy}  - h_2(\tilde v(2,t);x)\notag\\
 & + {q_2}{J_3}(x,1){h_4}(\tilde v(2,t))\notag\\
& =  - {c_1} \hat w(x,t) + h_2(\tilde v(2,t);x) + {q_1}\int_x^1 {K_3}(x,y){\hat z_x}(y,t)dy \notag\\
&+ \int_x^1 {{c_1}{K_3}(x,y)[\hat z(y,t) + \hat w(y,t)]dy} \notag\\
&  - \int_x^1 {{K_3}(x,y){h_2}(\tilde v(2,t);y)dy} \notag\\
& - {q_2}\int_x^1 {{J_3}(x,y){\hat w_x}(y,t)dy} \notag\\
& + \int_x^1 {{c_2}{J_3}(x,y)[\hat z(y,t) + \hat w(x,t)]dy} \notag\\
& - \int_x^1 {{J_3}(x,y){h_3}(\tilde v(2,t);y)dy} \notag\\
& - {q_1}\int_x^1 {K_{3x}}(x,y)\hat z(y,t)dy -  {q_1}\int_x^1 {{J_{3x}}(x,y)\hat w(y,t)dy} \notag\\
& + {q_1}{K_3}(x,x)\hat z(x,t) + {q_1}{J_3}(x,x)\hat w(x,t)\notag\\
& - \gamma (x){A_1}\hat Y(t) - \gamma (x){B_1}\hat z(1,t) - \gamma (x){\Gamma _1}\tilde v(2,t) - {q_1}\gamma '(x)\hat Y(t)\notag\\
& - {c_1}\int_x^1 {K_3}(x,y)\hat z(y,t)dy\notag\\
& -  {c_1}\int_x^1 {{J_3}(x,y)\hat w(y,t)dy - {c_1}} \gamma (x)\hat Y(t)\notag\\
&+ \gamma (x){\Gamma _1}\tilde v(2,t) + \int_x^1 {{J_3}(x,y){h_3}(\tilde v(2,t);y)dy}\notag \\
 &+ \int_x^1 {{K_3}(x,y){h_2}(\tilde v(2,t);y)dy}  - h_2(\tilde v(2,t);x)\notag\\
 & + {q_2}{J_3}(x,1){h_4}(\tilde v(2,t))\notag\\
& = \left[({q_2} + {q_1}){J_3}(x,x) - {c_1}\right]\hat w(x,t)\notag\\
& + \int_x^1 \big[{c_1}{K_3}(x,y) - {q_1}{J_{3x}}(x,y) \notag\\
&+ {q_2}{J_{3y}}(x,y) + ({c_2} - {c_1}){J_3}(x,y)\big]\hat w(y,t)dy \notag\\
& + \int_x^1 {\left[{c_2}{J_3}(x,y) - {q_1}{K_{3x}}(x,y) - {q_1}{K_{3y}}(x,y)\right]{\hat z}(y,t)dy} \notag\\
&+\left[{q_1}{K_3}(x,1) - {q_2}{J_3}(x,1)q - \gamma (x){B_1}\right]\hat z(1,t)\notag\\
& + \left[ - {q_1}\gamma'(x) - \gamma (x){A_1} - {c_1}\gamma (x) - {q_2}{J_3}(x,1){C_1}\right]\hat Y(t).\label{eq:s2}
\end{align}
For \eqref{eq:s2} to hold, conditions \eqref{eq:Kerc1}-\eqref{eq:c2L3}, \eqref{eq:Kerc1a} should be satisfied.

\emph{Step 3: }Inserting \eqref{eq:contran1a}-\eqref{eq:contran1b} into \eqref{eq:targ1},\eqref{eq:targ2} respectively, applying \eqref{eq:observer5}-\eqref{eq:observer6}, we obtain
\begin{align}
&\dot {\hat Y}(t) - {\bar A_1}\hat Y(t) - {B_1}\alpha (1,t) - {\Gamma _1}\tilde v(2,t)\notag\\
=&\dot {\hat Y}(t) - {A_1}\hat Y(t)+B_1F_1\hat Y(t) - {B_1}\hat z (1,t)+B_1\gamma(1)\hat Y(t)\notag\\
& - {\Gamma _1}\tilde v(2,t)\notag\\
=& B_1(F_1+\gamma(1))\hat Y(t) =0,\label{eq:3a}\\
&\beta (1,t) -q\alpha (1,t) - {h_4}(\tilde v(2,t))\notag\\
=&\hat w(1,t) - q\hat z(1,t) +(q\gamma(1)-\lambda(1))\hat Y(t)- {h_4}(\tilde v(2,t))\notag\\
=&(q\gamma(1)-\lambda(1)+C_1)\hat Y(t)=0.\label{eq:3b}
\end{align}
For \eqref{eq:3a}-\eqref{eq:3b} to hold, conditions \eqref{eq:gammaF1} and \eqref{eq:Kerc2a} should be satisfied.

\emph{Step 4: }Inserting \eqref{eq:contran1a}-\eqref{eq:contran1b} into \eqref{eq:targ5} and \eqref{eq:targ6}, applying \eqref{eq:observer2}, \eqref{eq:observer1} we have
\begin{align}
&\alpha (0,t)- p\beta (0,t) - \int_0^1 {{\bar K}_1}(x)\alpha (x,t)dx -  {\bar K}_3\hat Y(t)\notag\\
 &- \int_0^1 {{\bar K}_2}(x)\beta (x,t)dx - {C_0}\hat X(t)\notag\\
 =&\int_0^1 \bigg[ p{K_2}(0,x) - {K_3}(0,x) + \int_0^x {{\bar K}_1}(y){K_3}(y,x)dy\notag\\
  &+ \int_0^x {{{\bar K}_2}(y){K_2}(y,x)dy - {{\bar K}_1}(x)}  \bigg] z(x,t)dx \notag\\
 &- \int_0^1 \bigg[ p{J_2}(0,x) - {J_3}(0,x) - \int_0^x {{\bar K}_1}(y){J_3}(y,x)dy\notag\\
  &- \int_0^x {{\bar K}_2}(y){J_2}(y,x)dy + {{\bar K}_2}(x)  \bigg] w(x,t)dx \notag\\
 &+ \bigg[\int_0^1 {{{\bar K}_2}(x)\lambda (x)dx}  + \int_0^1 {{\bar K}_1}(x)\gamma (x)dx\notag\\
  &+ p\lambda (0) - \gamma (0) -  {{\bar K}_3}\bigg]Y(t) = 0,\label{eq:alpah0}\\
&\dot {\hat X}(t) - {A_0}\hat X(t) - {E_0}\beta (0,t) - \int_0^1 {{\bar K}_4}(x)\alpha (x,t)dx \notag\\
&- \int_0^1 {{{\bar K}_5}(x)\beta (x,t)d} x - {{\bar K}_6}\hat Y(t)- {B_0}U(t) - {h_1}(\tilde v(2,t))\notag\\
 =&  - \int_0^1 \bigg[ \int_0^x {{\bar K}_4}(y){K_3}(y,x)dy + \int_0^x {{\bar K}_5}(y){K_2}(y,x)dy\notag\\
 &- {E_0}{K_2}(0,x) - {{\bar K}_4}(x)  \bigg] z(x,t)dx \notag\\
 &- \int_0^1 \bigg[ \int_0^x {{\bar K}_4}(y){J_3}(y,x)dy + \int_0^x {{\bar K}_5}(y){J_2}(y,x)dy\notag\\
  &+ {E_0}{J_2}(0,x) - {{\bar K}_5}(x)  \bigg] w(x,t)dx \notag\\
& - \bigg[\int_0^1 {{{\bar K}_5}(x)\lambda (x)dx}  + \int_0^1 {{\bar K}_4}(x)\gamma (x)dx \notag\\
&+ {E_0}\lambda (0) -  {{\bar K}_6}\bigg]Y(t)=0.\label{eq:dX}
\end{align}
For \eqref{eq:alpah0}-\eqref{eq:dX} to hold, conditions \eqref{eq:barK1}-\eqref{eq:barK6} are obtained.
\subsection*{B. Calculations of \eqref{eq:wtwx} and \eqref{eq:ztzx}}\label{app:B}
\setcounter{equation}{0}
\renewcommand{\theequation}{B.\arabic{equation}}
Details of calculating  \eqref{eq:wtwx}:
\begin{align}
&{{\tilde w}_t}(x,t) - {q_2}{{\tilde w}_x}(x,t)+ {c_2}\tilde z(x,t)+ {c_2}\tilde w(x,t)\notag\\
& + {h_3(\tilde v(2,t);x)}\notag\\
 &= {{\tilde \beta }_t}(x,t) - \int_x^1 {\psi (x,y){{\tilde \alpha }_t}(y,t)dy} \notag\\
 &- {q_2}{{\tilde \beta }_x}(x,t) + {q_2}\int_x^1 {{{\psi }_x}(x,y)\tilde \alpha (y,t)dy} - {q_2}\psi (x,x)\tilde \alpha (x,t)\notag\\
 &+ {c_2}\tilde \alpha (x,t) - {c_2}\int_x^1 {\phi (x,y)\tilde \alpha (y,t)dy}  + {h_3(\tilde v(2,t);x)}\notag\\
 & + {c_2}\tilde \beta (x,t) - {c_2}\int_x^1 {\psi (x,y)\tilde \alpha (y,t)dy}\notag\\
 &= \int_x^1 {\bar N(x,y)\tilde \beta (y,t)dy}+{c_1}\int_x^1 {\psi (x,y)\tilde \alpha (y,t)dy} \notag\\
 &+ {q_1}\int_x^1 {\psi (x,y){{\tilde \alpha }_x}(y,t)dy + \int_x^1 {{c_1}\psi (x,y)\tilde \beta (y,t)dy} }\notag \\
 &- \int_x^1 {\int_x^y {\psi (x,\delta)\bar M(\delta,y)} d\delta\tilde \beta (y,t)dy} \notag\\
  &+ {q_2}\int_x^1 {{{\psi }_x}(x,y)\tilde \alpha (y,t)dy} \notag\\
 &- {q_2}\psi (x,x)\tilde \alpha (x,t)- {c_2}\int_x^1 {\psi (x,y)\tilde \alpha (y,t)dy}\notag\\
 &+ {c_2}\tilde \alpha (x,t) - {c_2}\int_x^1 {\phi (x,y)\tilde \alpha (y,t)dy}  + {h_3(\tilde v(2,t);x)}\notag\\
 &= [c_2 - ({q_1} + {q_2})\psi (x,x)]\tilde \alpha (x,t)\notag\\
 &+ \int_x^1 [ - {q_1}{{\psi }_y}(x,y) + {q_2}{{\psi }_x}(x,y) - {c_2}\phi (x,y)\notag\\
 &-(c_2-{c_1})\psi (x,y)]\tilde \alpha (x,t)dy \notag\\
 &+ \int_x^1 [{c_1}\psi (x,y) + \bar N(x,y) \notag\\
 &- \int_x^y {\psi (x,\delta)\bar M(\delta,y)} d\delta]\tilde \beta (y,t)dy \notag\\
 &+{q_1}\psi (x,1)\tilde\alpha (1,t) + {h_3(\tilde v(2,t);x)}\notag\\
 &= {q_1}\psi (x,1)\tilde\alpha (1,t) + {h_3(\tilde v(2,t);x)}=0.\label{eq:appB1}
\end{align}
Details of calculating \eqref{eq:ztzx}:
\begin{align}
&{{\tilde z}_t}(x,t) + {q_1}{{\tilde z}_x}(x,t)  + {c_1}\tilde w(x,t)+ {c_1}\tilde z(x,t)+ h_2(\tilde v(2,t);x)\notag\\
 &= {{\tilde \alpha }_t}(x,t) - \int_x^1 {\phi (x,y){{\tilde \alpha }_t}(y,t)dy}\notag \\
 &+ {q_1}{{\tilde \alpha }_x}(x,t) - {q_1}\int_x^1 {{{\phi }_x}(x,y)\tilde \alpha (y,t)dy}  + {q_1}\phi (x,x)\tilde \alpha (x,t)\notag\\
 &+ {c_1}\tilde \beta (x,t) - {c_1}\int_x^1 {\psi (x,y)\tilde \alpha (y,t)dy + h_2(\tilde v(2,t);x)}\notag \\
 &+ {c_1}\tilde \alpha (x,t) - {c_1}\int_x^1 {\phi (x,y)\tilde \alpha (y,t)dy}\notag\\
 &=  -{c_1}\tilde \beta (x,t) + {q_1}\int_x^1 {\phi (x,y){{\tilde \alpha }_x}(y,t)dy} \notag\\
 &+ \int_x^1 {[ - {q_1}{{\phi }_x}(x,y) - {c_1}\psi (x,y)]\tilde \alpha (y,t)dy} \notag\\
 &+ \int_x^1 {[ - \int_x^y {\phi (x,\delta)\bar M(y,\delta)} d\delta + \phi (x,y){c_1}]\tilde \beta (y,t)dy} \notag\\
 &+ {q_1}\phi (x,x)\tilde \alpha (x,t) +{c_1}\tilde \beta (x,t)\notag\\
 & +\int_x^1 {\bar M(x,y)\tilde \beta (y,t)dy}  + h_2(\tilde v(2,t);x)\notag\\
 &= \int_x^{1} {[ - {q_1}{{\phi }_x}(x,y) - {c_1}\psi (x,y) - {q_1}{{\phi }_y}(x,y)]\tilde \alpha (y,t)dy} \notag\\
 &+ \int_x^{1} {[\bar M(x,y) - \int_x^y {\phi (x,\delta)\bar M(\delta,y)} d\delta + {c_1}\phi (x,y)]\tilde \beta (y,t)dy}\notag \\
 &+{q_1}\phi (x,1)\tilde \alpha (1,t)+ h_2(\tilde v(2,t);x)\notag\\
 &={q_1}\phi (x,1)\tilde \alpha (1,t)+ h_2(\tilde v(2,t);x)=0.\label{eq:appB2}
\end{align}

\end{document}